\documentclass[10pt]{article}
\usepackage[margin=1in
]{geometry}

\usepackage{titlesec}
\titlespacing*{\paragraph}{0pt}{1ex plus .5ex minus .2ex}{1em}

\usepackage[T1]{fontenc}
\usepackage{textcomp}

\usepackage{libertinus}

\usepackage{url}
\usepackage{diagbox}

\usepackage{natbib}

\usepackage{amsthm,amsfonts,amsmath,amssymb,epsfig,color,float,graphicx,verbatim,bm,bbm}
\usepackage{enumerate}
\usepackage{enumitem}
\usepackage{wrapfig}
\usepackage{subcaption}
\usepackage[export]{adjustbox}
\usepackage{nicefrac}
\usepackage{hhline}
\usepackage{multicol}
\usepackage{multirow}
\usepackage[dvipsnames]{xcolor}

\usepackage[ruled, vlined, linesnumbered]{algorithm2e}

\usepackage{hyperref}
\hypersetup{
	colorlinks   = true, 
	urlcolor     = blue!75!black, 
	linkcolor    = blue!75!black, 
	citecolor   = blue!75!black 
}

\usepackage{cleveref}
\crefname{fact}{fact}{facts}
\Crefname{fact}{Fact}{Facts}
\crefname{claim}{claim}{claims}
\Crefname{claim}{Claim}{Claims}

\newtheorem{theorem}{Theorem}[section]

\newtheorem*{question*}{Question}

\newtheorem*{bigquestion*}{Big Question}

\newtheorem{lemma}[theorem]{Lemma}
\newtheorem{claim}[theorem]{Claim}

\newtheorem{definition}[theorem]{Definition}
\newtheorem*{definition*}{Definition}
\newtheorem{remark}[theorem]{Remark}

\newtheorem{fact}[theorem]{Fact}
\newtheorem*{keyquestion*}{Key Question}

\newcommand{\barbv}{\bar{\mathbf{v}}}
\newcommand{\barbu}{\bar{\mathbf{u}}}

\newcommand{\nnz}{\operatorname{nnz}}

\newcommand{\val}{\operatorname{val}}
\newcommand{\width}{\operatorname{width}}

\newcommand{\reals}{\mathbb{R}}

\newcommand{\argmin}[1]{{\operatorname{arg\,min}}_{#1}}

\newcommand{\Deltatr}{\Delta^{\rm tr}}
\newcommand{\epstr}{\epsilon_{\rm tr}}
\newcommand{\diag}{\operatorname{diag}}
\newcommand{\dd}{\mathrm{d}}

\newcommand{\bzero}{\mathbf{0}}
\newcommand{\bone}{\mathbf{1}}

\newcommand{\ba}{\mathbf{a}}
\newcommand{\bp}{\mathbf{p}}
\newcommand{\be}{\mathbf{e}}

\newcommand{\bx}{\mathbf{x}}
\newcommand{\bk}{\mathbf{k}}
\newcommand{\x}{\mathbf{x}}
\newcommand{\bw}{\mathbf{w}}

\newcommand{\bg}{\mathbf{g}}

\newcommand{\bb}{\mathbf{b}}
\newcommand{\bu}{\mathbf{u}}

\newcommand{\bv}{\mathbf{v}}

\newcommand{\bz}{\mathbf{z}}

\newcommand{\bc}{\mathbf{c}}
\newcommand{\bd}{\mathbf{d}}

\newcommand{\by}{\mathbf{y}}

\newcommand{\bq}{\mathbf{q}}

\newcommand{\bpi}{\boldsymbol{\Pi}}
\newcommand{\bgamma}{\boldsymbol{\Gamma}}

\newcommand{\mA}{\mathbf{A}}
\newcommand{\mAh}{\widehat{\mathbf{A}}}
\newcommand{\Ah}{\widehat{A}}
\newcommand{\mB}{\mathbf{B}}
\newcommand{\mF}{\mathbf{F}}

\newcommand{\bxh}{\widehat{\bx}}
\newcommand{\byh}{\widehat{\by}}
\newcommand{\bqh}{\widehat{\bq}}
\newcommand{\bzh}{\widehat{\bz}}

\newcommand{\Lcal}{\mathcal{L}}

\newcommand{\Ccal}{\mathcal{C}}

\newcommand{\Scal}{\mathcal{S}}

\newcommand{\Xcal}{\mathcal{X}}

\newcommand{\Ycal}{\mathcal{Y}}
\newcommand{\Zcal}{\mathcal{Z}}

\newcommand{\norm}[1]{\left\|#1\right\|}
\newcommand{\inner}[1]{\langle#1\rangle}

\newcommand{\out}{\mathrm{out}}

\newcommand{\relint}{\mathrm{rel\,int}}

\newcommand{\R}{\mathbb{R}}
\newcommand{\Z}{\mathbb{Z}}

\newcommand{\eps}{\epsilon}

\title{Nearly-Tight Oracle Complexity of \\Packing and Covering Linear Programs}
\author{Jelena Diakonikolas\\
UW-Madison\\
\texttt{jelena@cs.wisc.edu}
\and
Guy Kornowski\\
EPFL \& Institute for Advanced Study\\
\texttt{gkornowski@ias.edu}
}

\date{}

\begin{document}

\maketitle

\begin{abstract}
We settle the deterministic oracle complexity of packing and covering linear programs, as measured by the number of oracle queries to the constraint matrix modeled by left and right matrix-vector products, up to a logarithmic factor. We do so by providing nearly-matching upper and lower bounds. Our algorithm breaks the $1/\epsilon^2$ barrier of width-independent, $(1+\epsilon)$-multiplicative error parallel solvers from prior work, improving the scaling with the error $\epsilon >0$ to $1/\epsilon,$ with total number of arithmetic operations scaling near-linearly with $1/\epsilon$ and the number of nonzero entries of the constraint matrix. 
 \end{abstract}

\begin{center}
\large{
\emph{Dedicated to the memory of Michael B.\ Cohen.} 
}
\end{center}

\section{Introduction}

Packing and covering linear programs (LPs) are fundamental problems in theoretical computer science and operations research. They arise in network flow problems and the approximation algorithms built on them \cite{ShahrokhiMatula1990,LeightonRao1999,GargKonemann2007,Fleischer2000}; as fractional relaxations of computationally hard problems such as set cover, vertex cover, bin packing, and the metric traveling salesman problem \cite{Lovasz1975ratio,Chvatal1979greedy,Hochbaum1982,KarmarkarKarp1982,HobergRothvoss2017,HeldKarp1970,GoemansBertsimas1993,Goemans1995worstcase,KarlinKleinOveisGharan2022,TraubVygen2024}; in online algorithms for problems such as ski rental, paging, and ad allocation \cite{BuchbinderNaor2009FnT,BuchbinderNaor2009MOR,BuchbinderJainNaor2007,MehtaSaberiVazirani2007adwords,DevanurHayes2009adwords,BamasMaggioriSvensson2020}, and in online resource allocation more broadly \cite{AgrawalWangYe2014,BalseiroLuMirrokni2023}; as zero-sum games with nonnegative payoffs, which underlie boosting \cite{FreundSchapire1999,AHK12}; in computing densest subgraphs \cite{Charikar2000densest,BoobEtAl2020Flowless,ChekuriQuanrudTorres2022}; and as relaxations of packing and covering integer programs \cite{Raghavan1988packingIP,Srinivasan1999packingcovering}, among many other applications \cite{Vazirani2001,WilliamsonShmoys2011}. 

As linear programs, packing and covering LPs can be solved with standard LP solvers, in time polynomial in the size of the program (measured by the number of variables and the number of constraints) and logarithmic in the inverse error $1/\epsilon,$ $\epsilon > 0.$ Yet, in large-scale problem settings or when high-accuracy solvers are not needed, polynomial dependence on $1/\epsilon$ can be tolerated for solvers that are overall faster by having a near-linear scaling with the LP size, especially if they are amenable to parallelization with depth that is poly-logarithmic in the LP size and polynomial in the inverse error. For this reason, there is substantial literature on developing such parallelizable, near-linear-time (in the input size) solvers, developed over the past 30+ years; see, e.g., \cite{LN93,You01,KuhnMoscibrodaWattenhofer2006,awerbuch2008stateless,AO15,MRWZ16,DO17}. 

An important property of packing and covering LPs is that they are \emph{scale-invariant}: if all entries in the entrywise-nonnegative constraint matrix $\mA \in \R^{m \times n}$ associated with their formulation are scaled by the same positive constant, then their optimal solutions simply scale by the (inverse of the) same constant. This property suggests that multiplicative error is more appropriate for solving these positive LPs, as it is likewise scale-invariant. However, from an optimization perspective, iterative algorithms are typically more compatible with additive error, and simple rescaling to translate additive error to multiplicative error leads to the runtime scaling polynomially with the matrix \emph{width}, defined as the maximum ratio of the matrix's positive entries. Since matrix width can be arbitrarily large,
the design of efficient, near-linear-time algorithms for packing and covering LPs that are \emph{width-independent}---namely, whose runtime is independent of the matrix width---has been central to the literature on packing and covering LPs \cite{You01,AO15,MRWZ16,DO17}. 

Prior to this work, state-of-the-art width-independent near-linear time solvers, due to \cite{AO19,WRM16}, used $\widetilde{O}\big(\frac{\nnz(\mA)}{\epsilon}\big)$ arithmetic operations (a.k.a.\ work), where tilde notation hides polylogarithmic factors in its argument and $\nnz(\mA)$ counts the number of nonzero entries of $\mA.$ These solvers, however, are crucially sequential, as they rely on updating one randomly selected coordinate per iteration, so their parallel depth is linear in the number of variables. Among parallel solvers with poly($\log \nnz(\mA), 1/\epsilon$) depth, $\widetilde{O}\big(\frac{\nnz(\mA)}{\epsilon^2}\big)$ work has remained state-of-the-art for the past 10 years \cite{AO15,MRWZ16,DO17}.
At the same time, there have been nearly no informative lower bounds for this problem. An algorithm-specific lower bound by \cite{KleinYoung1999} applies only for so-called Dantzig-Wolfe-type methods, and is higher than the complexity of the discussed methods, as they are not of this type. On the other hand, the algorithm-independent lower bound in \cite{kyng2020packing} is conditional on hardness of sparse linear system solves, and only rules out $\widetilde{O}\big(\frac{\nnz(\mA)}{\epsilon^{0.165}}\big)$ total work. 
This has left open two core questions: (1) whether the $1/\epsilon$ gap between sequential and parallel algorithms is intrinsic, and (2) whether $\widetilde{O}(\frac{\nnz(\mA)}{\epsilon})$ is a barrier beyond which no further improvement on the $\epsilon$ dependence is possible for a broad class of algorithms. We address both questions in this work. 

\subsection{Problem Setup and Main Results}\label{sec:prob-and-main-results}

We study packing and covering LPs in their primal-dual form, expressible as the packing-covering game defined by 
\begin{equation}\label{eq:packing-covering-game}\tag{PCG}
    \min_{\bx\in\Delta_m}\max_{\by\in\Delta_n} \bx^\top \mA\by
\end{equation}
for an entrywise-nonnegative matrix $\mA\in\reals^{m\times n}_+,$ and where $\Delta_m$ and $\Delta_n$ denote the probability simplices of dimension $m$ and $n,$ respectively. For completeness, \Cref{sec:LP-reduction} shows how approximate solutions to \eqref{eq:packing-covering-game} translate into approximate solutions to packing and covering LPs of the same order, and vice-versa. 

We consider games \eqref{eq:packing-covering-game} that are \emph{nontrivial}, where $\mA$ has no all-zero rows or columns. If $\mA$ had an all-zero row, the value of the game would be zero, as the $\bx$-player could place all its mass on the zero row, forcing the value of the game to zero irrespective of the $\by$-player's choice. The solution to the game is thus trivial. If $\mA$ had a zero column, it could be removed from the matrix without loss of generality (fixing the corresponding coordinate of $\by$ to zero), as $\by$-player maximizes $\bx^\top\mA\by$ and no maximizing strategy would place non-zero weight on the all-zero column, as $\mA$ is nonnegative. Throughout the paper, whenever we refer to \eqref{eq:packing-covering-game} we thus assume it is defined for a nontrivial matrix $\mA.$ In this case, the value of the game, $\val(\mA):= \min_{\bx\in\Delta_m}\max_{\by\in\Delta_n} \bx^\top \mA\by,$ is guaranteed to be strictly positive.    

Our main results are summarized in the following informal theorem; for more formal statements, see \Cref{thm:main-upper-bound} and \Cref{thm: main lower}.

\begin{theorem}[Main theorem; informal]\label{thm:main-informal}
    There is a deterministic algorithm that, given any nontrivial instance \eqref{eq:packing-covering-game} and a target error $\epsilon \in (0, 1)$, outputs $\bx_{\rm out} \in \Delta_m, \by_{\rm out} \in \Delta_n$ such that
    \begin{equation}\label{eq:PCG-mult-guarantee}
         \max_{j \in [n]}(\mA^\top \bx_{\rm out})_j \leq (1 + \epsilon) \min_{i \in [m]}(\mA\by_{\rm out})_i.
    \end{equation}
    The algorithm accesses $\mA$ through  $O\big(\frac{\log(n)\log(mn/\epsilon)}{\epsilon}\big)$ left and right matrix-vector products, and uses $\widetilde{O}\big(\frac{\nnz(\mA)}{\epsilon}\big)$ arithmetic operations. It is implementable in parallel with depth $\frac{\mathrm{polylog}(\nnz(\mA), 1/\epsilon)}{\epsilon}$ and $\widetilde{O}\big(\frac{\nnz(\mA)}{\epsilon}\big)$  work.    

    Moreover, any deterministic algorithm for \eqref{eq:packing-covering-game} that only accesses $\mA$ through left and right matrix-vector product oracle queries must make at least $\Omega\big(\frac{\log(1/\epsilon)}{\epsilon}\big)$ such queries in order to output $\bx_{\rm out} \in \Delta_m, \by_{\rm out} \in \Delta_n$ that satisfy  \eqref{eq:PCG-mult-guarantee},
    in the worst case over instances of \eqref{eq:packing-covering-game} of size $m = n = \mathrm{poly}(1/\epsilon)$.
\end{theorem}

We make a few remarks regarding these results. 
In terms of the number of arithmetic operations that our algorithm takes, in addition to the $O\big(\frac{\log(n)\log(mn/\epsilon)}{\epsilon}\big)$ matrix-vector products (eight per algorithm iteration), the remaining work includes an auxiliary optimization procedure implemented with $O(n\log(n)\log^2(n/\epsilon))$ arithmetic operations per algorithm iteration. This term is of the lower order even for moderately dense matrices, thus the overall poly-logarithmic factor is rather modest. Moreover, we did not search for the most efficient implementation of the auxiliary procedure; only one that was near-linear-time.

In terms of the matrix-vector products, our lower and upper bounds agree up to a single logarithmic factor, hence both are nearly tight. For the considered class of problems, counting matrix-vector multiplications is a natural measure of algorithm complexity, as it dominates the computational effort of effectively all algorithms developed for general instances of this problem. Moreover, the matrix-vector products play the role of first-order or ``gradient'' information for the considered problem, thus they represent a natural oracle model for studying complexity of iterative methods, akin to information-based complexity of optimization \cite{nemirovsky1983problem,traub2003information}. Moreover, this model is in line with the oracle-based models used for studying the oracle complexity of linear system solving \cite{nemirovsky1991optimality,nemirovsky1992information} and the oracle complexity of solving general matrix games, with matrices that are not necessarily nonnegative \cite{KS26}. 

Finally, we highlight that how to even construct lower bounds for nonnegative problem instances is not clear a priori. Standard oracle lower bound constructions in optimization, pioneered in \cite{nemirovsky1983problem,yu2018lectures}, are based on revealing only a small number of useful directions of a target subspace containing the solution per oracle answer, forcing the algorithm to make many queries before it can learn sufficiently many target directions. This idea is not directly compatible with nonnegative instances, as it would effectively force the directions to have disjoint supports, which results in queries leaking too much information to the algorithm. Additionally, even without the nonnegativity requirement for the matrix, we note that nearly tight lower bounds for simplex-simplex matrix games were obtained only very recently \cite{KS26}.

\subsection{Packing-Covering Games Capture Packing and Covering LPs}\label{sec:LP-reduction}

We briefly comment on how the positive matrix game stated in \eqref{eq:packing-covering-game} captures packing and covering LPs in their standard form. In particular, without loss of generality (see, e.g., \cite{awerbuch2008stateless,AO19}), every packing LP and its dual covering LP can be stated in the following form:
   \begin{equation*}
\begin{aligned}
   \text{Packing LP:}\quad \max_{{\bx \geq \bzero_m}:{\mA^\top \bx \le \bone_n}}\;\; & \bone_m^\top \bx,
\end{aligned}
\qquad\qquad
\begin{aligned}
   \text{Covering LP:}\quad \min_{ \by \ge \bzero_n: \mA \by \ge \bone_m}\;\; & \bone_n^\top \by
\end{aligned}
\end{equation*}
where $\bone_{\cdot}$ and $\bzero_{\cdot}$ are the vectors of all ones and all zeros, of dimension indicated by the subscript. 

Given $\epsilon > 0,$ our upper bound result (Theorem~\ref{thm:main-upper-bound}) guarantees an output primal dual pair $(\bx_{\rm out}, \by_{\rm \out}) \in \Delta_m \times \Delta_n$ such that 
\begin{equation}\label{eq:approx-packing-covering}
    0 < \max_{j \in [n]}(\mA^\top \bx_{\rm out})_j \leq (1 + \epsilon) \min_{i \in [m]}(\mA\by_{\rm out})_i.
\end{equation}
Let $\bxh := \bx_{\rm out}/(\max_{j \in [n]}(\mA^\top \bx_{\rm out})_j)$ and $\byh := \by_{\rm out}/(\min_{i \in [m]}(\mA\by_{\rm out})_i).$ It is immediate that $\bxh$ is feasible for the packing LP with matrix $\mA$ and similarly $\byh$ is feasible for its dual covering LP with matrix $\mA.$ Let $\bx_*$ be the optimal solution to the primal packing LP and $\by_*$ be the optimal solution to its dual covering LP. Then \[\bone_m^\top \bxh \leq \bone_m^\top \bx_* = \bone_n^\top \by_* \leq \bone_n^\top \byh.\] 
Since $(\bx_{\rm out}, \by_{\rm \out}) \in \Delta_m \times \Delta_n$, \eqref{eq:approx-packing-covering} and the definitions of $\bxh, \byh$ imply 
\[
    \bone_m^\top \bxh = \frac{1}{\max_{j \in [n]}(\mA^\top \bx_{\rm out})_j} \geq \frac{1}{1+\epsilon}\cdot\frac{1}{\min_{i \in [m]}(\mA\by_{\rm out})_i} = \frac{1}{1+\epsilon} \bone_n^\top \byh. 
\]
Thus, $\bone_m^\top \bxh \geq \frac{1}{1+\epsilon}\bone_m^\top \bx_*$ and $\bone_n^\top \byh \leq (1+\epsilon) \bone_n^\top \by_*$, so we get that $\bxh$ and $\byh$ are $(1+\epsilon)$-approximate solutions to the packing LP and the covering LP, respectively. 

On the lower bound side, we prove a stronger guarantee than stated in the informal \Cref{thm:main-informal}; namely, that for any output $\by_{\rm out} \in \Delta_n$ (where the instance is with $m = n$) after up to $T = c\frac{\log(1/\epsilon)}{\epsilon}$ queries, for a universal constant $c > 0$, it must be (see \Cref{thm: main lower}):
\begin{equation}\label{eq:lb-guarantee-intro}
    \min_{i \in [m]}(\mA \by_{\rm out})_i <\frac{\min_{\bx \in \Delta_m}\max_{\by \in \Delta_n}\bx^\top \mA \by}{1+\epsilon}. 
\end{equation}
A covering LP algorithm that outputs $\byh$ such that  $\bone_n^\top \byh \leq (1+\epsilon) \bone_n^\top \by_*$  and $\mA\byh \geq \bone_m$ leads to $\by' = \frac{\byh}{\bone_n^\top \byh}$ such that 
\[
\min_{i \in [m]}(\mA \by')_i = \frac{\min_{i \in [m]}(\mA \byh)_i}{\bone_n^\top \byh}\geq \frac{1}{\bone_n^\top \byh} \geq \frac{1}{(1+\epsilon) \bone_n^\top \by_*} = \frac{\max_{j\in [n]}(\mA^\top \bx_*)_j}{(1+\epsilon)\bone_m^\top \bx_*} \geq \frac{1}{1 + \epsilon}\min_{\bx \in \Delta_m}\max_{\by \in \Delta_n} \bx^\top \mA \by,
\]
contradicting \eqref{eq:lb-guarantee-intro}. Hence, such an output cannot be constructed with fewer than $T$ queries to $\mA.$

It is possible to argue (with more technical work) that the claimed oracle lower bound also extends to packing LP solvers, as we show in Remark~\ref{rem:packing-lower-bound}.

\subsection{Technical Overview}

We review the main ideas for our upper bound (\Cref{thm:main-upper-bound}) and our lower bound (\Cref{thm: main lower}) and summarize our AI methodology at the end of this subsection. As summarized in \Cref{sec:prob-and-main-results}, both the upper bound and the lower bound are for the packing-covering game \eqref{eq:packing-covering-game} and, as discussed in \Cref{sec:LP-reduction}, also apply to general packing and covering LPs. Our choice to work with the simplex-simplex game version of the problem was motivated by the fixed scale of the problem and the ability to capture both the packing LP and its dual covering LP problem with one formulation.

\paragraph{Upper bound.} 

The first core idea in our approach is to reformulate \eqref{eq:packing-covering-game} to a problem that intrinsically captures its scale invariance. An immediate observation is that because the logarithm function turns products into sums, incorporating it into the problem would make the additive error invariant to scale. To preserve the convex-concave property of the min-max objective, which retains problem tractability, our first step is to take the logarithm of the ``covering constraints'' $\mA \by$, leading to the following linear-log game:
\begin{equation}\label{eq:linear-log-game-intro}
    \min_{\bx \in \Delta_m}\max_{\by \in \Delta_n} \bx^\top \log(\mA \by),
\end{equation}
where the logarithm is taken entrywise. It does not take much effort to then argue that both \eqref{eq:linear-log-game-intro} and \eqref{eq:packing-covering-game} are solved by the same primal-dual pair $(\bx_*, \by_*)$ (see \Cref{sec:reduction}). Additionally, if $\mA$ is rescaled by a constant $c>0,$ since $\bx \in \Delta_m,$ this only shifts the objective by an additive $\log(c),$ which has no impact on 
the additive error. Moreover, since taking the primal-dual gap of \eqref{eq:linear-log-game-intro} cancels out the $\log(c)$ shift, it is entirely irrelevant to any analysis that bounds the primal-dual gap or its surrogates like regret.  

The more intricate part of the argument is to relate the additive error for \eqref{eq:linear-log-game-intro} to the desired multiplicative guarantee stated in \eqref{eq:PCG-mult-guarantee}. Unfortunately, order-$\epsilon$ primal-dual gap of \eqref{eq:linear-log-game-intro} does \emph{not} imply a $1 + O(\epsilon)$ multiplicative gap for \eqref{eq:packing-covering-game}. To see this, consider the following $2\times 2$ example:
\begin{equation}\label{eq:2x2-example}
    \mA = \begin{pmatrix}
        2 & 1\\
        1 & 2
    \end{pmatrix},\quad
    \bx = \frac{1}{2}\begin{pmatrix}
        1 + \delta\\
        1 - \delta
    \end{pmatrix}, \quad 
    \by = \frac{1}{2}\begin{pmatrix}
        1\\
        1
    \end{pmatrix}, \quad
    \delta \in (0, 1/3].
\end{equation}
Elementary calculus bounds the primal-dual gap of \eqref{eq:linear-log-game-intro} at $(\bx, \by)$ as
\begin{equation}\notag
    \max_{\bv \in \Delta_2} \bx^\top \log(\mA \bv) - \min_{\bu \in \Delta_2} \bu^\top \log(\mA \by) = \frac{1 + \delta}{2}\log(1 + \delta) + \frac{1-\delta}{2}\log(1 - \delta) \leq \delta^2.
\end{equation}
On the other hand, the multiplicative error of the form stated in \Cref{thm:main-informal} is 
\[
    \frac{\max_{j \in [2]}(\mA^\top \bx)_j}{\min_{i \in [2]}(\mA \by)_i} = \frac{(3 + \delta)/2}{3/2} = 1 + \frac{\delta}{3}.
\]
Hence, taking $\delta = \sqrt{\epsilon}$ for $\epsilon \in (0, 1/9]$ gives primal-dual gap at most $\epsilon$ for $\eqref{eq:linear-log-game-intro}$ but only a $1 + \sqrt{\epsilon}/3 = 1 + \Theta(\sqrt{\epsilon})$ multiplicative error for \eqref{eq:packing-covering-game}. 

While direct translation of additive error for \eqref{eq:linear-log-game-intro} to multiplicative error of \eqref{eq:packing-covering-game} is insufficient, we argue that a primal-dual pair with target multiplicative error for \eqref{eq:packing-covering-game} can be inferred from a sequence of primal-dual pairs $(\bx_k, \by_k)_{0\leq k \leq K}$, with target error on the primal-dual \emph{regret}, which bounds above the primal-dual gap by using that the objective is convex (in fact, linear) in the first argument and concave in the second argument. In particular, letting $\Lcal(\bx, \by) = \bx^\top\log(\mA\by)$, we prove (in \Cref{lemma:PCG-to-Minty-gap}) that it suffices to argue that 
\begin{equation}\label{eq:min-max-regret}
    \sup_{\substack{\bu \in \Delta_m,\\ \bv \in \Deltatr_n}} \frac{1}{K+1}\sum_{k=0}^K \Big(\big(\nabla_{\bx}\Lcal(\bx_k, \by_k)\big)^\top (\bx_k - \bu) - \big( \nabla_{\by}\Lcal(\bx_k, \by_k)\big)^\top (\by_k - \bv) \Big) = O(\epsilon),
\end{equation}
where for technical reasons (namely, to prevent the logarithm from becoming too large in magnitude for small values of $\by$-coordinates), the supremum is taken over the truncated simplex $\Deltatr_n := \{\by \in \Delta_n: \by \geq \epstr \bone_n\}$ with $\epstr = \Theta(\eps/n),$ and we similarly ensure that the sequence $\by_k$ is also supported on $\Deltatr_n$. 

For notational convenience and easier application of results from prior work, we use a standard interpretation of \eqref{eq:min-max-regret} as the regret of a Minty variational inequality associated with a monotone operator $\mF$ that stacks $\nabla_{\bx}\Lcal(\bx_k, \by_k)$ and $-\nabla_{\by}\Lcal(\bx_k, \by_k)$ \cite{Nem04}; see \Cref{sec:prelims} for relevant definitions. 

The primary challenge in developing an iterative procedure that bounds \eqref{eq:min-max-regret} is that the operator $\mF$, although monotone, does not satisfy the usual Lipschitz continuity (at least not with a sufficiently small Lipschitz parameter) for a bounded-range, strongly convex distance generating function, both properties being crucial for standard approaches to lead to iterative methods with the number of iterations that not only scales linearly with $1/\epsilon,$ but has poly-logarithmic scaling with the problem size and no dependence on the numerical values associated with the problem formulation. Here, the distance generating function is used in defining the proximal (prox) updates in methods like the classical mirror-prox \cite{Nem04}. 

Our second core idea is to instead look for a distance generating function $\psi$ with respect to which the operator $\mF$ satisfies a less stringent regularity condition; namely, the relaxed relative Lipschitzness condition introduced by \cite{JT23}. Under this condition, \cite{JT23} have shown that a variant of classical mirror-prox method \cite{Nem04} converges to regret at most $\epsilon$ within $O\big(\frac{L \Omega_{\psi}(\bz_0)}{\epsilon}\big)$ iterations, where $L > 0$ is the parameter for which relaxed relative Lipschitzness condition holds for $\mF$ with respect to $\psi$ and $\Omega_{\psi}(\bz_0)$ is the range of Bregman divergence of $\psi$ anchored at $\bz_0$ over the feasible set (see \Cref{def:Bregman-divergence} for the definition of Bregman divergence). The main challenge now becomes to find a jointly convex in $\bx, \by$ function $\psi$ for which: (1) $\mF$ is $L$-relaxed relative Lipschitz with $L$  poly-logarithmically bounded in the input size, (2) $\Omega_{\psi}(\bz_0)$ is bounded poly-logarithmically in the input size, for $\bz_0$ that can be constructed efficiently, and (3) algorithm prox updates, which minimize the sum of $\psi$ and linear terms over the feasible set for $\bx, \by$, can be implemented efficiently. The core technical contribution lies in addressing (1)--(3). We remark here that for (3), it suffices to have an approximate solution to prox update problems, with error of the order $\epsilon/L;$ for completeness, we provide the necessary details in \Cref{sec:generic-VI}.  

To construct $\psi,$ we first note that the most commonly used distance generating function under the simplex constraints is the negative entropy $h(\bz) := \bz^\top\log(\bz).$ However, a simple sum of $h(\bx)$ and $h(\by)$ is insufficient here, as $\mF$ does not satisfy even the relaxed Lipschitz condition with respect to this function. We instead consider a close relative $\varphi(\bx, \by)$ that mimics the structure of the linear-log game objective in \eqref{eq:linear-log-game-intro}, given by the sum of $h(\bx)$ and the relative entropy $\bx^\top(\log(\bx) - \log(\mA\by)),$ which is jointly convex in $\bx, \by$ by standard arguments. We argue that $\mF$ satisfies relaxed relative Lipschitzness with respect to this function, with a parameter that is a universal constant. Unfortunately, the computation of the prox updates is unclear for this function; in particular, it is unclear that it is any easier than solving \eqref{eq:linear-log-game-intro}. 

Nevertheless, the considered entropy + relative entropy auxiliary function $\varphi$ provides useful insights for constructing a function $\psi$ with the required properties. A key insight is that the relative entropy heavily features in the optimal transport literature, where it appears in various instantiations of the classical Sinkhorn algorithm; see, e.g., \cite{CPSV18}. In particular, the same type of function arises from minimizing a coupling (a nonnegative matrix $\bpi \in \R^{m \times n}$ whose row and column marginals, $\bpi \bone_n$, $\bpi^\top\bone_m$ give the coupled probability vectors) with one of the marginals being fixed. This insight allows us to look for $\psi$ that represents the coupling of $\bx, \by$ using an auxiliary matrix $\bpi$. Our distance generating function $\psi$, constructed in \Cref{sec:dgf}, can be seen as a modification of the initial $\varphi$  function that regularizes minimization over $\bpi$ and replaces $h(\bx)$ by $h(\by)$ to sufficiently decouple $\bx$ and $\by$, which is key in ensuring that prox updates can be implemented efficiently. We argue in \Cref{sec:dgf} that the constructed $\psi$ satisfies (1) with $L = O(\log(n))$.

For property (2), we use the first-order characterization of minimizing couplings $\bpi$ appearing in the definition of $\psi$ to construct an initial vector $\bz_0 = (\bx_0, \by_0)$ for which both $\bz_0$ and its corresponding minimizing coupling $\bpi_0$ are expressible in closed form. We then argue in \Cref{sec:initialization} that for constructed initial $\bz_0$, which can be computed with a constant number of (namely, three) matrix-vector products and with work $O(\nnz(\mA)),$ we have $\Omega_{\psi}(\bz_0) = O(\log(mn/\epsilon))$ (here, $\epsilon$ appears due to the $\epstr$ truncation of the simplex $\Delta_n$). 

Finally, for property (3), in \Cref{sec:implementation} we argue that the required prox iterations can be implemented with two matrix-vector products per prox update and in near-linear time in the problem size. The key step for this part is to argue that minimization of $\psi(\bx, \by)$ and linear functions in $\bx, \by$, which involves simultaneous minimization over $\bx, \by, \bpi,$ can be equivalently written as minimization of a convex objective $f(\bq, \by)$ for $\bq$ supported on the simplex $\Delta_n$ and $\by$ supported on the truncated simplex $\Deltatr_n.$ This minimization problem has $2n$ variables and it is highly structured. In particular, $\min_{\by \in \Deltatr_n}f(\bq, \by)$ can be computed to error $\delta >0 $ with $n\,  \mathrm{polylog}(n/\delta)$ arithmetic operations using known, parallelizable procedures, while $\widehat{f}(\bq) := \min_{\by \in \Deltatr_n}f(\bq, \by)$ possesses ``relative'' smoothness and strong convexity properties with respect to $h(\bq)$ that lead to geometrically fast convergence of the standard entropic mirror descent (a.k.a.\ the multiplicative weights update) algorithm, based on known results \cite{bauschke2017descent,lu2018relatively}. 

\paragraph{Lower bound.} A starting point for our lower bound is the additive error lower bound for general matrix games from \cite{KS26}. At a high level, the lower bound from \cite{KS26} is a type of a resisting oracle lower bound whose hard family of instances is expressible as the sum of rank-1 matrices $\bv_t \bw^{\top}_t.$ In particular, for each query $\bx_t, \by_t$, $t \ge 0,$ the oracle returns $(\bx_t^\top \mA_t, \mA_t \by_t),$ where $\mA_t = \sum_{\tau = 0}^t \bv_{\tau}\bw_\tau^\top$ and $\bv_t$ is orthogonal to all left queries $(\bx_0, \dots, \bx_t)$ in hindsight and all $\bv_{\tau},$ $\tau < t,$ while $\bw_t$ is orthogonal to all right queries $(\by_0, \dots, \by_t)$ and all $\bw_{\tau},$ $\tau < t.$ This orthogonality property ensures that the returned queries remain consistent with the final matrix after $T$ queries, $\mA := \mA_T$, while forcing the algorithm to ask $T$ queries before it can reveal sufficient information about $\mA$ to construct a solution with  primal-dual gap at most $\epsilon$. Of course, additional technical work goes into arguing that for the selected choice of vectors $\bv_t, \bw_t$ (in \cite{KS26}, projected Gaussians) and a sufficiently large $T = T(\epsilon)$, the primal-dual gap is at least $\epsilon.$

The first instinct when trying to construct a lower bound for \eqref{eq:packing-covering-game} may be to use a similar construction. However, it does not take long to see that using the same or similarly orthogonal instance family is futile here, due to the nonnegativity requirement for $\mA.$ In particular, a single dense query (e.g., $\bx = (1/m)\bone_m, \by = (1/n)\bone_n$) would force $\bv$ and $\bw$ to be zero vectors (to enforce orthogonality), making the matrix $\mA$ (and thus the problem) trivial for any $\epsilon > 0, m \geq 1, n \geq 1$.

Another, perhaps less obvious idea, is to ``spike'' the construction from \cite{KS26} by adding a rank-one matrix (e.g., $\eta \bone_m\bone_n^\top$ for sufficiently large $\eta$) to ensure the constructed matrix (and all ``hindsight'' matrices $\mA_t$) is entrywise nonnegative. It is possible to argue that an adaptation of the argument from \cite{KS26} does lead to an informative lower bound in this case, albeit of the lower order $\Omega(\epsilon^{-2/3}).$ 

The key difficulty in strengthening this lower bound is the dependence between $\eta$ and the scale of $\bv_t, \bw_t$  that is needed to ensure that the constructed matrix is nonnegative while keeping the primal-dual gap sufficiently large relative to $\val(\mA).$ Recalling the scale invariance property discussed at the beginning of this section, what is at play here is that the relevant quantity---the matrix \emph{width}---is insufficiently large to enforce higher oracle complexity. This is the key insight: to have the hope of obtaining a stronger lower bound, the width of the constructed instance should grow with the inverse error $1/\eps$.

There are additional challenges when attempting to construct $\bv_t, \bw_t$ as subgaussian vectors common to prior work on oracle lower bounds in optimization; see, e.g., \cite{nemirovski1994parallel,diakonikolas2020lower,KS26}. On one hand, the randomness is useful because of known concentration and anticoncentration properties for standard choices of random vectors (e.g., i.i.d.\ Gaussian or i.i.d.\ Rademacher vectors) that are needed for both hiding the information (ensuring orthogonality or near-orthogonality) and for ensuring that the primal-dual gap remains sufficiently large over a sequence of queries. On the other hand, it becomes hard to control how negative some of the entries of $\sum_{t=0}^T \bv_t \bw_t^\top$ may become, forcing a large correction that brings us back to the small-width issue. However, what these difficulties reveal is that sparse, bounded-entry vectors $\bv_t, \bw_t$ should be useful here, at the very least to not cause large corrections for ensuring nonnegativity.

To simultaneously achieve all desired properties of the hard instance family, we abandon the prior approach of using random vectors. Instead, we construct $\bv_t$, $\bw_t$ as sparse vectors of known support, where each vector is supported on a contiguous small block of coordinates. 
To ensure these vectors are bounded below and maintain orthogonality to queries (and to each other, if supported on the same block), we explicitly enforce these properties as a feasibility condition on a polytope (defined as the intersection of $[-1, 1]^n,$ the block subspace, and the orthogonality conditions; observe that orthogonality to a set of vectors defines a linear subspace). The correction to the resulting ``sum of orthogonal rank ones'' matrix to ensure nonnegativity now becomes trivial by simply adding a block-type matrix with constant-valued entries.

Finally, to ensure that the constructed instance corresponds to a game with a nontrivial value that remains well-separated from possible algorithm outputs, we additionally spike the constructed matrix by adding a simple nonnegative matrix that ensures these properties. This leads to our hard instance family provided in \Cref{sec:instance-fam}. \Cref{sec:lb-proof} then formalizes the outlined intuition and observations for the constructed family through a sequence of lemmas that provide the necessary ingredients for the lower bound proof. A notable feature of our lower bound is that it does not require the queries $(\bx_t, \by_t)$ to be from $\Delta_m \times \Delta_n$ or even nonnegative; hence, designing algorithms that potentially utilize such more general queries cannot improve the oracle complexity. Finally, although our lower bound is stated for the packing-covering game \eqref{eq:packing-covering-game}, it is possible to show that it similarly implies an oracle lower bound (with matrix-vector product oracles) for individual packing LP and covering LP, as discussed in \Cref{sec:LP-reduction}.

\paragraph{AI methodology.} The main technical results in this paper were initially obtained through a long series of prompts to OpenAI GPT models, with authors supplying ideas, proof skeletons, and pointers to relevant techniques. The authors then fully wrote all the proofs from scratch, adding the intuition and connections to related work, improving the logical flow of technical arguments, and refining and tightening the technical results. AI was used in this latter process to test ideas and audit the proofs. The prose in the paper was written by the authors, with AI assisting in the literature search, copy-editing, and the $2\times 2$ example construction stated in \eqref{eq:2x2-example} and discussed earlier in this subsection.  

\subsection{Related Work}\label{sec:related}

\paragraph{Width-dependent algorithms.}
Classical Lagrangian-relaxation algorithms such as   \cite{PST95} and multi\-plicative-weights methods surveyed in \cite{AHK12} give approximate
packing and covering solvers whose iteration counts polynomially scale with the matrix width in the worst case. Other optimization methods applicable to matrix games such as Nesterov's smoothing~\cite{Nes05} and Nemirovski's mirror-prox \cite{Nem04} lead to additive error $\epsilon$ with $O\big(\frac{\nnz(\mA)\log(m n)\max_{ij}|A_{ij}|}{\epsilon}\big)$ work. Translated to multiplicative approximation, the resulting work bound scales polynomially  with $m, n, \width(\mA)$ in the worst case. 

\paragraph{Width-independent algorithms with nearly linear work.} \cite{LN93} initiated the line of work on width-independent algorithms with a parallel algorithm of $\frac{\mathrm{polylog}(m, n, 1/\epsilon)}{\epsilon^4}$ depth, which \cite{BBR04} then extended to distributed computation,
while \cite{AK09} developed stateless distributed
algorithms. \cite{AO15} introduced optimization-based techniques that improved the parallel accuracy dependence, and subsequent work reached $(1/\epsilon^2)$-dependence \cite{MRWZ16,DO17}. It is worth noting that the approach from \cite{AO15}, presented as scaling with $(1/\epsilon^3)$, could already give the $(1/\epsilon^2)$ scaling with step size larger by a factor $1/\epsilon$ than presented, without changing the analysis.

In the sequential setting, \cite{KY14} obtained
a randomized primal-dual solver with work
$O(\nnz(\mA)+(m+n)\log(\nnz(\mA))/\epsilon^2)$, giving a different tradeoff between input size and accuracy.
Subsequent sequential algorithms improved the accuracy dependence: \cite{AO19} obtained
$\widetilde{O}(\nnz(\mA)/\epsilon)$ work for packing and
$\widetilde{O}(\nnz(\mA)/\epsilon^{3/2})$ for covering, and \cite{WRM16} subsequently\footnote{\cite{AO19} is the journal version of a STOC 2015 paper.} obtained
$\widetilde{O}(\nnz(\mA)/\epsilon)$ work for covering as well. These results are based on randomized coordinate descent, and their depth crucially scales linearly with the size of the vector of variables, in addition to the $1/\epsilon$ scaling. 

Width-independent algorithms also exist for mixed packing and covering feasibility problems, which allow both types of constraints on the same vector of variables and are more general than solving a packing LP and its covering
dual. We refer to \cite{You01,You14,MRWZ16,CQ18,Qua20,OV26} for results in this line of work. None of these existing results gives  parallel depth scaling that is better than $1/\epsilon^2.$ 

\paragraph{Width-independent algorithms with additional dimension or sparsity
dependence.}
As emphasized in~\cite{AO19}, removing coefficient-width dependence
alone does not ensure nearly linear work. \cite{BI06} adapted Nesterov's smoothing to obtain a width-independent $1/\epsilon$ rate, with an additional polynomial
factor depending on dimension and row sparsity. For mixed packing-covering feasibility, \cite{BSW19} used Sherman's
area-convexity framework~\cite{She17}
to obtain parallel depth scaling with $(s/\epsilon)$ and
$\widetilde{O}(s \nnz(\mA)/\epsilon)$ total work, where $s$ is the maximum number of nonzero coefficients in a constraint. 

\paragraph{Oracle lower bounds and general matrix games.}
In the context of packing and covering LPs, \cite{KY15} proved iteration lower bounds for the 
class of Dantzig-Wolfe and Lagrangian-relaxation algorithms based on
linear optimization oracles. Those bounds concern a restricted algorithmic
framework and scale with $1/\epsilon^2$ and the matrix width. They are not applicable to any of the width-independent packing and covering LP solvers, which access the constraints through matrix-vector products. On the other hand, \cite{kyng2020packing} establish \emph{conditional} work lower bounds $\widetilde{\Omega}\big(\frac{\nnz(\mA)}{\epsilon^{0.165}}\big)$ for solving packing LPs under linear-system hardness assumptions. Our oracle lower bound is unconditional; it applies to every deterministic
algorithm with two-sided matrix-vector access and it is nearly tight in the considered matrix-vector access model. 

For general payoff matrices with $|A_{ij}|\le1$, the deterministic
matrix-vector oracle complexity of additive approximation has recently
been determined up to logarithmic factors. \cite{KS26} prove a $\widetilde{\Omega}(\epsilon^{-2/3})$ lower bound,
matched by the $\widetilde{O}(\epsilon^{-2/3})$ upper bound in~\cite{KOS26}. A recent preprint \cite{WF26} reports a $\widetilde{\Omega}(\epsilon^{-2/3})$ lower bound
for randomized algorithms as well. These results are for  additive
primal-dual error. Even for nonnegative payoff matrices, converting an
additive error to a multiplicative one requires error on the scale
of the game's value, which can be arbitrarily small after normalizing the
entries. Our lower bound instead characterizes
multiplicative approximation for nonnegative games without coefficient
width dependence.

\paragraph{Independent concurrent work.} 

\cite{CCGS26} independently obtained a
deterministic, width-independent algorithm for packing LPs with $\widetilde{O}\big(\frac{\nnz(\mA)}{\epsilon}\big)$ work and $\frac{\mathrm{polylog}(\nnz(\mA), 1/\epsilon)}{\epsilon}$ parallel depth. Their approach reduces packing to composite $\ell_{1,\infty}$
regression and uses parallel solvers for structured residual problems, building on techniques for positive $\ell_p$ regression. Their main theorem returns an approximate packing LP solution only. They also improve the oracle complexity of a composite $\ell_{1,\infty}$
regression problem introduced in \cite{chen2025accelerated}, which is not considered in our work. Our algorithm constructs approximate solutions to both packing and covering LP by solving a scale-invariant variational inequality that we introduce. We also prove a nearly matching deterministic matrix-vector oracle lower bound, a question not addressed in \cite{CCGS26}.

\subsection{Paper Organization}

The rest of the paper is organized as follows. \Cref{sec:prelims} introduces the notation and basic definitions and facts that are needed for our technical results. Our main upper bound result, including the core statements and proofs, is provided in \Cref{sec:main-upper}, while some of the proofs of auxiliary lemmas and claims are deferred to Appendix \ref{sec:appx}. Our main lower bound result is stated and proved in \Cref{sec:lower}. We conclude in \Cref{sec:conclusion}. 

\section{Notation and Preliminaries}\label{sec:prelims}

We use bold lowercase letters for vectors and bold uppercase letters for matrices. For vectors $\bx, \by,$ $\inner{\bx, \by} = \bx^\top \by$ denotes the standard inner product. For matrices $\mA, \mB,$ $\inner{\mA, \mB} = \mathrm{Trace}(\mA^\top \mB)$. We use $\Delta_n := \{\bx \in \R_+^n: \bone_n^\top\bx = 1\}$ to denote the probability simplex, where $\R_+$ denotes the nonnegative orthant. For a set $\Ccal,$ $\relint(\Ccal)$ denotes its relative interior; in particular, $\relint(\Delta_n) = \{\bx \in \R^n: \bx > \bzero_n, \bone^\top \bx = 1\}$, where $'>'$ is understood entrywise. 
We use $\log$ to denote the natural logarithm. For a vector $\bx = (x_1, \dots, x_n)^\top,$ $\log(\bx) := (\log(x_1),\dots, \log(x_n))^\top,$ and the same notation is used for entrywise application to matrices. When the context is clear, we omit explicitly writing the transpose, and simply write $\bx = (x_1, \dots, x_n)$ for a column vector.  Throughout the paper, the notation $h(\bx) = \bx^\top \log(\bx)$ is reserved for the (negative) entropy function, where $\bx$ is entrywise nonnegative and $0 \log 0 = 0$, by convention. When applied to an entrywise nonnegative matrix $\bpi \in \R_+^{m \times n},$ $h$ is understood as acting on the vectorization of $\bpi,$ so that $h(\bpi) = \sum_{i,j}\Pi_{ij}\log(\Pi_{ij} ) = \inner{\bpi, \log(\bpi)}$. For two vectors $\bx = (x_1, \dots, x_n)$ and $\by = (y_1, \dots, y_n)$ of the same dimension, we define entrywise product $\odot$ and division $\oslash$ by $\bx \odot \by := (x_1 y_1, \dots, x_n y_n)$ and $\bx \oslash \by := (x_1/y_1, \dots, x_n /y_n).$  The gradient of a  function $f$, assuming differentiable, is denoted by $\nabla f.$ Bregman divergence of a function is defined in the usual way, as follows.

\begin{definition}[Bregman divergence]\label{def:Bregman-divergence}
    Let $f$ be a function that is differentiable on the interior of its domain. Then, the Bregman divergence of $f$ at any $\by \in \mathrm{dom}(f)$ and $\bx \in \mathrm{int\,dom}(f)$ is denoted by $D_f(\by, \bx)$ and defined by:
    \[
    D_f(\by, \bx) = f(\by) - f(\bx) - \inner{\nabla f(\bx), \by - \bx}. 
    \]
\end{definition}
Bregman divergence of a function $f$ is the difference between the function itself ($f(\by)$) and its affine approximation at $\bx$ ($f(\bx) + \inner{\nabla f(\bx), \by - \bx}$). It is thus immediate that the Bregman divergence of an affine function is zero and of a convex function is nonnegative. A differentiable function $f$ is $\mu$-strongly convex w.r.t.\ a norm $\norm{\cdot}$ if and only if $D_f(\by, \bx) \geq \frac{\mu}{2}\norm{\by - \bx}^2$ for all $\bx, \by.$ It is a standard fact (see, e.g., \cite[Example 5.27]{beck2017first}) that $h$ is $1$-strongly convex w.r.t.\ the $\ell_1$ norm on the probability simplex.

Bregman divergence of a function is also useful for defining generalizations of properties such as standard Lipschitzness, smoothness (gradient Lipschitzness), and  strong convexity. Relative smoothness and relative strong convexity are summarized in the following definition.

 \begin{definition}[Relative smoothness and relative strong convexity]\label{def:rel-smooth-rel-str-cvx}
     Let functions $f$ and $g$ be well-defined on the same set $\Ccal$ and differentiable on its relative interior. Given $L > 0,$ we say that $f$ is $L$-relatively smooth with respect to $g$ on $\Ccal$ if for all $\by \in \Ccal,$ $\bx \in \relint(\Ccal),$ $D_f(\by, \bx) \leq L D_g(\by, \bx).$ Similarly, given $\mu > 0,$ we say that $f$ is $\mu$-relatively strongly convex with respect to $g$ if for all $\by \in \Ccal,$ $\bx \in \relint(\Ccal),$ we have $D_f(\by, \bx) \geq \mu D_g(\by, \bx).$
 \end{definition}

A standard fact about Bregman divergence of the entropy---a.k.a.\ Kullback-Leibler (KL) divergence---is that it is jointly convex in both arguments, a property that does not necessarily hold for Bregman divergences of other (strongly) convex functions. A closely related relative entropy function (equal to KL divergence up to linear terms) is likewise convex.

\begin{fact}[E.g., Example 3.19 in \cite{vandenberghe2004convex}]\label{fact:convexity-of-KL}
    Let $\bx \geq \bzero,$ $\bu > \bzero.$ Then, the function $D_h(\bx, \bu) = h(\bx) - h(\bu) - \nabla h(\bu)^\top(\bx - \bu)$ is jointly convex in its arguments, and so is the relative entropy defined by $\bx^\top \log (\bx \oslash \bu)$. 
\end{fact}

Additionally, we will make use of the standard three-point inequality that arises in the proofs of mirror descent-style methods, summarized in the following fact. The proof is an immediate consequence of the first-order optimality condition for the stated problem and the definition of a Bregman divergence.

\begin{fact}[Three-point inequality]\label{fact:3PI}
    Let $f$ be a convex function whose domain contains a closed nonempty convex set $\Ccal \subseteq \R^n$ and let 
    \begin{equation}\label{eq:prox-map-generic}
        \bxh \in \argmin{\by \in \Ccal}\{\inner{\ba, \by} + D_f(\by, \bx)\},
    \end{equation}
    where $\ba \in \R^n$ is fixed. 
    Then, for any $\by \in \Ccal,$ 
    \[
    \inner{\ba, \by - \bxh} \geq D_f(\by, \bxh) - D_f(\by, \bx) + D_f(\bxh, \bx).  
    \]    
\end{fact}
Given $\epsilon > 0,$ we define the notation $\approx_\epsilon$ via 
\begin{equation}\label{eq:approx-3PI}
   \bxh \approx_\epsilon \argmin{\by \in \Ccal}\{\inner{\ba, \by} + D_f(\by, \bx)\} \;\Leftrightarrow\; \inner{\ba, \by - \bxh} \geq D_f(\by, \bxh) - D_f(\by, \bx) + D_f(\bxh, \bx) - \epsilon,
   \quad\forall\by\in \Ccal. 
\end{equation}
 In \eqref{eq:prox-map-generic}, $\bx$ is often referred to as the prox center, and the update $\bxh$ as the prox update. Note that when $\Ccal$ is in addition bounded and $f$ is strictly convex, the argmin in \eqref{eq:prox-map-generic} is a singleton set, so $\bxh$ is uniquely defined. Additionally, when $\Ccal = \Delta_n$ and $f = h,$ $\bxh \in \relint(\Delta_n).$ 

We will further rely on well-known properties of functions expressible as minima of functions of two arguments. First, the infimum over one of the arguments of a jointly convex function of two arguments minimized over a convex set is convex.

\begin{fact}[E.g., Section 3.2.5 in  \cite{vandenberghe2004convex}]\label{fact:inf-convexity}
    Let $f(\bx, \by)$ be a jointly convex function in both arguments and let $\Ccal$ be a nonempty convex set. Then $g(\bx) = \inf_{\by \in \Ccal}f(\bx, \by)$ is convex in $\bx$, so long as $g(\bx) >-\infty$ for all $\bx.$ 
\end{fact}

The following fact is a consequence of the classical theorem of Danskin \cite{danskin2012theory}. The version stated below is a special case of the more general theorem provided in \cite[Theorem 4.1]{bonnans1998optimization} that is nevertheless sufficient for our application. 

\begin{fact}[Danskin's Theorem, e.g., Theorem 4.1 in \cite{bonnans1998optimization}]\label{fact:danskin}
    Let $\Xcal \subseteq \R^d,\, \Ycal \subseteq \R^k$ be nonempty compact convex sets and let $f : \Xcal \times \Ycal \to \R$ be continuous. Suppose that:
    \begin{enumerate}[label = (\arabic*)]
        \item For every $\by \in \Ycal,$ the function $f(\cdot, \by)$ is differentiable on $\relint(\Xcal)$ and $(\bx, \by) \mapsto \nabla_\bx f(\bx, \by)$ is continuous on $\relint(\Xcal) \times \Ycal.$ 
        \item For every $\bx \in \relint(\Xcal),$ the minimizer $\by^*(\bx) = \argmin{\by \in \Ycal}f(\bx, \by)$ is unique.
    \end{enumerate}
    Then the minimizing function $g(\bx):= \min_{\by \in \Ycal}f(\bx, \by)$ is continuously differentiable on $\relint(\Xcal)$ and
    \[
    \nabla g(\bx) = \nabla_\bx f(\bx, \by^*(\bx)). 
    \]
\end{fact}

We will also make use of the following elementary fact.

\begin{fact}\label{fact:shifted-entropy-expression-convexity}
    Let $p \geq 0, t > 0,$ and $s > 0.$ Then:
    \[
    (p + t)\log(p + t) - p \log(s) + s\log(s) = D_h(p+t, s) + (s + t)\log(s) + p + t - s.
    \]
    Moreover, the function in the display above is convex for $s \ge t.$
\end{fact}

\paragraph{Minty variational inequality.} The analysis of our algorithm reduces to solving a Minty variational inequality (MVI). Given a nonempty closed convex set $\Zcal \subseteq \R^d$ and an operator $\mF: \Zcal \to \R^d,$ Minty variational inequality asks for $\bz_* \in \Zcal$ such that
\begin{equation}\label{eq:MVI}\tag{MVI}
    \inner{\mF(\bz), \bz_* - \bz} \leq 0, \quad \forall \bz \in \Zcal. 
\end{equation}
Given $\epsilon > 0,$ we say that $\bz_\epsilon$ is an $\epsilon$-approximate solution to \eqref{eq:MVI} if for all $\bz \in \Zcal,$ $\inner{\mF(\bz), \bz_\epsilon - \bz} \leq \epsilon.$ For any $\bw \in \Zcal,$ we refer to $\sup_{\bz \in \Zcal}\inner{\mF(\bz), \bw - \bz}$ as the MVI gap of $\bw.$
We say that an operator $\mF$ is monotone on $\Zcal$ if for all $\bz, \bw \in \Zcal,$ $\inner{\mF(\bz) - \mF(\bw), \bz - \bw} \geq 0.$ Our argument is built for relaxed relatively Lipschitz operators, as defined below. 

\begin{definition}[Definition 1 in \cite{JT23}]\label{def:relaxed-rel-Lip}
We say that an operator $\mF: \Zcal \to \Zcal^*$ is $L$-relaxed relatively Lipschitz with respect to a function $\psi: \Zcal \to \R$ if for any triple $(\bu, \bz, \bw) \in \Zcal \times \Zcal \times \Zcal,$ we have
\begin{equation}\notag
    (\mF(\bw) - \mF(\bz))^\top(\bw - \bu) \leq L(D_\psi(\bw, \bz) + D_\psi(\bu, \bw) + D_\psi(\bu, \bz)). 
\end{equation}
\end{definition}

Finally, standard analyses of algorithms for MVIs typically produce a sequence of iterates $\bz_k,$ $k \geq 0$, and bound regret-like quantities: for $a_k > 0$, $k \geq 0,$ and $A_k := \sum_{i=0}^k a_i,$ the analyses will typically bound
 $   \frac{1}{A_k}\sum_{i=0}^k a_i \inner{\mF(\bz_i), \bz_i - \bz}$ 
for an arbitrary but fixed $\bz \in \Zcal,$ and then use monotonicity of $\mF$ to bound the MVI gap of $\bzh_k := \frac{1}{A_k}\sum_{i=0}^k a_i\bz_i$, since, by monotonicity of $\mF,$ 
\[\inner{\mF(\bz), \bzh_k - \bz} = \frac{1}{A_k}\sum_{i=0}^k a_i \inner{\mF(\bz), \bz_i - \bz} \leq \frac{1}{A_k}\sum_{i=0}^k a_i \inner{\mF(\bz_i), \bz_i - \bz}.\]
The bound on the MVI gap then follows by taking the supremum over $\bz \in \Zcal.$ 
For operators that are monotone (but do not satisfy stronger properties than monotonicity like strong monotonicity), $a_k$ is typically constant and the analyses of iterative methods bound the usual regret 
\begin{equation}\label{eq:MVI-regret}
    \frac{1}{k + 1}\sum_{i=0}^k \inner{\mF(\bz_i), \bz_i - \bz}.
\end{equation}
The guarantee of our algorithm will be based on a regret bound for a nonlinear monotone operator $\mF$ constructed based on \eqref{eq:packing-covering-game}. 

\section{Efficient Iterative Algorithm}\label{sec:main-upper}

In this section, we provide our main upper bound result stated in \Cref{thm:main-upper-bound}. 
Throughout this
section, we assume without loss of generality that $n\geq 2$. Our main algorithm, which includes definitions of all relevant parameters, is stated in \Cref{alg:packing-covering-mirror-prox}. In the rest of the section, we prove \Cref{thm:main-upper-bound} through the following line of argument. In \Cref{sec:reduction}, 
we first argue that the target multiplicative error $(1+\epsilon)$ for the packing-covering game stated in \eqref{eq:packing-covering-game} can be reduced to the additive error $O(\epsilon)$ of an MVI with a monotone operator $\mF$ that we construct. In \Cref{sec:generic-VI}, we argue that, to solve the stated MVI with an iterative algorithm (relaxed mirror-prox), it suffices to construct a distance generating function (dgf) $\psi$, with respect to which the operator $\mF$ is $L$-relaxed relative Lipschitz (\Cref{def:relaxed-rel-Lip}). The relaxed mirror-prox then uses the number of iterations that scales linearly with $L/\epsilon$ and $\sup_{\bz \in \Zcal}D_{\psi}(\bz, \bz_0),$ where $\bz_0 = (\bx_0, \by_0)$ is the initial vector and $\Zcal \subseteq \Delta_m \times \Delta_n$ is the feasible set for the MVI. In \Cref{sec:dgf}, we construct $\psi$ for which $L  = O(\log(n))$.
In \Cref{sec:initialization}, we provide explicit initialization $\bz_0$ for which  $\sup_{\bz \in \Zcal}D_{\psi}(\bz, \bz_0) = O(\log(mn/\epsilon)).$ As a final ingredient, in \Cref{sec:implementation}, we provide efficient implementation of the mirror-prox steps. \Cref{sec:final-upper} then combines all these ingredients to complete the proof of \Cref{thm:main-upper-bound}.

\begin{theorem}\label{thm:main-upper-bound}
    Given an instance of \eqref{eq:packing-covering-game} with a matrix $\mA \in \R_+^{m \times n}$ with no all-zero rows or columns and target error $\epsilon \in (0, 1)$, \Cref{alg:packing-covering-mirror-prox} outputs $\bxh_K\in\Delta_m$ and $\byh_K\in\Delta_n$ such that
\[\max_{j\in[n]}(\mA^\top\bxh_K)_j
\leq(1+\epsilon)\min_{i\in[m]}(\mA\byh_K)_i.\]
The algorithm uses $O\big(\frac{\log(n)\log(m n/\epsilon)}{\epsilon}\big)$ matrix-vector multiplications and \[
O\Big(\frac{\log(n)\log(m n/\epsilon)}{\epsilon}\big(\nnz(\mA) + n \log(n)\log^2(n/\epsilon)\big)\Big) = \widetilde{O}\Big(\frac{\nnz(\mA)}{\epsilon}\Big)
\] 
arithmetic operations. It is implementable in parallel, in depth $\frac{\mathrm{polylog}(mn/\epsilon)}{\epsilon}$ with  $\widetilde{O}\big(\frac{\nnz(\mA)}{\epsilon}\big)$ total work. 
\end{theorem}

\begin{algorithm}[t!]
\caption{Approximate relaxed mirror-prox for packing-covering games}
\label{alg:packing-covering-mirror-prox}
\DontPrintSemicolon
\small
\KwIn{A matrix $\mA\in\R_+^{m\times n}$ with no zero rows or columns,
$n\geq2$, and accuracy $\epsilon\in(0,1)$.}
\KwOut{$\bxh_K\in\Delta_m$ and $\byh_K\in\Delta_n$ such that
$\max_{j\in[n]}(\mA^\top\bxh_K)_j
\leq(1+\epsilon)\min_{i\in[m]}(\mA\byh_K)_i$.}

\tcp{Set parameters.}
Set $\epstr\gets\epsilon/(8n)$,
$\beta\gets1/(1+\log n)$,
$L\gets7(3+4e)/\beta$, and $\xi\gets\epsilon/(12L)$\;
Set $\displaystyle
K\gets\left\lceil\frac{CL}{\epsilon}\log\frac{m}{\epstr}\right\rceil$,
where $C>0$ is a sufficiently large absolute constant. Initialize $\bqh_K \gets \bzero_m, \byh_K\gets \bzero_n.$\;

\tcp{Define the feasible set, operator, and prox function.}
Define $\Zcal:=\Delta_m\times\Deltatr_n$ where $\Deltatr_n:=\{\by\in\Delta_n:\by\geq\epstr\bone_n\}$\;
Define
$\displaystyle
\Delta_\mA:=
\left\{
    \bpi\in\R_+^{m\times n}:
    \bone_m^\top\bpi\bone_n=1,\quad
    A_{ij}=0\Longrightarrow\Pi_{ij}=0
\right\}
$\;
Define $\mF(\bx,\by):= \big(\log(\mA\by), 
       -\mA^\top(\bx\oslash(\mA\by))\big)$\; 
Define the prox function $\psi(\bx,\by):=
    h(\by)+
    \min_{\substack{\bpi\in\Delta_\mA\\\bpi\bone_n=\bx}}
    \big\{
        \sum_{i,j:A_{ij}>0}
        \Pi_{ij}\log\frac{\Pi_{ij}}{A_{ij}y_j}
        +\beta h(\bpi^\top\bone_m+\epstr\bone_n)
    \big\}$

Compute initialization $\bz_0=(\bx_0,\by_0)$ as in
\Cref{lemma:initialization}\;

\tcp{Implement the two following prox steps using \Cref{sec:implementation}.}
\For{$k=0,\ldots,K$}{
    \tcp{Relaxed inexact mirror-prox update. See \Cref{algo:approx-prox} for a procedure that implements approximate minimization in Lines 9 and 11.
    }

    Find $\bw_k=(\bu_k,\bv_k)\in\Zcal$, with $\bu_k>\bzero$,
    such that, for every $\bz\in\Zcal$,
    \(
        \big(
            \frac{1}{L}\mF(\bz_k)
            +\nabla\psi(\bw_k)-\nabla\psi(\bz_k)
        \big)^\top(\bw_k-\bz)\leq\xi.
    \)\;

    Update outputs $\bqh_K\gets \bqh_K + \frac{1}{K+1}\bu_k\oslash(\mA\bv_k)$,\quad $\byh_K \gets \byh_K + \bv_k/(K+1)$\;

    Find $\bz_{k+1}=(\bx_{k+1},\by_{k+1})\in\Zcal$,
    with $\bx_{k+1}>\bzero$,
    such that, for every $\bz\in\Zcal$,
    \(
        \big(
            \frac{1}{2L}\mF(\bw_k)
            +\nabla\psi(\bz_{k+1})-\nabla\psi(\bz_k)
        \big)^\top(\bz_{k+1}-\bz)\leq\xi.
    \)\;
}

\Return{$\displaystyle
    \bxh_K=\frac{\bqh_K}{\bone_m^\top\bqh_K},
    \quad
    \byh_K
$}\;
\end{algorithm}

\subsection{Reduction to a Variational Inequality}\label{sec:reduction}

As a first step in our analysis, we show how to reduce the original packing-covering game to an MVI with a nonlinear monotone operator, with the property that the MVI regret (see \eqref{eq:MVI-regret}) being bounded by $\epsilon > 0$ implies $(1 + O(\epsilon))$-multiplicative error for the original packing-covering game \eqref{eq:packing-covering-game}. 

Consider first the following linear-log game:
\begin{equation}\label{eq:log-game}
    \min_{\bx \in \Delta_m}\max_{\by \in \Delta_n} \Lcal(\bx, \by), \;\text{ where }\; \Lcal(\bx, \by) := \bx^\top \log(\mA\by).
\end{equation}

The insight behind the linear-log game in \eqref{eq:log-game} is that it is scale invariant under \emph{additive} error. This invariance becomes more useful for the analysis of iterative methods, which typically bound additive error. By contrast, \eqref{eq:packing-covering-game} is \emph{scale-invariant}, so its invariance is with respect to the \emph{multiplicative} error. 

In particular, in the original problem \eqref{eq:packing-covering-game}, rescaling all elements of $\mA$ by the same constant $c > 0$ changes the value of \eqref{eq:packing-covering-game} by the same constant, but has no impact on the solutions to the game or the multiplicative error. However, this rescaling \emph{does} affect the additive error when solving \eqref{eq:packing-covering-game}, which is the natural notion of error that first-order methods track. By contrast, such a rescaling has no impact on the additive error when solving \eqref{eq:log-game}. In particular, $\bx^\top \log(c\mA\by) = \bx^\top\log(\mA\by) + \log c,$ which shifts the value of the program and does not have any impact on the additive error when solving it. 

Moreover, \eqref{eq:log-game} is solved by the same primal-dual pair as \eqref{eq:packing-covering-game}. To see this, suppose $(\bx_*, \by_*)$ solves \eqref{eq:packing-covering-game}. We claim that for an arbitrary primal-dual pair $\bu \in \Delta_m, \bv \in \Delta_n$, the primal-dual gap $\Lcal(\bx_*, \bv) - \Lcal(\bu, \by_*)$ is non-positive, and thus $(\bx_*, \by_*)$ solves \eqref{eq:log-game}. First, by Jensen's inequality (since $\bx_* \in \Delta_m$ and log is concave), we have 
\[
\Lcal(\bx_*, \bv) \leq \log(\bx_*^\top \mA\bv) \leq \max_{\bv' \in \Delta_n} \log(\bx_*^\top \mA\bv') =  \log(\max_{\bv' \in \Delta_n} \bx_*^\top \mA\bv') = \log(\bx_*^\top \mA \by_*),
\]
where we have used the strict monotonicity of logarithm (to exchange $\max$ and $\log$). On the other hand, 
\[
\Lcal(\bu, \by_*) = \bu^\top \log(\mA\by_*) \geq \min_{\bu' \in \Delta_m} (\bu')^{\top} \log(\mA\by_*) = \min_i \log((\mA\by_*)_i) = \log(\min_{i}(\mA\by_*)_i) = \log(\bx_*^\top \mA \by_*).
\]
Thus $\Lcal(\bx_*, \bv) - \Lcal(\bu, \by_*) \leq 0.$ 

Instead of working directly with the linear-log game in \eqref{eq:log-game}, it will be more convenient to work with the MVI associated with it. In particular, we will control the multiplicative error for the original problem in \eqref{eq:packing-covering-game} by the MVI regret (see \eqref{eq:MVI-regret}), which is an \emph{upper bound} on the primal-dual gap of \eqref{eq:log-game}. To define our MVI, as is standard, we stack $\bx, \by$ to define $\bz$ via $\bz := (\bx, \by)$, and use the monotone operator associated with the convex-concave problem \eqref{eq:log-game}, defined by
\begin{equation}\label{eq:mon-op-def}
    \mF(\bz) = \mF(\bx, \by) := \begin{pmatrix}
        \nabla_\bx \Lcal(\bx, \by)\\ -\nabla_\by \Lcal(\bx, \by)
    \end{pmatrix} = \begin{pmatrix}
        \log(\mA \by)\\ - \mA^\top (\bx \oslash \mA \by)
    \end{pmatrix}.
\end{equation}
For technical reasons (namely, to prevent entries of $\bx \oslash \mA \by$ in the definition of $\mF$ from blowing up), we consider MVI associated with \eqref{eq:mon-op-def} on $\Delta_m \times \Deltatr_n,$ where $\Deltatr_n$ is the \emph{truncated} simplex that ensures all coordinates are sufficiently positive; namely:
\begin{equation}\label{eq:truncated simplex}
    \Deltatr_n := \big\{\by \in \R_+^n: y_i \geq \epstr \text{ for } i \in [n] \text{ and } \bone^\top\by = 1\big\}.
\end{equation}
Recall that $\epstr=\epsilon/(8n)$.

A key result is summarized in the following lemma.
\begin{lemma}\label{lemma:PCG-to-Minty-gap}
    Consider \eqref{eq:packing-covering-game} and the associated monotone operator $\mF$ defined by \eqref{eq:mon-op-def}. Given $K \geq 0,$ fix any sequence $\{\bz_k\}_{k=0}^K,$ $\bz_k = (\bx_k, \by_k) \in \Delta_m \times \Deltatr_n$ and define   $\bxh_K := \frac{\sum_{k=0}^K \bx_k \oslash (\mA \by_k)}{\bone^\top (\sum_{k=0}^K \bx_k \oslash (\mA \by_k))}$ and $\byh_K = \frac{1}{K+1}\sum_{k=0}^K \by_k$. Then:
    \begin{equation}\notag
        \log\Big(\frac{\max_{j \in [n]}(\mA^\top \bxh_K)_j}{\min_{i \in [m]}(\mA\byh_K)_i}\Big) \leq \sup_{\bw \in \Delta_m \times \Deltatr_n} \frac{1}{K+1}\sum_{k=0}^K \mF(\bz_k)^\top(\bz_k - \bw) - \log(1 - n\epstr). 
    \end{equation}
\end{lemma}
\begin{proof}
    Write $\log\big(\frac{\max_{j \in [n]}(\mA^\top \bxh_K)_j}{\min_{i \in [m]}(\mA\byh_K)_i}\big) = \log(\max_{j \in [n]}(\mA^\top \bxh_K)_j) - \log(\min_{i \in [m]}(\mA\byh_K)_i)$. We carry out the proof by bounding $\log(\max_{j \in [n]}(\mA^\top \bxh_K)_j)$ above and bounding $\log(\min_{i \in [m]}(\mA\byh_K)_i)$ below.

    Define $\bq_k := \bx_k \oslash (\mA \by_k)$ and $\bqh_K := \frac{1}{K+1}\sum_{k=0}^K \bq_k$ so that $\bxh_K = \bqh_K/(\bone^\top \bqh_K).$ We have: 
    \begin{align}
        \log(\max_{j \in [n]}(\mA^\top \bxh_K)_j) &= \log(\max_{\bv \in \Delta_n} \bv^\top \mA^\top \bxh_K)\notag \\
        &\leq \log(\max_{\bv \in \Deltatr_n} \bv^\top \mA^\top \bxh_K/(1-n\epstr))\notag\\
        &= \max_{\bv \in \Deltatr_n}\log(\bv^\top \mA^\top \bqh_K) - \log(\bone^\top \bqh_K) - \log(1-n\epstr), \label{eq:log-max-1}
    \end{align}
    where in the inequality, we have used that $\max_{\bv \in \Deltatr_n} \bv^\top \mA^\top \bxh_K \geq (1-n \epstr) \max_{j \in [n]} (\mA^\top \bxh_K)_j.$ 
    We bound the first two terms in \eqref{eq:log-max-1} separately. First, observe that $\by^\top \mA^\top (\bx \oslash (\mA \by)) = 1$ and so $\by_k^\top\mA^\top\bq_k=1$.  Recall that $\log(t) \leq t -1$, for $t  > 0$. We thus have, for any $\bv \in \Delta_n$: 
    \begin{align}
        \log(\bv^\top \mA^\top \bqh_K) &\leq \bv^\top \mA^\top \bqh_K -1 \notag\\
        &= \frac{1}{K+1} \sum_{k=0}^K (\bv^\top \mA^\top \bq_k -1)\notag\\
        &= \frac{1}{K+1} \sum_{k=0}^K (\bv - \by_k)^\top \mA^\top (\bx_k \oslash (\mA \by_k)).\label{eq:log-max-2}
    \end{align}
    On the other hand, recalling that $\bx \in \Delta_m$ and using Jensen's inequality, we have
    \begin{align}\label{eq:log-max-3}
        \log(\bone^\top \bqh_K) &= \log\Big(\frac{1}{K+1}\sum_{k=0}^K\sum_{j=1}^m x_{k,j} \frac{1}{(\mA\by_k)_j}\Big) \geq \frac{1}{K+1}\sum_{k=0}^K\sum_{j=1}^m x_{k, j} \log(1/(\mA\by_k)_j)\notag\\
        &= - \frac{1}{K+1}\sum_{k=0}^K \bx_k^\top \log(\mA\by_k). 
    \end{align}
    Combining \eqref{eq:log-max-1}--\eqref{eq:log-max-3}, we thus have:
    \begin{equation}\label{eq:log-max-final}
         \log(\max_{j \in [n]}(\mA^\top \bxh_K)_j) \leq \max_{\bv \in \Deltatr_n}\frac{1}{K+1}\sum_{k=0}^K\big[(\bv - \by_k)^\top \mA^\top (\bx_k \oslash (\mA \by_k)) + \bx_k^\top \log(\mA\by_k)\big] - \log(1-n\epstr). 
    \end{equation}
    Another two applications of Jensen's inequality similarly give
    \begin{equation}\label{eq:log-min}
        \log(\min_{i \in [m]}(\mA\byh_K)_i) = \min_{\bu \in \Delta_m}  \log(\bu^\top\mA \byh_K) \geq  \min_{\bu \in \Delta_m}  \bu^\top\log(\mA \byh_K) \geq \min_{\bu \in \Delta_m}  \frac{1}{K+1}\sum_{k=0}^K \bu^\top\log(\mA \by_k). 
    \end{equation}
    It remains to subtract \eqref{eq:log-min} from \eqref{eq:log-max-final}, use the definitions of $\mF$ and $\bz_k,$ group the like terms, and set $\bw^\top:= (\bu, \bv)$. 
\end{proof}
The rest of the upper bound proof for \eqref{eq:packing-covering-game} now boils down to arguing that we can efficiently bound the MVI regret of $\mF$ on $\Zcal := \Delta_m \times \Deltatr_n$; namely, that we can get $\sup_{\bw \in \Delta_m \times \Deltatr_n} \frac{1}{K+1}\sum_{k=0}^K \mF(\bz_k)^\top(\bz_k - \bw) \leq \epsilon$ in $\propto(1/\epsilon)$ iterations, with each iteration taking a constant number of matrix-vector products and near-linear in $\nnz(\mA)$ number of arithmetic operations.

\subsection{Solving the Constructed Variational Inequality}\label{sec:generic-VI}

It is immediate from \Cref{lemma:PCG-to-Minty-gap} that to obtain a $(1 + O(\epsilon))$ solution to \eqref{eq:packing-covering-game}, it suffices to use that $\epstr = O(\eps/n)$ and approximately solve the MVI associated with $\mF$ to get $\sup_{\bw \in \Zcal} \frac{1}{K+1}\sum_{k=0}^K \mF(\bz_k)^\top(\bz_k - \bw) = O(\epsilon),$ where $\Zcal := \Delta_m \times \Deltatr_n.$ In this subsection, we quickly argue that it is possible to get such a guarantee using $O(L\sup_{\bz \in \Zcal} D_\psi(\bz, \bz_0)/\epsilon)$ iterations, each using a relaxed mirror-prox update with a dgf $\psi$ with respect to which $\mF$ is $L$-relaxed relatively Lipschitz (see \Cref{def:relaxed-rel-Lip}), initialized at some $\bz_0 \in \Zcal$. The subsequent subsections will then introduce such a dgf, argue that there is an explicit (and thus efficiently implementable) initialization for which $\sup_{\bz \in \Zcal} D_\psi(\bz, \bz_0) = O(\log(mn/\epsilon))$, and argue that the required prox steps can be implemented in near-linear time in $\nnz(\mA),$ using a constant number of matrix-vector products. 

To obtain our result, we rely on the guarantee for a relaxed mirror-prox method with asymmetrical step sizes (see \cite{She17,diakonikolas2021efficient} for related variants), which was analyzed under relaxed relative Lipschitzness and exact prox updates in \cite{JT23}. We briefly summarize this result and provide an elementary extension to inexact prox updates, which are needed in our algorithm. 

\begin{fact}[Proposition 3 in \cite{JT23}]\label{fact:op-extrapolation-convergence}
    Let $\Zcal$ be a closed convex set and let $\mF$ be an operator that is $L$-relaxed relative Lipschitz with respect to some function $\psi,$ where $L > 0.$ Let $\bz_0 \in \Zcal$ and for $k \geq 0,$ let 
    \begin{equation}\label{eq:relaxed-MP}
    \begin{aligned}       
        \bw_{k} &\approx_{\epsilon_k} \argmin{\bz \in \Zcal} \Big\{ \frac{1}{L}\mF(\bz_k)^\top \bz + D_\psi(\bz, \bz_k) \Big\},\\
        \bz_{k+1} &\approx_{\epsilon_k} \argmin{\bz \in \Zcal} \Big\{ \frac{1}{2L}\mF(\bw_k)^\top \bz + D_\psi(\bz, \bz_k) \Big\},
    \end{aligned}
    \end{equation}
    where we recall the $\approx_{\epsilon_k}$ notation from \eqref{eq:approx-3PI}. 
    Then, for any $\bz \in \Zcal$ and any $K \geq 0:$
    \begin{equation}\notag
        \frac{1}{K+1}\sum_{k=0}^K \mF(\bw_k)^\top(\bw_k - \bz) \leq \frac{2L}{K+1} D_\psi(\bz, \bz_0) + \frac{3L}{K+1}\sum_{k=0}^K \epsilon_k. 
    \end{equation}
    In particular, if $\sup_{\bz \in \Zcal}D_{\psi}(\bz, \bz_0) < \infty,$ then for any $\epsilon \in (0, 1),$  $\eps_k = \eps/(12L)$, and $K = \lceil \frac{8 L \sup_{\bz \in \Zcal}D_{\psi}(\bz, \bz_0)}{\epsilon} \rceil$, we have  
    \begin{equation*}
        \sup_{\bz \in \Zcal}\frac{1}{K+1}\sum_{k=0}^K \mF(\bw_k)^\top(\bw_k - \bz) \leq \epsilon/2.
    \end{equation*}
\end{fact}
\begin{proof}[Proof Sketch]
    All that is needed to extend the proof of Proposition 3 in \cite{JT23} is to replace the standard prox inequality from \Cref{fact:3PI} (which is the inequality below with $\epsilon_k = 0$) with the approximate version
    \[
        \frac{1}{L}\mF(\bz_k)^\top (\bw_k - \bz_{k+1}) \leq D_\psi(\bz_{k+1}, \bz_k) - D_\psi(\bz_{k+1}, \bw_k) - D_\psi(\bw_k, \bz_k) + \epsilon_k
    \]
    and, similarly, 
    \[
        \frac{1}{2L}\mF(\bw_k)^\top(\bz_{k+1} - \bz) \leq D_\psi(\bz, \bz_k) - D_\psi(\bz, \bz_{k+1}) - D_\psi(\bz_{k+1}, \bz_k) + \epsilon_k.
    \]
    The rest of the proof follows from  \cite{JT23}, accounting for the
additive error
$3L\epsilon_k$ per iteration. 
\end{proof}

\subsection{Distance Generating Function}\label{sec:dgf}

To construct the needed dgf, we first start with proving that $\mF$ is $7$-relaxed relatively Lipschitz with respect to a function that mimics its structure through a similar linear-log coupling. Namely, the function is:
\begin{equation}\label{eq:aux-dgf}
    \varphi(\bx, \by) :=  \bx^\top\log(\bx) + \bx^\top\log(\bx \oslash\mA\by).
\end{equation}
It is immediate that this function is jointly convex in $\bx, \by$ (by \Cref{fact:convexity-of-KL}). Further, it satisfies the needed property for the dgf in \Cref{fact:op-extrapolation-convergence}, as shown in the following lemma, which is proved in Appendix~\ref{appx:proof-lemma:relaxed-rel-Lip-of-aux-reg}.
\begin{lemma}\label{lemma:relaxed-rel-Lip-of-aux-reg}
    Let $\mF$ be defined by \eqref{eq:mon-op-def} and $\varphi$ be defined by \eqref{eq:aux-dgf}. Then $\mF$ is $7$-relaxed relative Lipschitz with respect to $\varphi$ on $\Zcal = \Delta_m \times \Deltatr_n.$ 
\end{lemma}
While function $\varphi$ establishes the desired relaxed relative Lipschitzness condition needed for applying \Cref{fact:op-extrapolation-convergence}, it is not immediately usable for our algorithm, as it is unclear how to efficiently implement the prox steps in \eqref{eq:relaxed-MP} because of the term $- \bx^\top \log(\mA\by)$. In particular, the prox term computation requires minimizing the sum of linear functions in $\bx$ and $\by$ and $\varphi(\bx, \by).$ Minimization over $\bx$ is easy, as, treating $\by$ as a constant, it boils down to the prox function of the entropy, which is expressible in closed form on the probability simplex, and the problem can be reduced to a minimization over $\by \in \Deltatr_n$ only. However, the resulting minimization problem appears nontrivial to solve in near-linear time in $\nnz(\mA),$ without losing additional $\mathrm{poly}(1/\epsilon)$ factors to solve it approximately. 

Instead, we construct a different dgf $\psi$ inspired by similar functions appearing in the optimal transport literature and prove that $\varphi$ is $(L/7)$-relatively smooth with respect to $\psi$, with a parameter $L = \Theta(\log(n)).$ As a consequence, $\mF$ is $L$-relaxed relative Lipschitz with respect to $\psi.$ We then argue  later in this section that the prox steps required by \eqref{eq:relaxed-MP} can be computed efficiently, in near-linear time in $\nnz(\mA).$ 

The key insight in constructing the dgf is to notice that $\varphi(\bx, \by)$ can be written as the sum of the negative entropy $h(\bx)$ and the KL projection with the fixed $\bx$-marginal, the latter being an elementary operation in the optimal transport literature (see, e.g., \cite{CPSV18}). To formalize this intuition, we need to introduce the notation for the couplings of $\bx$ and $\by$ compatible with the structure of $\mA$. To do so, define 
\begin{equation}\label{eq:A-simplex}
    \Delta_\mA := \{\bpi \in \R_+^{m \times n}: A_{ij} = 0 \Rightarrow \Pi_{ij} = 0 \text{ and } \bone_m^\top \bpi \bone_n = 1\},
\end{equation}
which is the probability simplex over the non-zero entries of $\mA.$ Then $\bpi \in \Delta_\mA$ are the ``allowed'' couplings we consider next. 
Now observe that:
\begin{equation}\label{eq:varphi-OT-interpretation}
    \varphi(\bx, \by) = h(\bx) + \min_{\bpi \in \Delta_{\mA}: \bpi \bone = \bx} \Big\{\sum_{i=1}^m \sum_{j=1}^n \Pi_{ij}\log\Big(\frac{\Pi_{ij}}{A_{ij}y_j}\Big)\Big\}.
\end{equation}
 The minimizer of the above problem defines a ``transport plan'' with the specified marginal $\bx$; see, e.g., \cite{CPSV18}. Similar minimization problems appear in the definition of the classical Sinkhorn algorithm and its generalizations. This interpretation of $\varphi$ through couplings $\bpi$ motivates the search for similar dgfs whose prox mappings are easier to compute.

Observe first that $\sum_{i=1}^m \sum_{j=1}^n \Pi_{ij}\log\Big(\frac{\Pi_{ij}}{A_{ij}y_j}\Big) = \inner{\bpi, \log(\bpi\oslash \mA)} - \inner{\bpi^\top \bone_m, \log(\by)}.$ Formulation in \eqref{eq:varphi-OT-interpretation} suggests designing the dgf through its coupling formulation. In a prox computation, which minimizes the dgf plus a linear function over $\bx$ and $\by,$ eliminating \(\bx=\bpi\bone_n\) turns the outside term \(h(\bx)\) into \(h(\bpi\bone_n)\), which couples different columns of \(\bpi\). We therefore omit this term, allowing minimization over \(\bpi\) to be performed independently in each column once the column marginal \(\bpi^\top\bone_m\) is fixed. To obtain a suitable geometry, we add \(h(\by)\) and regularize the column marginal by \(\beta h(\bpi^\top\bone_m+\epstr \bone_n)\). The column-marginal penalty is constant during the columnwise minimization over \(\bpi\), so it preserves its closed form, while \(h(\by)\) provides uniform curvature for the conditional minimization over \(\by\). Together, these additions also supply sufficient curvature to control the Bregman divergence of \(\varphi\) with a width-independent factor \(O(1/\beta)\), as established in \Cref{lemma:dgf-properties} below. The shift by $\epstr \bone_n$ keeps the marginal-entropy gradient finite even when a column marginal $\bpi^\top\bone_m$ vanishes. This leads to our final dgf:
\begin{equation}\label{eq:final-dgf}
    \psi(\bx, \by) := h(\by) + \min_{\bpi \in \Delta_{\mA}: \bpi \bone_n = \bx}\Big\{\inner{\bpi, \log(\bpi\oslash \mA)} - \inner{\bpi^\top \bone_m, \log(\by)} + \beta\, h\big(\bpi^\top \bone_m + \epstr \bone_n\big)\Big\},
\end{equation}
where $\beta := \frac{1}{1 + \log(n)}$ is the regularization parameter. We first argue that the minimizer of the problem in \eqref{eq:final-dgf} can be written in an equivalent form that will get repeatedly used in our proofs.

\begin{lemma}\label{lemma:Pi_*}
    Given $\bx \in \Delta_m,\, \by \in \relint(\Delta_n),$ let 
    \begin{equation}\notag
        \bpi_* \in \argmin{\bpi \in \Delta_{\mA}: \bpi \bone_n = \bx}\Big\{\inner{\bpi, \log(\bpi\oslash \mA)} - \inner{\bpi^\top \bone_m, \log(\by)} + \beta\, h\big(\bpi^\top \bone_m + \epstr \bone_n\big)\Big\}.
    \end{equation}
    Then $\bpi_*$ is unique and the following equality characterizes it, entrywise:
    \begin{equation}\label{eq:minimizer-Pi_*}
        \Pi_{*, ij} = \frac{x_i A_{ij}y_j([\bpi_*^\top \bone_m]_j + \epstr)^{-\beta}}{\sum_{j'=1}^n A_{ij'}y_{j'}([\bpi_*^\top \bone_m]_{j'} + \epstr)^{-\beta}}.
    \end{equation}
\end{lemma}
\begin{proof}
    Uniqueness follows by strong convexity: the first term in the minimum is the entropy of $\bpi$ (plus a linear term in $\bpi$), the second is linear in $\bpi,$ and the last one is convex. \eqref{eq:minimizer-Pi_*} follows from the first-order optimality conditions of the minimization problem that defines $\bpi_*.$
\end{proof}

We now argue the properties of $\psi$ that are needed for applying \Cref{fact:op-extrapolation-convergence}.

\begin{lemma}\label{lemma:dgf-properties}
    Let function $\psi: \Delta_m \times \relint(\Delta_n) \to \R$ be defined by \eqref{eq:final-dgf}. Then $\psi$ satisfies all of the following properties:
    \begin{enumerate}[label=(\roman*)]
        \item $\psi$ is convex.
        \item $\psi$ is differentiable for $\bx \in \relint(\Delta_m),\, \by \in \relint(\Delta_n).$ Moreover,
        \begin{equation}\notag
            \begin{aligned}
                (\nabla_\bx \psi(\bx, \by))_i &= 1 + \beta + \log\Big(\frac{\Pi_{*, ij}}{A_{ij}y_j}\Big) + \beta\log([\bpi_*^\top \bone_m]_j + \epstr), \quad i \in [m],~j\in[n]:A_{ij}\neq 0;\\
                \nabla_\by \psi(\bx, \by) &= \log(\by) + \bone_n - (\bpi_*^\top \bone_m) \oslash \by,
            \end{aligned}
        \end{equation}
        where $\bpi_*$ is the solution to the minimum in \eqref{eq:final-dgf}, given $\bx \in \relint(\Delta_m),\, \by \in \relint(\Delta_n)$. Additionally, for $\bu \in \Delta_m, \bv \in \relint(\Delta_n),$ and $\bpi_*'$ their corresponding minimizer in the definition of $\psi(\bu, \bv),$ we have
        \begin{equation}\label{eq:D-psi}
            D_\psi((\bu, \bv), (\bx, \by)) = D_h(\bv, \by) + D_h(\bpi_*', \bpi_*\diag(\bv\oslash \by)) + \beta D_h\big((\bpi_*')^\top \bone_m + \epstr \bone_n,\, \bpi_*^\top \bone_m + \epstr \bone_n\big). 
        \end{equation}
        \item For any $(\bu, \bv) \in \Zcal$ and $(\bx, \by) \in \Zcal$ with $\bx>\bzero$, $D_\psi((\bu, \bv), (\bx, \by)) \geq \frac{\beta}{3 + 4 e}D_\varphi((\bu, \bv), (\bx, \by)).$ As a result, $\mF$ defined by \eqref{eq:mon-op-def} is $\frac{7(3 + 4e)}{\beta}$-relaxed relatively Lipschitz with respect to $\psi.$
        \item For $(\bx, \by) \in \Zcal,$ $\bx > 0$ and $(\bu, \bv) \in \Zcal,$ $D_\psi((\bu, \bv), (\bx, \by)) \geq \frac{\beta}{2}D_h(\bu, \bx) + \frac{1}{2}D_h(\bv, \by).$ 
    \end{enumerate}
\end{lemma}
\begin{proof}
Let $\Phi(\bpi, \by) := \inner{\bpi, \log(\bpi\oslash \mA)} - \inner{\bpi^\top \bone_m, \log(\by)} + \beta\, h\big(\bpi^\top \bone_m + \epstr \bone_n\big)$ be the function under the minimum in \eqref{eq:final-dgf}, so that we can write $\psi(\bx, \by) = h(\by) + \min_{\bpi \in \Delta_{\mA}: \bpi \bone_n = \bx} \Phi(\bpi, \by)$. Observe that $\Phi(\bpi, \by) = \inner{\bpi, \log(\bpi\oslash(\mA \diag(\by)))} +\beta h\big(\bpi^\top \bone_m + \epstr \bone_n\big)$. The first term is the relative entropy between $\bpi$ and a linear transformation of $\by$, which is jointly convex in $\bpi$ and $\by$ by \Cref{fact:convexity-of-KL}. The second term, $h\big(\bpi^\top \bone_m + \epstr \bone_n\big),$ is the shifted entropy of $\bpi^\top \bone_m,$ which is likewise convex in $\bpi.$ Thus, $\Phi$ is jointly convex in $\bpi, \by$, and, by \Cref{fact:inf-convexity}, $\min_{\bpi \in \Delta_{\mA}: \bpi \bone_n = \bx} \Phi(\bpi, \by)$ is convex in $\bx, \by$ (convexity in $\bx$ follows also from \Cref{fact:inf-convexity}, by defining $\Phi$ to be a function of $\bx$ set to $\infty$ when $\bpi\bone_n = \bx$ is violated; see \cite[Example 3.17]{vandenberghe2004convex}). We thus conclude that $\psi$ is convex, completing the proof of Part $(i)$.

For Part $(ii),$ we apply \Cref{fact:danskin} to differentiate $\min_{\bpi \in \Delta_{\mA}: \bpi \bone_n = \bx}\Phi(\bpi, \by)$. We first consider the differentiation with respect to $\by.$ To apply \Cref{fact:danskin}, we note that the differentiability of $\Phi$ with respect to $\by \in \relint(\Delta_n)$ and continuity of
$\nabla_{\by}\Phi(\bpi,\by)$ on $\relint(\Delta_\mA)$ are immediate, while the uniqueness of $\argmin{\bpi \in \Delta_{\mA}: \bpi \bone_n = \bx}\Phi(\bpi, \by)$ follows by the strong convexity of $h$. Given $\bx \in \Delta_m, \by \in \relint(\Delta_n)$, let $\bpi_* := \argmin{\bpi \in \Delta_{\mA}: \bpi \bone_n = \bx}\Phi(\bpi, \by)$; it thus follows that:
\begin{equation}\label{eq:grad-psi-wrt-y}
    \nabla_\by\psi(\bx, \by) = \log(\by) + \bone_n - (\bpi_*^\top \bone_m) \oslash \by. 
\end{equation}
Now consider the differentiation with respect to $\bx.$ Because $\bx$ does not appear as an argument of $\Phi,$ but in the minimization constraint, we first need to reformulate $\min_{\bpi \in \Delta_{\mA}: \bpi \bone_n = \bx}\Phi(\bpi, \by)$ to apply \Cref{fact:danskin} more directly. In particular, consider the following box-constrained marginal set:
\begin{equation}\notag
    \Scal_{\mA} := \big\{\bgamma \in \R_+^{m \times n}: A_{ij} = 0 \Rightarrow \Gamma_{ij} = 0 \text{ and } \bgamma \bone_n = \bone_m\big\}. 
\end{equation}
Observe that the matrices $\bpi\in\Delta_{\mA}$ satisfying $\bpi\bone_n=\bx$
are exactly those of the form
$\bpi=\diag(\bx)\bgamma$ with $\bgamma\in \Scal_{\mA}$.
As a result,
\begin{equation}\notag
    \min_{\bpi \in \Delta_{\mA}: \bpi \bone_n = \bx}\Phi(\bpi, \by) \quad \Leftrightarrow \quad \min_{\bgamma \in \Scal_{\mA}}\Phi(\diag(\bx) \bgamma, \by). 
\end{equation}
We can now use \Cref{fact:danskin} to differentiate $\min_{\bgamma \in \Scal_{\mA}}\Phi(\diag(\bx)\bgamma, \by)$ with respect to $\bx$ (the assumptions required by the fact are immediate). Given $\bx \in \relint(\Delta_m)$, $\by \in \relint(\Delta_n),$ let $\bgamma_* := \argmin{\bgamma \in \Scal_{\mA}}\Phi(\diag(\bx)\bgamma, \by)$ and recall $\bpi_* = \diag(\bx) \bgamma_*.$ Since $\nabla_{\bx} \psi(\bx, \by) = \nabla_{\bx}\Big(\min_{\bgamma \in \Scal_{\mA}}\Phi(\diag(\bx) \bgamma, \by)\Big),$ we get
\begin{equation}\notag
    (\nabla_\bx\psi(\bx, \by))_i = \sum_{j=1}^n \Gamma_{*,ij}\Big(\log\Big(\frac{x_i\Gamma_{*,ij}}{A_{ij} y_{j}}\Big) + 1 + \beta \log([\bgamma_*^\top \bx]_j + \epstr) + \beta\Big), \quad i \in [m].
\end{equation}
The expression inside the big parentheses is $\frac{1}{x_i}\frac{\partial}{\partial \Gamma_{ij}}\Phi(\diag(\bx)\bgamma, \by)\mid_{\bgamma = \bgamma_*}.$ Note that at $\bgamma = \bgamma_*,$ it must be $\frac{\partial}{\partial \Gamma_{ij}}\Phi(\diag(\bx) \bgamma, \by)\mid_{\bgamma = \bgamma_*} = \frac{\partial}{\partial \Gamma_{ij'}}\Phi(\diag(\bx) \bgamma, \by)\mid_{\bgamma = \bgamma_*}$ for all $j, j'$ such that $A_{ij} \neq 0, A_{ij'} \neq 0$. The reason is that, at optimality, it must be $\bgamma_{*,ij} > \bzero$
whenever $A_{ij}\neq 0,$
so if there were $j, j'$ with  $\frac{\partial}{\partial \Gamma_{ij}}\Phi(\diag(\bx)\bgamma, \by)\mid_{\bgamma = \bgamma_*} \neq \frac{\partial}{\partial \Gamma_{ij'}}\Phi(\diag(\bx)\bgamma, \by)\mid_{\bgamma = \bgamma_*}$, there would exist a small redistribution of mass between entries $ij$ and $ij'$ that would keep $\bgamma \in \Scal_{\mA},$ but decrease the value of $\Phi,$ leading to a contradiction. 
Now, since $\sum_{j=1}^n \Gamma_{*, ij} = 1$, $\Pi_{*, ij} = x_i \Gamma_{*, ij}$, and $\bgamma_*^\top \bx = \bpi_*^\top \bone_m$,
we equivalently have:
\begin{equation}\label{eq:grad-psi-wrt-x}
    (\nabla_\bx\psi(\bx, \by))_i = 1 + \beta + \log\Big(\frac{\Pi_{*, ij}}{A_{ij}y_j}\Big) + \beta\log([\bpi_*^\top \bone_m]_j + \epstr), \quad i \in [m],~j\in[n]:A_{ij}\neq 0.
\end{equation}
Moreover, if for $\bx' \in \Delta_m, \by' \in \relint(\Delta_n)$, we denote $\bpi_*' := \argmin{\bpi \in \Delta_{\mA}: \bpi \bone_n = \bx'}\Phi(\bpi, \by')$, then
(since $\Phi(\bpi_*, \by)$ depends on $\bx$ only through $\bpi_* \bone_n = \bx$), we have
\begin{equation}\notag
    D_\psi\big((\bx', \by'), (\bx, \by)\big) = D_h(\by', \by) + D_h(\bpi_*', \bpi_*\diag(\by'\oslash \by)) + \beta D_h\big((\bpi_*')^\top \bone_m + \epstr \bone_n,\, \bpi_*^\top \bone_m + \epstr \bone_n\big). 
\end{equation}

We now argue Part $(iv).$ First, in Part $(i)$, we have already argued that $\min_{\bpi \in \Delta_{\mA}: \bpi \bone_n = \bx} \Phi(\bpi, \by)$ is convex in $\bx, \by$. Since $\psi(\bx, \by) = h(\by) + \min_{\bpi \in \Delta_{\mA}: \bpi \bone_n = \bx} \Phi(\bpi, \by)$, it immediately follows that $D_{\psi}((\bu, \bv), (\bx, \by)) \geq D_h(\bv, \by).$ Thus, to complete the proof of Part $(iv)$, it suffices to argue that also $D_{\psi}((\bu, \bv), (\bx, \by)) \geq \beta D_h(\bu, \bx).$ Equivalently, it suffices to argue that $\psi(\bx, \by) - \beta h(\bx)$ is convex on the defined domain. Fix any $\bpi \in \Delta_\mA$ with $\bpi \bone_n = \bx.$ Then, using elementary algebra, plugging in $\bx = \bpi \bone_n,$ and grouping terms, we get:
\begin{equation}\label{eq:psi-rel-convex-h(x)}
    \begin{aligned}
        h(\by) + \Phi(\bpi, \by) - \beta h(\bx) =\; & (1-\beta) \Big(h(\bpi) - \inner{\bpi^\top \bone_m, \log(\by)} + h(\by)\Big) \\
        &+ \beta \Big(h\, (\bpi^\top \bone_m + \epstr \bone_n) - \inner{\bpi^\top \bone_m, \log(\by)} + h(\by)\Big)\\
        &+ \beta(h(\bpi) - h(\bpi \bone_n)) - \inner{\bpi, \log(\mA)}.
    \end{aligned}
\end{equation}
Minimization of the left-hand side of \eqref{eq:psi-rel-convex-h(x)} over $\bpi \in \Delta_\mA: \bpi \bone_n = \bx$ gives $\psi(\bx, \by) - \beta h(\bx)$. Thus, using \Cref{fact:inf-convexity}, it suffices to argue that the right-hand side of \eqref{eq:psi-rel-convex-h(x)} is jointly convex in $\bpi, \by$. We do so by arguing convexity of each line on the right-hand side of \eqref{eq:psi-rel-convex-h(x)}. For the first line, the convexity follows immediately, since $h(\by)$ is convex and $h(\bpi) - \inner{\bpi^\top \bone_m, \log(\by)}=\langle\bpi,\log(\bpi\oslash(\mathbf{1}_m\by^\top))\rangle$ is the relative entropy function, thus convex by \Cref{fact:convexity-of-KL}. The expression in the second line, when written in its coordinate-wise sum form, is coordinate-wise equivalent to the function from \Cref{fact:shifted-entropy-expression-convexity}, which is convex by the same fact for $\by \in \Deltatr_n$. Finally, for the third line, $- \inner{\bpi, \log(\mA)}$ is linear in $\bpi$ (thus convex) and the remaining conditional entropy term $h(\bpi) - h(\bpi \bone_n)$ is convex, as, after elementary algebra, it can be written as $h(\bpi) - h(\bpi \bone_n) = \sum_{i=1}^m\sum_{j=1}^n \Pi_{ij}\log(\frac{\Pi_{ij}}{(\bpi\bone_n)_i})$, which is convex by \Cref{fact:convexity-of-KL}. 

It remains to prove Part $(iii).$ Fix any $(\bu, \bv) \in \Zcal$, $(\bx, \by) \in \Zcal$, and $t \in [0, 1]$, and define:
\begin{equation}
    \begin{aligned}
        \varphi_1(t) := \varphi(t(\bu, \bv) + (1-t)(\bx, \by)),\\
        \psi_1(t) := \psi(t(\bu, \bv) + (1-t)(\bx, \by)).
    \end{aligned}
\end{equation}
Observe that
\[
D_{\varphi}((\bu, \bv), (\bx, \by)) = \varphi_1(1) - \varphi_1(0) - \varphi_1'(0) = \int_0^1 (1-t)\varphi_1''(t)\mathrm{d}t,
\]
where the last equality holds by Taylor's Theorem, 
and we have an equivalent expression for $D_\psi.$ Thus, to prove that $D_\varphi((\bu, \bv), (\bx, \by)) \leq M D_{\psi}((\bu, \bv), (\bx, \by))$ for some $M > 0,$ it suffices to prove that $\varphi_1''(t) \leq M \psi_1''(t)$ for all $t \in (0, 1)$. To differentiate $\varphi_1$ and $\psi_1$, let us define $(\bx_t, \by_t) := t(\bu, \bv) + (1-t)(\bx, \by).$ We further use the dot notation for differentiating with respect to $t$; namely, $\dot{\bx}_t = \frac{\dd}{\dd t}\bx_t$ and $\ddot{\bx}_t = \frac{\dd^2}{\dd t^2}\bx_t$. Observe that, by definition: $\dot{\bx}_t = \bu - \bx,$ $\ddot{\bx}_t = \bzero,$ and, similarly, $\dot{\by}_t = \bv - \by,$ $\ddot{\by}_t = \bzero.$ It thus follows that:
\begin{equation}\label{eq:2nd-derivative}
\begin{aligned}
    \varphi_1''(t) &= \Big\langle{\frac{\partial}{\partial t}\nabla_\bx \varphi(\bx_t, \by_t), \dot{\bx}_t}\Big\rangle + \Big\langle{\frac{\partial}{\partial t}\nabla_{\mA\by} \varphi(\bx_t, \by_t), \dot{\mA\by}_t}\Big\rangle,\\
    \psi_1''(t) &= \Big\langle{\frac{\partial}{\partial t}\nabla_\bx \psi(\bx_t, \by_t), \dot{\bx}_t}\Big\rangle + \Big\langle{\frac{\partial}{\partial t}\nabla_\by \psi(\bx_t, \by_t), \dot{\by}_t}\Big\rangle.
\end{aligned}
\end{equation}
A direct calculation now leads to:
\begin{align}
    \varphi_1''(t) &= \sum_{i=1}^m\Big(2 \frac{\dot{x}_{t, i}^2}{x_{t, i}} - 2 \frac{\dot{x}_{t_i}(\dot{\mA\by_t})_i}{(\mA\by_t)_i} + \frac{x_{t, i}(\dot{\mA\by_t})_i^2}{(\mA\by_t)_i^2}\Big)\notag\\
    &= \sum_{i'=1}^m \frac{\dot{x}_{t, i'}^2}{x_{t, i'}} + \sum_{i=1}^m x_{t, i}\Big(\frac{\dot{x}_{t,i}}{x_{t, i}} - \frac{(\dot{\mA\by_t})_i}{(\mA\by_t)_i}\Big)^2. \label{eq:varphi_1''}
\end{align}
To derive $\psi_1''(t)$, we use \eqref{eq:grad-psi-wrt-y}, \eqref{eq:grad-psi-wrt-x}, and define $\bpi_t := \argmin{\bpi \in \Delta_{\mA}: \bpi \bone_n = \bx_t}\Phi(\bpi, \by_t).$ Observe that $\dot{(\bpi_t \bone_n)} = \dot{\bx}_t,$ but $\dot{(\bpi_t^\top \bone_m)} = \dot{\bpi}_t^\top \bone_m$ does not have such an explicit formula in terms of $\bx_t$ and/or $\by_t$. First, differentiating \eqref{eq:grad-psi-wrt-y} with respect to $t,$ we get:
\begin{equation}\label{eq:psi-partial-t-y}
    \Big\langle{\frac{\partial}{\partial t}\nabla_\by \psi(\bx_t, \by_t), \dot{\by}_t}\Big\rangle = \sum_{j=1}^n \Big(\frac{\dot{y}_{t, j}^2}{y_{t, j}} - \frac{(\dot{\bpi}_t^\top \bone_m)_j\dot{y}_{t, j}}{y_{t, j}} + \frac{(\bpi_t^\top \bone_m)_j\dot{y}_{t, j}^2}{y_{t, j}^2}\Big).
\end{equation}
Second, differentiating \eqref{eq:grad-psi-wrt-x} with respect to $t$ and recalling that $\sum_{j=1}^n \dot{\Pi}_{t, ij} =  \sum_{j\in [n]:A_{ij} \neq 0} \dot{\Pi}_{t, ij} = \dot{x}_{t, i}$ we get:
\begin{equation}\label{eq:psi-partial-t-x}
    \Big\langle{\frac{\partial}{\partial t}\nabla_\bx \psi(\bx_t, \by_t), \dot{\bx}_t}\Big\rangle = \sum_{\substack{(i,j)\in[m]\times[n]:\\A_{ij}\neq 0}}\frac{\dot{\Pi}_{t, ij}^2}{\Pi_{t, ij}} - \sum_{j=1}^{n}\frac{[\dot{\bpi}_{t}^\top \bone_m]_j \dot{y}_{t, j}}{y_{t, j}} + \beta \sum_{j=1}^{n}\frac{[\dot{\bpi}_{t}^\top \bone_m]_j^2}{[{\bpi}_{t}^\top \bone_m]_j + \epstr}. 
\end{equation}
Summing \eqref{eq:psi-partial-t-y} and \eqref{eq:psi-partial-t-x} and simplifying, we get:
\begin{equation}\label{eq:psi_1''}
    \psi_1''(t) = \sum_{ij \in [m]\times[n]: A_{ij}\neq 0}\frac{(\dot{\Pi}_{t, ij} - \Pi_{t, ij} \dot{y}_{t, j}/y_{t, j})^2}{\Pi_{t, ij}} + \beta \sum_{j'=1}^n \frac{[\dot{\bpi}_{t}^\top \bone_m]_{j'}^2}{[{\bpi}_{t}^\top \bone_m]_{j'} + \epstr} + \sum_{j''=1}^n \frac{\dot{y}_{t, j''}^2}{y_{t, j''}}. 
\end{equation}
It now remains to compare $\varphi_1''(t)$ and $\psi_1''(t).$ We first observe that $\sum_{i'=1}^m \frac{\dot{x}_{t, i'}^2}{x_{t, i'}} = h_1''(t),$ where $h_1(t) = h(x_t)$. We have already argued in Part $(iv)$ that $\psi(\bx, \by) - \beta h(\bx)$ is convex, thus it immediately follows that
\begin{equation}\label{eq:psi''-vs-h''}
    \psi_1''(t) \geq \beta h_1''(t). 
\end{equation}
On the other hand, we can bound $\varphi_1''(t)$ above as the sum of a constant times $h_1''(t)$ plus a term involving $\by$ that we can hope to bound separately. In particular, using Young's inequality:
\begin{equation}\label{eq:varphi''-ub-1}
    \varphi_1''(t) \leq 3 h_1''(t) + 2\sum_{i=1}^m x_{t, i}\frac{(\mA \dot{\by}_t)_i^2}{(\mA \by_t)_i^2}. 
\end{equation}
We simplify \eqref{eq:varphi''-ub-1} further by writing $\frac{(\mA \dot{\by}_t)_i}{(\mA \by_t)_i} = \sum_{j=1}^n \frac{A_{ij}y_{t, j}}{(\mA \by_t)_i}\frac{\dot{y}_{t, j}}{y_{t,j}}$ and using Jensen's inequality to get:
\begin{equation}\label{eq:varphi''-ub-2}
    \varphi_1''(t) \leq 3 h_1''(t) + 2\sum_{j=1}^n\Big(\sum_{i=1}^m x_{t, i}\frac{A_{ij}y_{t, j}}{(\mA \by_t)_i}\Big) \frac{\dot{y}_{t, j}^2}{y_{t, j}^2}. 
\end{equation}
We have the following claim (with the proof deferred to Appendix~\ref{appx:proof-claim:xy-vs-pi-relationship}):
\begin{claim}\label{claim:xy-vs-pi-relationship}
    $\sum_{i=1}^m x_{t, i}\frac{A_{ij}y_{t, j}}{(\mA \by_t)_i} \leq e([\bpi_t^\top \bone_m]_j + \epstr).$
\end{claim}
As a result, \eqref{eq:varphi''-ub-2} further simplifies to:
\begin{equation}\label{eq:varphi''-ub-3}
    \varphi_1''(t) \leq 3 h_1''(t) + 2 e\sum_{j=1}^n \big([\bpi_t^\top \bone_m]_j + \epstr\big) \frac{\dot{y}_{t, j}^2}{y_{t, j}^2}. 
\end{equation}
To complete the proof, we use the following technical claim (proved in Appendix~\ref{appx:proof-claim:psi_1''-ub-y}).
\begin{claim}\label{claim:psi_1''-ub-y}
   $\psi_1''(t) \geq \frac{\beta}{1+\beta} \sum_{j=1}^n \big([\bpi_t^\top \bone_m]_j + \epstr\big) \frac{\dot{y}_{t, j}^2}{y_{t, j}^2}$. 
\end{claim}
 Using \Cref{claim:psi_1''-ub-y}, observe that by \eqref{eq:psi''-vs-h''} and \eqref{eq:varphi''-ub-3}, we now have
\begin{equation}\notag
    \phi_1''(t) \leq \frac{3}{\beta}\psi_1''(t) + \frac{2e(1+\beta)}{\beta} \psi_1''(t) \leq \frac{3 + 4e}{\beta}\psi_1''(t),
\end{equation}
which is the relationship claimed in Part $(iii).$
\end{proof}

\subsection{Bounded-Range Initialization}\label{sec:initialization}

In addition to having a suitable dgf $\psi$ constructed in the previous subsection, for \Cref{fact:op-extrapolation-convergence} to provide a bound that ultimately leads to near-linear work, both $L$ and $\sup_{\bz \in \Zcal}D_\psi(\bz, \bz_0)$ need to be at most poly-logarithmic in $m, n , 1/\epsilon.$ We have already shown (in \Cref{lemma:dgf-properties}) that $L = O(1/\beta) = O(\log(n)),$ so to get the desired iteration complexity, we also need to argue that we can choose $\bz_0 \in \Zcal$ such that $\sup_{\bz \in \Zcal}D_\psi(\bz, \bz_0) = \mathrm{polylog}(m, n, 1/\epsilon).$ We establish such a bound in this subsection.

We first note that what appears like the most obvious initialization (setting $\bx_0 = (1/m)\bone_m,$ $\by_0 = (1/n) \bone_n$) is insufficient for our purposes. In particular, it is possible to show that in the worst case, such an initialization can lead to $\sup_{\bz \in \Zcal}D_\psi(\bz, \bz_0) = \Omega(n).$ Instead, we use the characterization of the transport plan (minimizing $\bpi$) from \Cref{lemma:Pi_*} to construct $\bz_0 = (\bx_0, \by_0)$ whose corresponding transport plan $\bpi_0$ is given by \eqref{eq:minimizer-Pi_*}. This then gives an explicit formula for $\bpi_0$ that allows us to suitably bound $\sup_{\bz \in \Zcal}D_\psi(\bz, \bz_0).$ 

\begin{lemma}[Initialization]\label{lemma:initialization}
    Let $\mAh := \mA \oslash (\mA \bone_n \bone_n^\top)$ be the row-normalized (with rows summing to one) transformation of matrix $\mA.$ Let $\bpi_0 \in \Delta_{\mA}$ be defined via:
    \begin{equation}\label{eq:initial-coupling}
        \Pi_{0, ij} \propto \Ah_{ij}\Big(\frac{(\mAh^\top\bone_m)_j^{1/\beta}}{\sum_{j'=1}^n (\mAh^\top\bone_m)_{j'}^{1/\beta}} + 2\epstr\Big)^{1-\beta},
    \end{equation}
    where the expression above omits normalization by the sum of the terms on the right-hand side, so that $\bone_m^\top \bpi_0 \bone_n = 1.$

    Let $\bz_0 = (\bx_0, \by_0)$, where $\bx_0, \by_0$ are defined by
    \begin{equation}\label{eq:initial-vectors}
        \bx_0 = \bpi_0 \bone_n, \quad \by_0 \propto {\bpi_0^\top \bone_m}\oslash{\mAh^\top \bone_m} \odot \big(\bpi_0^\top \bone_m + \epstr \bone_n\big)^{\beta},
    \end{equation}
    where, again, the $\propto$ notation omits the normalization that ensures $\bone_n^\top \by_0 = 1.$ 
    Then $\bpi_0$ satisfies \eqref{eq:minimizer-Pi_*}, $\bz_0 \in \Zcal$, and 
    \[
        \sup_{\bz \in \Zcal}D_\psi(\bz, \bz_0) = O\big(\log({m}/{\epstr})\big) = O\big(\log(mn/\epsilon)\big). 
    \]
    Moreover, $\bz_0$ can be computed with $O(1)$ matrix-vector multiplications and $O(\nnz(\mA))$ total  arithmetic operations. 
\end{lemma}
\begin{proof}
    Let $\bpi_0,$ $\bx_0,$ and $\by_0$ be defined as in the lemma statement. 
    
    Let us first argue that $\bz_0$ can be computed with $O(1)$ matrix-vector multiplications, with $O(\nnz(\mA))$ arithmetic operations. To see this, notice that $\mAh = \diag(\bone_m\oslash (\mA\bone_n)) \mA$. Since $\mA\bone_n$ takes one matrix-vector multiplication (which takes $O(\nnz(\mA))$ time) and writing $\diag(\bone_m\oslash (\mA\bone_n))$ explicitly takes $O(\nnz(\mA))$ time, we can simulate matrix-vector products involving $\mAh$ using matrix-vector products with $\mA$ and rescaling by $\diag(\bone_m\oslash (\mA\bone_n))$. Computing $\bx_0$ and $\by_0$ from $\bpi_0$ further requires computing $\bone_m^\top \mAh$ and $\mAh \bc$, with entries of $\bc$ defined by \eqref{eq:Cj-def} below. These are two additional matrix-vector products requiring $O(\nnz(\mA))$ arithmetic operations; all remaining arithmetic operations take $O(m + n) = O(\nnz(\mA)).$

    Let $\Pi'_{ij} = \frac{x_{0,i} A_{ij}y_{0, j}([\bpi_0^\top \bone_m]_j + \epstr)^{-\beta}}{\sum_{j'=1}^n A_{ij'}y_{0, j'}([\bpi_0^\top \bone_m]_{j'} + \epstr)^{-\beta}}$ be the right-hand side of \eqref{eq:minimizer-Pi_*} for $\bpi_* = \bpi_0$ and $\bx = \bx_0, \by = \by_0.$ Since the matrix row normalization  applies the same rescaling on the numerator and the denominator in $\Pi'_{ij},$ we also have
    \[
        \Pi'_{ij} = \frac{x_{0,i} \Ah_{ij}y_{0, j}([\bpi_0^\top \bone_m]_j + \epstr)^{-\beta}}{\sum_{j'=1}^n \Ah_{ij'}y_{0, j'}([\bpi_0^\top \bone_m]_{j'} + \epstr)^{-\beta}}.
    \]
    By the lemma statement, $\x_{0, i} = (\bpi_0 \bone_n)_i$ and $y_{0, j} \propto \frac{(\bpi_0^\top \bone_m)_j}{(\mAh^\top \bone_m)_j}([\bpi_0^\top \bone_m]_j + \epstr)^{\beta}$. Let 
    \begin{equation}\label{eq:Cj-def}
    C_j := \Big(\frac{(\mAh^\top\bone_m)_j^{1/\beta}}{\sum_{j'=1}^n (\mAh^\top\bone_m)_{j'}^{1/\beta}} + 2\epstr\Big)^{1-\beta},
    \end{equation} 
    so that $\Pi_{0, ij} \propto \Ah_{ij}C_j.$ Then $y_{0, j} \propto C_j ([\bpi_0^\top \bone_m]_j + \epstr)^{\beta}.$ 
    Plugging into $\Pi_{ij}'$, we get
    \begin{equation}\notag
        \Pi_{ij}' = \frac{(\bpi_0 \bone_n)_i \Ah_{ij} C_j}{\sum_{j'=1}^n \Ah_{ij'} C_{j'}} = \frac{(\bpi_0 \bone_n)_i \Pi_{0, ij}}{\sum_{j' = 1}^n \Pi_{0, ij'}} = \Pi_{0, ij}.
    \end{equation}
    Thus, $\bpi_0$ satisfies \eqref{eq:minimizer-Pi_*} and, by construction, $\bpi_0 \in \Delta_\mA, \bpi_0 \bone_n = \bx_0.$ It follows by \Cref{lemma:Pi_*} that $\bpi_0$ is the minimizer in the definition of $\psi.$ 

    Let us now argue that $\by_0 \in \Deltatr_n$ ($\bx_0 \in \relint(\Delta_m)$ is immediate, by construction). Due to normalization, we have $\by_0 \in \Delta_n,$ so we only need to argue that $\by_0 \geq \epstr \bone_n.$ Note that we have already argued $y_{0, j} \propto C_j ([\bpi_0^\top \bone_m]_j + \epstr)^{\beta},$ where $C_j$ is defined by \eqref{eq:Cj-def}. It is immediate that $C_j \geq (2\epstr)^{1-\beta}$ and  $([\bpi_0^\top \bone_m]_j + \epstr)^{\beta} \geq \epstr^{\beta},$ thus $C_j ([\bpi_0^\top \bone_m]_j + \epstr)^{\beta} \geq 2^{1-\beta} \epstr.$ Since $n \geq 2$,  we have $\beta = \frac{1}{1 + \log(n)} \leq 2/3,$ so  $2^{1-\beta} > 5/4$. Thus, to prove that $y_{0, j} \geq \epstr,$ it suffices to argue that the normalization constant $\sum_{j=1}^n C_j ([\bpi_0^\top \bone_m]_j + \epstr)^{\beta} \leq 5/4.$ Recalling the definition of $C_j$ from \eqref{eq:Cj-def} and using the weighted AM-GM inequality (by which $a^{1-\beta}b^\beta \leq (1-\beta)a + \beta b$ for $a, b > 0$), we have the following bound on the normalization factor for $\by_0$:
    \begin{align}
       \sum_{j=1}^n C_j ([\bpi_0^\top \bone_m]_j + \epstr)^{\beta} &\leq (1-\beta)\sum_{j=1}^n \Big(\frac{(\mAh^\top\bone_m)_j^{1/\beta}}{\sum_{j'=1}^n (\mAh^\top\bone_m)_{j'}^{1/\beta}} + 2\epstr\Big)  + \beta\sum_{j''=1}^n ([\bpi_0^\top \bone_m]_{j''} + \epstr)\notag\\
        &= (1-\beta)(1 + 2n\epstr) + \beta(1 + n\epstr)\notag\\
        &= 1 + (2-\beta)n\epstr \leq 5/4, \label{eq:Y-norm}
    \end{align}
    where we have used $n \epstr \leq 1/8.$ Thus, $\by_0 \in \Deltatr_n.$

    Fix any $\bx \in \Delta_m, \by \in \Deltatr_n$, $\bz := (\bx, \by)$, and let $\bpi_*$ be the minimizing coupling in the definition of $\psi(\bx, \by)$. Then
from \Cref{lemma:dgf-properties}$(ii)$, Eq.~\eqref{eq:D-psi}, we have
\begin{equation}\label{eq:breg-div-of-dgf}
    D_\psi\big(\bz, \bz_0\big) = D_h(\by, \by_0) + D_h(\bpi_*, \bpi_0\diag(\by\oslash \by_0)) + \beta D_h\big(\bpi_*^\top \bone_m + \epstr \bone_n,\, \bpi_0^\top \bone_m + \epstr \bone_n\big). 
\end{equation}
Since $\by \in \Deltatr_n,$ it is immediate that $D_h(\by, \by_0) \leq \log(1/\epstr),$ and, similarly, $\beta D_h\big(\bpi_*^\top \bone_m + \epstr \bone_n,\, \bpi_0^\top \bone_m + \epstr \bone_n\big) = O(\beta \log(1/\epstr)) = O(\log(1/\epstr)).$ Thus, we only need to bound the middle Bregman divergence term on the right-hand side of \eqref{eq:breg-div-of-dgf}. To do so, we use that $\by \in \Deltatr_n$ (so $y_j \in [\epstr, 1)$, $j\in [n]$) and bound $\Pi_{0, ij}/y_{0, j}.$ In particular, recalling the definition of $C_j$ from \eqref{eq:Cj-def} and $y_{0, j} \propto C_j ([\bpi_0^\top \bone_m]_j + \epstr)^{\beta},$ we have:
\begin{align*}
    \frac{\Pi_{0, ij}}{y_{0, j}} &= \frac{\Ah_{ij}C_j}{C_j([\bpi_0^\top \bone_m]_j + \epstr)^{\beta}}\cdot \frac{\sum_{j'=1}^n C_{j'} ([\bpi_0^\top \bone_m]_{j'} + \epstr)^{\beta}}{\sum_{j''=1}^n (\mAh^\top \bone_m)_{j''} C_{j''}}.
\end{align*}
On the other hand, using \Cref{lemma:Pi_*} (and, again, that row normalization rescales the numerator and the denominator by the same quantity), we have
\begin{equation}\notag
    \frac{\Pi_{*, ij}}{y_j} = \frac{x_i \Ah_{ij}([\bpi_*^\top \bone_m]_j + \epstr)^{-\beta}}{\sum_{j'=1}^n \Ah_{ij'}y_{j'}([\bpi_*^\top \bone_m]_{j'} + \epstr)^{-\beta}}.
\end{equation}
Thus, we have
\begin{align}\notag
    \frac{\Pi_{*, ij}}{y_j} \cdot \frac{y_{0, j}}{\Pi_{0, ij}} = \frac{x_i ([\bpi_*^\top \bone_m]_j + \epstr)^{-\beta}([\bpi_0^\top \bone_m]_j + \epstr)^{\beta}}{\sum_{j'=1}^n \Ah_{ij'}y_{j'}([\bpi_*^\top \bone_m]_{j'} + \epstr)^{-\beta}}\cdot  \frac{\sum_{j''=1}^n (\mAh^\top \bone_m)_{j''} C_{j''}}{\sum_{j'=1}^n C_j ([\bpi_0^\top \bone_m]_j + \epstr)^{\beta}}.
\end{align}
Observe that $x_i \in [0, 1],$ $[\bpi_*^\top \bone_m]_j + \epstr, [\bpi_0^\top \bone_m]_j + \epstr \in [\epstr, 1 + \epstr],$ $\mAh\bone_n = \bone_m$ (by construction), $\sum_{j'=1}^n \Ah_{ij'}y_{j'} \geq \epstr (\mAh\bone_n)_i = \epstr,$ and $C_j \in [(2\epstr)^{1-\beta}, (1 + 2\epstr)^{1-\beta}].$ As a result, we have
\begin{align} \label{eq: Pi 1st term}
    \frac{\Pi_{*, ij}}{y_j} \cdot \frac{y_{0, j}}{\Pi_{0, ij}} \leq \frac{\epstr^{-2\beta}}{\epstr(1+\epstr)^{-\beta}}\frac{m (1 + 2\epstr)^{1-\beta}}{ n 2^{1-\beta}\epstr} =:B.
\end{align}
Since $\beta \in (0, 1),$ the expression $B$ above is bounded as $B=O\big(\frac{m}{\epstr^{2 + 2\beta}}\big)= m\,  \mathrm{poly}(1/\epstr)$.

We also need to bound the mass of
$\bpi_0\diag(\by\oslash\by_0)$.
Write $p_j=[\bpi_0^\top\mathbf{1}_m]_j$ and
$a_j=C_j^{1/(1-\beta)}-2\epstr$.
Then $\sum_j a_j=1$, $C_j=(a_j+2\epstr)^{1-\beta}$, and
\[
p_j=\frac{a_j^\beta C_j}{\sum_\ell a_\ell^\beta C_\ell}
\le a_j^\beta C_j
\le a_j+2\epstr,
\]
since the denominator is at least $\sum_\ell a_\ell=1$.
Using \eqref{eq:Y-norm}, we get
\[
y_{0,j}\ge \frac45 C_j(p_j+\epstr)^\beta
=\frac45(a_j+2\epstr)^{1-\beta}(p_j+\epstr)^\beta
\ge \frac45p_j.
\]
Thus the additional mass term in
$D_h(\bpi_*,\bpi_0\operatorname{diag}(\by\oslash\by_0))$ satisfies
\[
\sum_{i,j}\Pi_{0,ij}\frac{y_j}{y_{0,j}}
-\sum_{i,j}\Pi_{*,ij}
=\sum_j\frac{p_jy_j}{y_{0,j}}-1
\le \frac14.
\]
Overall,
\[
D_h(\bpi_*, \bpi_0\diag(\by\oslash \by_0)) \leq
\log B + \frac{1}{4}
=O(\log(m) + \log(1/\epstr)) = O(\log(mn/\epsilon)),
\]
completing the proof. 
\end{proof}

\subsection{Efficient Implementation of Prox Solves}\label{sec:implementation}

We now show that each prox update in \eqref{eq:relaxed-MP} can be
implemented using $O(1)$ matrix-vector multiplications. 
For error parameter $\xi > 0,$ a prox center $\bz\in\Zcal$ whose $\bx$-component is strictly positive,
and a vector $\bg\in\R^{m+n}$, the required approximation is a vector 
$\bw\in\Zcal$ satisfying
\begin{equation}\label{eq:prox-three-point}
    \bg^\top(\bw-\bu)
    \leq D_\psi(\bu,\bz)-D_\psi(\bu,\bw)-D_\psi(\bw,\bz)+\xi,
    \qquad \forall \bu\in\Zcal,
\end{equation}
which, using the definition of Bregman divergence, is equivalent to
\begin{equation}\label{eq:prox-first-order}
    \bigl(\bg+\nabla\psi(\bw)-\nabla\psi(\bz)\bigr)^\top
    (\bw-\bu)\leq\xi,
    \qquad \forall \bu\in\Zcal.
\end{equation}
Our implementation stores an approximation of the gradient at a
current ``center''. We first compute it at the initial point, and then
show how to maintain it during the updates.

\paragraph{Computing the first gradient.}
Let $\bz_0=(\bx_0,\by_0)$ and $\bpi_0$ be the point and coupling
constructed in \Cref{lemma:initialization}.
Then $\bx_0>\bzero$, $\bpi_0$ is the minimizing coupling in
\eqref{eq:final-dgf} at $\bz_0$, and its column sums
$\bpi_0^\top\bone_m$ are already known (and, by \Cref{lemma:initialization}, computed with $O(1)$ matrix-vector products). 
By \eqref{eq:grad-psi-wrt-x}, \eqref{eq:grad-psi-wrt-y},
and \eqref{eq:minimizer-Pi_*}, it holds that
\begin{equation}\label{eq:prox-initial-gradient}
    \nabla\psi(\bz_0)=
    \begin{pmatrix}
        (1+\beta)\bone_m+\log(\bx_0)
        -\log\!\left(
            \mA\bigl(
                \by_0\odot
                (\bpi_0^\top\bone_m+\epstr\bone_n)^{-\beta}
            \bigr)
        \right)\\
        \log(\by_0)+\bone_n
        -(\bpi_0^\top\bone_m)\oslash\by_0
    \end{pmatrix}.
\end{equation}
This requires one additional matrix-vector multiplication.
Set $\widehat{\nabla\psi}(\bz_0):=\nabla\psi(\bz_0)$.
Thus the initial gradient has zero approximation error, and its computation uses only $O(1)$ matrix-vector
multiplications.

The following lemma computes a prox update using the stored gradient
approximation at its center, and returns a gradient approximation
at its output, to be used by the next iterate.

\begin{lemma}\label{lem:efficient-prox-solves}
Fix $\xi\in(0,1)$, a center $\bz\in\Zcal$ whose $\bx$-component
is strictly positive, a vector $\bg\in\R^{m+n}$, and suppose we are given a vector $\widehat{\nabla\psi}(\bz)$ satisfying
\begin{equation}\label{eq:prox-center-gradient}
    \bigl\|\widehat{\nabla\psi}(\bz)-\nabla\psi(\bz)\bigr\|_\infty
    \leq \xi/16.
\end{equation}
From these inputs, we can compute a vector $\bw\in\Zcal$ satisfying
\eqref{eq:prox-three-point}, with a strictly positive $\bx$-component,
and a vector $\widehat{\nabla\psi}(\bw)$ satisfying $\|\widehat{\nabla\psi}(\bw)-\nabla\psi(\bw)\|_\infty\leq \xi/16$, using two matrix-vector multiplications, \[
O\Big(\nnz(\mA)+n\log(n)\log^2\Big(\frac{n}{\epsilon \xi}\Big)\Big)\]
arithmetic operations, and polylogarithmic parallel depth in $m, n, 1/\epsilon, 1/\xi$.
\end{lemma}

\begin{proof}
Throughout the proof, we denote $\bw = (\bx, \by)$ and $(\bb,\bd):=\bg-\widehat{\nabla\psi}(\bz)$ where $\bb\in\reals^m,~\bd\in\reals^n$.
 
We first argue that the required prox computation can be reduced to a block-separable problem, which we then argue can be solved efficiently. 
Expanding the Bregman divergence $D_\psi(\bw,\bz)$ gives
\[
    \bg^\top\bw + D_\psi(\bw,\bz)
    =
    \psi(\bw)+(\bg-\nabla\psi(\bz))^\top\bw
    +\nabla\psi(\bz)^\top\bz-\psi(\bz).
\]
We therefore approximately minimize 
$\psi(\bx,\by)+\bb^\top\bx+\bd^\top\by$ over $\Zcal$ (the rest of the expression is independent of $\bw$) 
and account for the gradient error $(\widehat{\nabla\psi}(\bz) - \nabla \psi(\bz))^\top \bw$ using \eqref{eq:prox-center-gradient}. 

Since $\psi$ itself involves a minimization problem and since $\bx$ is coupled with the rest of the variables only through $\bpi\bone_n = \bx$, we can instead simply require that $\bpi\bone \in \Delta_m$ (which is already ensured by $\bpi \in \Delta_\mA$), substitute  $\bpi\bone_n$ for $\bx$ to eliminate it from minimization, and assign $\bx = \bpi\bone_n$ at an approximate solution. In particular, we use
\begin{align*}
   &\; \min_{\bx \in \Delta_m, \by \in \Deltatr_n}\big\{\psi(\bx,\by)+\bb^\top\bx+\bd^\top\by\big\}\\
  \Leftrightarrow &\; \min_{\bx \in \Delta_m, \by \in \Deltatr_n, \bpi \in \Delta_A: \bpi \bone_n = \bx} \big\{  \bb^\top\bx+\bd^\top\by + \inner{\bpi, \log(\bpi\oslash \mA)} - \inner{\bpi^\top \bone_m, \log(\by)} + \beta\, h\big(\bpi^\top \bone_m + \epstr \bone_n\big) + h(\by)\big\} \\
  \Leftrightarrow &\; \min_{\by \in \Deltatr_n, \bpi \in \Delta_\mA} \big\{\bb^\top\bpi \bone_n +\bd^\top\by + \inner{\bpi, \log(\bpi\oslash \mA)} - \inner{\bpi^\top \bone_m, \log(\by)} + \beta\, h\big(\bpi^\top \bone_m + \epstr \bone_n\big) + h(\by)\big\}.
\end{align*}
This is now an optimization problem that only involves $\bpi$ and $\by.$ Let us introduce an auxiliary variable $\bq := \bpi^\top \bone_m$ to fix the column marginals. Then the above problem is equivalent to
\begin{align}\label{eq:prox-min-yqpi}
     \min_{\by \in \Deltatr_n, \bq \in \Delta_n, \bpi \in \Delta_\mA: \bpi^\top\bone_m = \bq} \big\{\bb^\top\bpi \bone_n +\bd^\top\by + \inner{\bpi, \log(\bpi\oslash \mA)} - \inner{\bq, \log(\by)} + \beta\, h\big(\bq + \epstr \bone_n\big) + h(\by)\big\}.
\end{align}
Now minimization over $\bpi$ can be performed in closed form, using first-order optimality conditions, to get
\begin{equation}\label{eq:prox-recovery}
    \Pi_{ij}=e^{-b_i}A_{ij}\frac{q_j}{(\mA^\top e^{-\bb})_j}, \quad \bpi\bone_n = e^{-\bb}\odot \mA(\bq\oslash (\mA^\top e^{-\bb})).
\end{equation}

Substituting \eqref{eq:prox-recovery} back into \eqref{eq:prox-min-yqpi} reduces the problem to
\begin{equation}\label{eq:separable-prob}
    \min_{\bq\in\Delta_n,\ \by\in\Deltatr_n} f(\bq, \by), \quad\text{ where }\quad
    f(\bq, \by) := 
        \sum_{j=1}^n q_j\log\frac{q_j}{c_jy_j}
        +\beta\, h(\bq+\epstr\bone_n)+h(\by)+\bd^\top\by,
\end{equation}
where we define $c_j := (\mA^\top e^{-\bb})_j,$ for notational convenience, and where $e^{-\bb} := (e^{-b_1}, \dots, e^{-b_m})$. 

The problem stated in \eqref{eq:separable-prob} has a block-separable feasible set, with blocks corresponding to $\bq$ and $\by.$ It further has other strong structural properties that make it efficiently solvable. We state the guarantee of such an algorithm in the claim below, and defer the proof of the claim to Appendix~\ref{appx:alt-min}.

\begin{claim}[Splitting solver]\label{claim:splitting-solve}
    Consider the problem stated in \eqref{eq:separable-prob} and let $\bq_*, \bv_*$ be its optimal solution. There is a deterministic algorithm that given \eqref{eq:separable-prob} and $\delta \in (0, \epstr)$ outputs feasible $\bq, \by$ such that $f(\bq, \by) - f(\bq_*, \bv_*) \leq \delta.$ The algorithm uses $O(\frac{n}{\beta}\log^2(\frac{n}{ \delta })) = O(n\log(n)\log^2(\frac{n}{\delta}))$ arithmetic operations and it is implementable in parallel with poly-logarithmic depth in $n$ and $1/\delta.$   
\end{claim}

To prove \Cref{lem:efficient-prox-solves}, we argue that the following choice of $\bw$ and $\widehat{\nabla\psi}(\bw)$ satisfies the lemma claim. Let $\bq, \by$ be the output of the algorithm from \Cref{claim:splitting-solve}, where $\delta > 0$ is chosen suitably small (to be specified later in the proof). Let $\widehat{q}_j \propto \frac{c_j y_j}{(q_j + \epstr)^\beta},$ where $\propto$ omits the normalization factor $Q := \sum_{j=1}^n \frac{c_j y_j}{(q_j + \epstr)^\beta}.$ Then define $\widehat{\bpi}$ via $\widehat{\Pi}_{ij} = A_{ij}e^{-b_i}\frac{\widehat{q}_j}{c_j}$, let $\bx = \widehat{\bpi}\bone_n$, and let $\bw := (\bx, \by)$. Since \(\mA\) has no zero columns, \(\bc > \bzero\); hence \(\bqh > \bzero\). Since \(\mA\) has no zero rows, \(\bx = e^{-\bb}\odot \mA(\bqh \oslash \bc) > \bzero\).

For analysis only, let $\bpi_\bw$ denote the minimizing coupling associated with $\bw.$ Observe that constructing this output $\bw$ requires only two matrix-vector products, for computing $\bc = \mA^\top e^{-\bb}$ and for multiplying $\widehat{\bpi}$ by $\bone_n$ ($\widehat{\bpi}$ does not need to be explicitly constructed). For $\widehat{\nabla\psi}(\bw),$ we choose the following vector:
\begin{equation}\label{eq:output-prox-center}
    \widehat{\nabla\psi}(\bw) := \begin{pmatrix}
        -\bb + (1 + \beta - \log (Q))\bone_m\\
        \log (\by) + \bone_n - \bqh\oslash \by
    \end{pmatrix}. 
\end{equation}
Thus, for the constructed output, the claimed bounds on the matrix-vector multiplications, arithmetic operations, and parallel depth follow from what has already been argued and \Cref{claim:splitting-solve} for $\delta = \mathrm{poly}(\xi, \epstr)$. 
It remains to argue that the constructed $\bw$ and $\widehat{\nabla\psi}(\bw)$ attain the claimed approximation errors. 

The reason for the specific choice of $\widehat{\bpi}$ that we made above is that, as defined, it satisfies the following for all $ij$ such that $A_{ij} > 0$:
\begin{equation}
    \log\Big(\frac{\widehat{\Pi}_{ij}}{A_{ij}y_j}\Big) + \beta\log(q_j + \epstr) = -b_i - \log(Q).
\end{equation}
As a consequence, using the definition of $\bx,$ $\widehat{\nabla \psi}(\bw)$ admits the following form, which makes it more directly comparable to $\nabla \psi(\bw)$:
\begin{equation}\label{eq:output-prox-center-equiv}
      \widehat{\nabla\psi}(\bw) := \begin{pmatrix}
         (1 + \beta)\bone_m + \log(\bx) - \log(\mA(\by \odot (\bq + \epstr\bone_n)^{-\beta})) \\
        \log (\by) + \bone_n - \bqh\oslash \by
    \end{pmatrix}. 
\end{equation}
Recall that \Cref{lemma:dgf-properties}(ii) gives the expression for $\nabla \psi(\bw)$; subtracting it from \eqref{eq:output-prox-center-equiv} and taking the $\ell_\infty$ norm, we get:
\begin{align}
    \norm{\widehat{\nabla \psi}(\bw) - \nabla \psi(\bw)}_\infty =\;& \norm{\log\Big(\frac{\mA(\by\odot(\bpi_\bw^\top \bone_m + \epstr\bone_n)^{-\beta})}{\mA(\by\odot(\bq + \epstr\bone_n)^{-\beta})}\Big), (\bpi_\bw^\top \bone_m - \bqh)\oslash \by}_\infty\notag\\
    \leq \;& \max\Big\{\beta \norm{\log(\bpi_\bw^\top \bone_m + \epstr\bone_n) - \log(\bq + \epstr\bone_n)}_\infty,\,  \norm{(\bpi_\bw^\top \bone_m - \bqh)\oslash \by}_\infty \Big\} \notag\\
    \leq \; & \max\Big\{ \frac{\beta}{\epstr} \norm{\bpi_\bw^\top \bone_m - \bq}_1,\, \frac{1}{\epstr}\norm{\bpi_\bw^\top \bone_m - \bqh}_1 \Big\}, \label{eq:prox-center-grad-error-1}
\end{align}
where in the last inequality we used that the argument of the logarithm is at least $\epstr,$ thus it is  $\frac{1}{\epstr}$-Lipschitz and that $\by \in \Deltatr_n$ (so $\by \geq \epstr \bone_n$). 

Similar to the proof of \Cref{lemma:dgf-properties}, let $\Phi(\bpi, \by)$ be the function under the minimum in the definition of $\psi(\bx, \by)$ so that $\psi(\bx, \by) = h(\by) + \min_{\bpi \in \Delta_\mA: \bpi \bone_n = \bx}\Phi(\bpi, \by).$ Observe that for a fixed $\by,$ $D_\Phi((\widehat{\bpi}, \by), (\bpi_\bw, \by)) \geq D_h(\widehat{\bpi}, \bpi_\bw) \geq \frac{1}{2}\norm{\widehat{\bpi} - \bpi_\bw}_1^2 \geq \frac{1}{2}\norm{\bqh - \bpi_\bw^\top \bone_m}_1^2$. Thus $\Phi$ is 1-strongly convex w.r.t.\ the $\ell_1$ norm of the column marginals. A direct consequence is that
\begin{align}
    \norm{\bqh - \bpi_\bw^\top \bone_m}_1^2 &\leq \inner{\nabla_{\bpi} \Phi(\widehat{\bpi}, \by) - \nabla_{\bpi} \Phi(\bpi_\bw, \by), \widehat{\bpi} - \bpi_{\bw}}\notag\\
    &= \beta \Big\langle \log\Big(\frac{\bqh + \epstr \bone_n}{\bq + \epstr \bone_n}\Big), \bqh - \bpi_\bw^\top \bone_m \Big\rangle \notag\\
    &\leq \frac{\beta}{\epstr}\|\bqh - \bq\|_1 \norm{\bqh - \bpi_\bw^\top \bone_m}_1.\notag
\end{align}
We thus deduce from the last inequality that
\begin{equation}\notag
    \norm{\bqh - \bpi_\bw^\top \bone_m}_1 \leq \frac{\beta}{\epstr}\|\bqh - \bq\|_1.
\end{equation}
Let $\bq_*, \bv_*$ be the optimal solution to \eqref{eq:separable-prob}, as in \Cref{claim:splitting-solve}. 
Plugging the last inequality back into \eqref{eq:prox-center-grad-error-1}, using the triangle inequality, and simplifying, we have
\begin{equation}\label{eq:prox-center-grad-error-corrected}
    \norm{\widehat{\nabla \psi}(\bw) - \nabla \psi(\bw)}_\infty \leq \frac{2\beta}{\epstr^2}\norm{\bqh - \bq}_1 \leq \frac{2\beta}{\epstr^2} \big(\norm{\bqh - \bq_*}_1 + \norm{\bq_* - \bq}_1\big).
\end{equation}
To bound $\norm{\bq_* - \bqh}_1,$ we observe that both $\bqh$ and $\bq_*$ take a similar, normalized form. In particular, the first-order optimality conditions for $\bq_*$ in minimization of $f(\bq, \by)$ imply $q_{*, j} = \frac{c_j v_{*, j}/(q_{*, j} + \epstr)^\beta}{\sum_{j'=1}^n c_{j'} v_{*, j'}/(q_{*, j'} + \epstr)^\beta}$. Recall that, by construction, $\widehat{q}_j = \frac{c_j y_j/(q_j + \epstr)^\beta}{\sum_{j'=1}^n c_{j'} y_{j'}/(q_{j'} + \epstr)^\beta}$. Using that the map $\ba \mapsto e^{\ba}/(\bone^\top e^{\ba})$ is $1$-Lipschitz from $\ell_\infty$ to $\ell_1$, as it is the gradient of the log-sum-exp function (see, e.g., \cite[Example 5.15]{beck2017first}), we have:
\begin{align}
    \norm{\bq_* - \bqh}_1 &\leq \norm{\log(\bv_*) - \log(\by) - \beta(\log(\bq_* + \epstr \bone_n) - \log(\bq + \epstr \bone_n) )}_\infty\notag\\
    &\leq \frac{1}{\epstr}\norm{\bv_* - \by}_1 + \frac{\beta}{\epstr} \norm{\bq - \bq_*}_1,\label{eq:q*-qhat}
\end{align}
where we have again used that the logarithm is $\frac{1}{\epstr}$-Lipschitz since its argument is at least $\epstr.$ 
To bound the error terms arising from approximately minimizing $f(\bq, \by)$ and obtain a final bound on $\norm{\widehat{\nabla \psi}(\bw) - \nabla \psi(\bw)}_\infty,$ we will use that $f$ is relatively strongly convex with respect to the entropy function of its individual arguments, summarized in the following claim.
\begin{claim}\label{claim:f-rel-str-cvx}
    Let $f$ be defined by \eqref{eq:separable-prob}. Then  $f(\bq, \by)$ is both $\beta$-relatively strongly convex with respect to $h(\bq)$ and $1$-relatively strongly convex with respect to $h(\by).$  
\end{claim}
\begin{proof}
    For the first claim, it suffices to argue that $f(\bq, \by) - \beta h(\bq)$ is jointly convex in $\bq, \by.$ To do so, observe that:
    \begin{equation}\notag
        \begin{aligned}
            f(\bq, \by) - \beta h(\bq) =\; & (1-\beta)\big(h(\bq) - \bq^\top \log(\by) + h(\by)\big)\\
            &+ \beta \big(h(\bq + \epstr \bone_n) - \bq^\top \log(\by) + h(\by)\big)\\
            &- \bq^\top \log(\bc) + \bd^\top \by.
        \end{aligned}
    \end{equation}
    Comparing to \eqref{eq:psi-rel-convex-h(x)}, convexity of $f(\bq, \by) - \beta h(\bq)$ follows by the same line of argument.

Convexity of $f(\bq, \by) - h(\by)$ is immediate, as all terms appearing in $f(\bq, \by)$ (which includes $h(\by)$ that gets subtracted) are jointly convex in $\bq, \by.$     
\end{proof}

 By \Cref{claim:splitting-solve,claim:f-rel-str-cvx}, we have that:
\begin{equation}\notag
  \max\big\{\beta D_h(\bq, \bq_*), D_h(\by, \bv_*)\} \leq f(\bq, \by) - f(\bq_*, \bv_*) \leq \delta.
\end{equation}
Since $h$ is 1-strongly convex with respect to the $\ell_1$ norm, this implies:
\begin{equation}\label{eq:dist-qq*-yv*}
    \norm{\bq - \bq_*}_1 \leq \sqrt{2\delta/\beta} \quad \text{ and }\quad \norm{\by - \bv_*}_1 \leq \sqrt{2\delta}.
\end{equation}
Plugging \eqref{eq:q*-qhat} and \eqref{eq:dist-qq*-yv*} back into \eqref{eq:prox-center-grad-error-corrected} and simplifying, we finally have
\begin{equation}\notag
     \norm{\widehat{\nabla \psi}(\bw) - \nabla \psi(\bw)}_\infty \leq \frac{2\sqrt{2\delta}(\sqrt{\beta}(\epstr + \beta) + \beta)}{\epstr^3} \leq \frac{6\sqrt{2\delta}}{\epstr^3}.
\end{equation}
Thus, to have $\|\widehat{\nabla \psi}(\bw) - \nabla \psi(\bw)\|_\infty \leq \xi/16,$ it suffices that $\delta \leq {2^{-16} \epstr^{6}}\xi^2,$ so we choose $\delta = {2^{-16} \epstr^{6}}\xi^2.$

It remains to argue the approximate first-order condition corresponding to the three-point inequality, which was stated in \eqref{eq:prox-first-order}. To do so, we decompose $\bigl(\bg+\nabla\psi(\bw)-\nabla\psi(\bz)\bigr)^\top
    (\bw-\bu)$ and then bound individual terms using what we have proved so far. In particular,
\begin{equation}\label{eq:approx-first-order-cond}
    \begin{aligned}
        \bigl(\bg+\nabla\psi(\bw)-\nabla\psi(\bz)\bigr)^\top
    (\bw-\bu) =\;&  \bigl(\nabla\psi(\bw) - \widehat{\nabla\psi}(\bw)\bigr)^\top
    (\bw-\bu) 
    + \bigl(\widehat{\nabla\psi}(\bz) - \nabla\psi(\bz)\bigr)^\top
    (\bw -\bu)\\
    &+ \bigl(\bg+\widehat{\nabla\psi}(\bw)-\widehat{\nabla\psi}(\bz)\bigr)^\top
    (\bw-\bu). 
    \end{aligned}
\end{equation}
For the terms in the first line, we immediately get the following bound from \eqref{eq:prox-center-gradient}, $\|\widehat{\nabla \psi}(\bw) - \nabla \psi(\bw)\|_\infty \leq \xi/16$, and $\norm{\bw - \bu}_1 \leq 4:$
\begin{equation}\label{eq:approx-old-prox-center-contribution}
    \bigl(\nabla\psi(\bw) - \widehat{\nabla\psi}(\bw)\bigr)^\top
    (\bw-\bu) 
    + \bigl(\widehat{\nabla\psi}(\bz) - \nabla\psi(\bz)\bigr)^\top
    (\bw -\bu) \leq {\xi}/{2}.
\end{equation}
For the remaining term, observe that $\bg+\widehat{\nabla\psi}(\bw)-\widehat{\nabla\psi}(\bz) = (\bb, \bd) + \widehat{\nabla\psi}(\bw)$. Let $\bu = (\bu_\bx, \bu_\by)$. Since
$\bb+\widehat{\nabla_{\bx}\psi}(\bw)$ is constant and the $\bx$-components of both $\bw, \bu$ are on the probability simplex, we have that $\inner{\bb + \widehat{\nabla_{\bx}\psi}(\bw), \bx - \bu_\bx} = 0.$ Thus, we have the following bound for the third term in \eqref{eq:approx-first-order-cond}:
\begin{align}
    \bigl(\bg+\widehat{\nabla\psi}(\bw)-\widehat{\nabla\psi}(\bz)\bigr)^\top
    (\bw-\bu) &= \bigl(\bd +\widehat{\nabla_{\by}\psi}(\bw)\bigr)^\top
    (\by-\bu_\by)\notag\\
    &= \bigl({\nabla_{\by} f}(\bqh, \by)\bigr)^\top
    (\by-\bu_\by)\notag\\
   &= \bigl({\nabla_{\by} f}(\bqh, \by) - {\nabla_{\by} f}(\bq_*, \bv_*)\bigr)^\top
    (\by-\bu_\by) + \big({\nabla_{\by} f}(\bq_*, \bv_*)\big)^\top(\by-\bu_\by)\notag\\
    &\leq 2\norm{{\nabla_{\by} f}(\bqh, \by) - {\nabla_{\by} f}(\bq_*, \bv_*)}_\infty + \big({\nabla_{\by} f}(\bq_*, \bv_*)\big)^\top(\by-\bu_\by), \label{eq:third-term-bnd-1}
\end{align}
where we used H\"{o}lder's inequality along with $\norm{\by-\bu_{\by}}_1 \leq 2$.  

Since $\nabla_\by f(\bq, \by) = \log(\by) + \bone_n + \bd - \bq\oslash \by,$ using the same Lipschitz arguments as before,
\begin{align}\label{eq:grad-error-f}
    \norm{{\nabla_{\by} f}(\bqh, \by) - {\nabla_{\by} f}(\bq_*, \bv_*)}_\infty &\leq \frac{1}{\epstr}\norm{\by - \bv_*}_1 + \frac{1}{\epstr}\norm{\bqh - \bq_*}_1 + \frac{1}{\epstr^2}\norm{\by - \bv_*}_1\notag\\
    &\leq \frac{\sqrt{2\delta}(2 + \epstr + \sqrt{\beta})}{\epstr^2} \leq \frac{3\xi}{16},
\end{align}
where we have used already derived bounds \eqref{eq:q*-qhat}, \eqref{eq:dist-qq*-yv*} and the choice of $\delta.$

Finally, from the optimality of $\bv_*$ and $\delta$-optimality gap at $\bq, \by,$ we can infer the following two inequalities:
\begin{gather*}
    \inner{\nabla_\by f(\bq_*, \bv_*), \by - \bv_*} = \inner{\nabla_\bq f(\bq_*, \bv_*), \bq - \bq_*} +  \inner{\nabla_\by f(\bq_*, \bv_*), \by - \bv_*} \leq f(\bq, \by) - f(\bq_*, \bv_*) \leq \delta,\\
    \inner{\nabla_\by f(\bq_*, \bv_*), \bv_* - \bu_\by} \leq 0,
\end{gather*}
where in the first line we have used that $\inner{\nabla_\bq f(\bq_*, \bv_*), \bq - \bq_*} = 0,$ which holds because $\bq_* > \bzero$ and so simplex stationarity makes all coordinates of $\nabla_\bq f(\bq_*, \bv_*)$ equal.  
Adding these two inequalities, we get
\begin{equation}\label{eq:last-one-yay!}
    \big({\nabla_{\by} f}(\bq_*, \bv_*)\big)^\top(\by-\bu_\by) \leq \delta \leq \xi/2^{16}. 
\end{equation}
Thus, $\bigl(\bg+\widehat{\nabla\psi}(\bw)-\widehat{\nabla\psi}(\bz)\bigr)^\top
    (\bw-\bu) \leq \xi/2,$ which, together with \eqref{eq:approx-first-order-cond}, \eqref{eq:approx-old-prox-center-contribution} leads to \eqref{eq:prox-first-order}. 
\end{proof}
\paragraph{Applying the oracle to the mirror-prox updates.}
Set $\xi=\epsilon/(12L)$ and start with the exact gradient
$\widehat{\nabla\psi}(\bz_0)=\nabla\psi(\bz_0)$ from
\eqref{eq:prox-initial-gradient}.
At iteration $k$ of \eqref{eq:relaxed-MP}, apply
\Cref{lem:efficient-prox-solves} first with center $\bz_k$ and
$\bg=\mF(\bz_k)/L$, obtaining $\bw_k$, and then with the same
center and $\bg=\mF(\bw_k)/(2L)$, obtaining $\bz_{k+1}$.
Both calls use the stored approximation
$\widehat{\nabla\psi}(\bz_k)$.
Keep the approximation $\widehat{\nabla\psi}(\bz_{k+1})$
returned by the second call for the next iteration.
The lemma ensures inductively that every center has a strictly
positive $\bx$-component and satisfies
\eqref{eq:prox-center-gradient} with the same tolerance $\xi$.
Thus both updates satisfy \eqref{eq:prox-three-point}, and
\Cref{fact:op-extrapolation-convergence} bounds the contribution
of the prox errors to the MVI regret by
$3L\xi=\epsilon/4$.

Each prox call uses two matrix-vector multiplications.
Computing $\mF(\bx,\by)$ in \eqref{eq:mon-op-def}
requires the two products $\mA\by$ and
$\mA^\top(\bx\oslash(\mA\by))$.
Hence each iteration uses at most eight matrix-vector
multiplications and
\[
O\Big(\nnz(\mA)+n\log(n)\log^2\Big(\frac{n}{\epsilon}\Big)\Big)
\]
arithmetic operations, where we have used 
$L=\Theta(\log n)$, by \Cref{lemma:dgf-properties}$(iii)$ and $\beta = \Theta(1/\log(n))$.

Finally, for completeness, we provide the pseudocode for implementing approximate prox solves in \Cref{lem:efficient-prox-solves}. The auxiliary procedure (\Cref{algo:alt-min}) can be found in Appendix~\ref{appx:alt-min}.

\begin{algorithm}[ht]
    \caption{Approximate Prox Solve}\label{algo:approx-prox}
    \DontPrintSemicolon
\small
\KwIn{$\xi, \bg, \widehat{\nabla \psi}(\bz), n, \beta, \epstr$, matrix-vector multiplication oracle for $\mA$.}
\KwOut{$\bw$ that satisfies \eqref{eq:prox-first-order} and $\widehat{\nabla \psi}(\bw)$ that satisfies \eqref{eq:prox-center-gradient} (with respect to $\bw$)}

Set $(\bb, \bd):= \bg - \widehat{\nabla \psi}(\bz)$, $\bc = \mA^\top\exp(-\bb)$ \;

Call \Cref{algo:alt-min} with parameters $\epstr, n, \beta, \bc, \bd, \delta = 2^{-16}\epstr^6 \xi^2$ and let $\bq, \by$ be its output\;

Set
$\widehat{q}_j = \frac{c_j y_j}{(q_j + \epstr)^\beta}$ for $j \in [n],$ $Q \gets (\bone_n^\top \bqh),$ $\bqh \gets \bqh/ Q,$ $\bx = e^{-\bb}\odot \mA(\bqh \oslash \bc)$\;

$\bw \gets (\bx, \by)$, $\widehat{\nabla\psi}(\bw) \gets (
        -\bb + (1 + \beta - \log (Q))\bone_m, \, 
        \log (\by) + \bone_n - \bqh\oslash \by)$

\Return{$\bw$, $\widehat{\nabla \psi}(\bw)$}
\end{algorithm}

\subsection{Putting Everything Together}\label{sec:final-upper}

We now combine all the results from this section to complete the proof of \Cref{thm:main-upper-bound}.

\begin{proof}[Proof of \Cref{thm:main-upper-bound}]
    Let \(c_0>0\) be an absolute constant such that \Cref{lemma:initialization} gives \(\sup_{\bz\in \Zcal}D_\psi(\bz,\bz_0)\le c_0\log(m/\epstr)\), and set \(C = 8c_0\). With the parameters of \Cref{alg:packing-covering-mirror-prox}, \Cref{lemma:dgf-properties}$(iii)$ ensures that \(\mF\) is \(L = O(\log(n))\)-relaxed relatively Lipschitz with respect to \(\psi\), while \Cref{lem:efficient-prox-solves}, initialized using \eqref{eq:prox-initial-gradient}, implements both prox updates with tolerance \(\xi=\epsilon/(12L)\). Thus, \Cref{fact:op-extrapolation-convergence} bounds the MVI regret by \(2 L c_0\log(m/\epstr)/(K+1)+3L\xi\le\epsilon/4+\epsilon/4=\epsilon/2\). Applying \Cref{lemma:PCG-to-Minty-gap} to the sequence \(\{\bw_k\}_{k=0}^K\), with \(\epstr = \epsilon/(8n)\), gives the desired multiplicative guarantee, since \(e^{\epsilon/2}/(1-\epsilon/8)\le 1+\epsilon\) for \(\epsilon\in(0,1)\). Finally, \(K+1=O(\epsilon^{-1}\log n\log(mn/\epsilon))\); multiplying this iteration count by the per-iteration bounds established in \Cref{sec:implementation}, and including the initialization cost from \Cref{lemma:initialization}, gives the stated matrix-vector, arithmetic, parallel work, and depth bounds.
\end{proof}

\section{Nearly-Matching Oracle Lower Bound}\label{sec:lower}

In this section, we argue that the bound provided in \Cref{thm:main-upper-bound} is nearly tight, up to a single logarithmic factor, in terms of oracle complexity measured by the matrix-vector multiplications. Our main lower bound result is stated in \Cref{thm: main lower} below. The rest of this section first constructs the ``hard'' instance family in \Cref{sec:instance-fam} with basic properties of the instances proved in the same subsection, and then proceeds with the proof of \Cref{thm: main lower}.

\begin{theorem} \label{thm: main lower}
For any deterministic algorithm and any $T\geq 4$, let $K=\Theta(\frac{T}{\log(T)})$ be the smallest odd integer such that $K\lfloor\log_4 K\rfloor>T$. There exists a
positive matrix $\mA\in\reals_{>0}^{n\times n}$ with $n=O(T^3)$ such that
\[
\val(\mA):=
\min_{\bx\in\Delta_n}\max_{\by\in\Delta_n}
\bx^\top \mA\by
\geq
\frac{9}{K}+\frac{1}{K^2}~,
\]
yet after $T$ rounds of access to a two-sided matvec oracle of $\mA$, the algorithm
returns $\by_{T+1}\in\Delta_n$ such that
\[
\min_{\bx\in\Delta_n}\bx^\top \mA\by_{T+1}
=
\min_{i\in[n]}(\mA\by_{T+1})_i
\leq
\frac{9}{K}+\frac{1}{2K^2}~.
\]
Consequently,
\[
\frac{\val(\mA)}
{\min_{\bx\in\Delta_n}\bx^\top \mA\by_{T+1}}
\geq
1+\frac{1}{18K+1}
=1+\Theta\left(\frac{\log(T)}{T}\right)
~,
\]
and so achieving $(1+\eps)$-multiplicative error requires $\Omega(\frac{\log(1/\eps)}{\eps})$ oracle calls.
\end{theorem}

A couple of remarks are in order here.

\begin{remark}[Multiplicative error.]
    The bound on the multiplicative error in \Cref{thm: main lower} is stated in a slightly different form than what we used for the upper bound in \Cref{thm:main-upper-bound}. This is primarily for proof convenience. In fact, the stated lower bound also implies $\max_{j\in [n]}(\mA^\top \bx_{T+1})_j \geq \big(1 + \Theta\big(\frac{\log(T)}{T}\big)\big)\min_{i\in[n]}(\mA\by_{T+1})_i$ for any primal-dual pair $(\bx_{T+1}, \by_{T+1})\in \Delta_m \times \Delta_n$ that can be output by an oracle-based algorithm after $T$ queries. To see this, observe that the stated bound in \Cref{thm: main lower} is equivalent to
    \begin{equation}\notag
        \Big(1 + \Theta\Big(\frac{\log(T)}{T}\Big)\Big)\min_{i\in[n]}(\mA \by_{T+1})_i \leq \min_{\bx \in \Delta_m} \max_{\by \in \Delta_n}\bx^\top\mA\by \leq \min_{\bx \in \Delta_m}\max_{j\in [n]}(\mA^\top\bx)_j \leq \max_{j\in [n]}(\mA^\top\bx_{T+1})_j,\quad \forall \bx_{T+1} \in \Delta_m.
    \end{equation}
\end{remark}

\begin{remark}[Matrix width.]\label{rem:width}
The matrix $\mA$ constructed in the proof of Theorem~\ref{thm: main lower} has width 
\[
\mathrm{width}(\mA)\leq \frac{8+2/K}{1/K}=8K+2=\Theta\Big(\frac{T}{\log(T)}\Big),
\]
as established by Lemma~\ref{lem: A entry} in \Cref{sec:lb-proof}. The large (possibly super-polynomial in $m, n, 1/\epsilon$) width regime is what motivated the whole line of work on width-independent packing and covering LP solvers, started in the 1990s. Our results show that width-independence becomes important already at $\widetilde{O}(1/\epsilon)$ scale, since the matrix-vector oracle complexity of the lower bound is nearly matched by our upper bound from \Cref{sec:main-upper}. 
\end{remark}

\subsection{Hard Instance Family}\label{sec:instance-fam}

We now construct our hard instance family and prove basic properties for it. 
Let $\ell:=\lfloor\log_4 K\rfloor,~L:=K\ell,~M:=8(L+1)4^\ell,~n:=KM$. We partition $[n]$ into $K$ blocks of size $M$, and associate indices $[K]\times[M]\cong[n]$ as $(k,m)\sim (k-1)M+m$, namely $k$ indicates the block number and $m$ indicates the position inside the block.
Denote by $\bb_k\in\reals^n$ the indicator of block $k$, namely $(\bb_k)_{(i,j)}=\bone_{\{i=k\}}$,
and let $\bb_{K+1}:=\bb_{1}$.

We construct a hard matrix
by letting $\bu_t,\bv_t\in\mathrm{span}(\be_{\lceil t/\ell\rceil,1},\dots,\be_{\lceil t/\ell\rceil,M})\subset \reals^n$, $t\in[L]$, which will be specified soon according to a resisting oracle strategy,
and setting
\[
\mA:=\frac{1}{K}\bone\bone^{\top}+4\sum_{k\in [K]}\bb_k(\bb_k+\bb_{k+1})^\top+\frac{1}{2K}\bone(\bb_1+\bu_1)^\top+\sum_{t\in[L-1]}\bv_t(\bu_{t+1}-\bu_t)^\top~.
\]
Equivalently, by letting $\barbu_t,\barbv_t\in\reals^M$ be the restrictions of $\bu_t,\bv_t$ to block $\lceil t/\ell\rceil$, we can write $\mA$ explicitly as
\begin{align}
\mA
&=
\frac{1}{K}
\begin{pmatrix}
\bone_{M\times M} & \bone_{M\times M} &  \cdots & \bone_{M\times M}\\
\bone_{M\times M} & \bone_{M\times M} &  \cdots & \bone_{M\times M}\\
\vdots & \vdots & \ddots & \vdots\\
\bone_{M\times M} & \bone_{M\times M} & \cdots & \bone_{M\times M}
\end{pmatrix}+
4\begin{pmatrix}
\bone_{M\times M} & \bone_{M\times M} & 0 & \cdots & 0\\
0 & \bone_{M\times M} & \bone_{M\times M} & \cdots & 0\\
\vdots & & \ddots & \ddots & \vdots\\
0 & \cdots & 0 & \bone_{M\times M} & \bone_{M\times M}\\
\bone_{M\times M} & 0 & \cdots & 0 & \bone_{M\times M}
\end{pmatrix}
\label{eq: B explicit}
\\&~~~~~
+
\frac{1}{2K}
\begin{pmatrix}
\bone_M(\bone_M+\barbu_1)^\top & 0 & \cdots & 0\\
\bone_M(\bone_M+\barbu_1)^\top & 0 & \cdots & 0\\
\vdots & \vdots & & \vdots\\
\bone_M(\bone_M+\barbu_1)^\top & 0 & \cdots & 0
\end{pmatrix}
+
\begin{pmatrix}
H_1 & \barbv_\ell\barbu_{\ell+1}^\top & 0 & \cdots & 0\\
0 & H_2 & \barbv_{2\ell}\barbu_{2\ell+1}^\top & \cdots & 0\\
\vdots & \ddots & \ddots & \ddots & \vdots\\
0 & \cdots & 0 & H_{K-1}
  & \barbv_{(K-1)\ell}\barbu_{(K-1)\ell+1}^\top\\
0 & \cdots & \cdots & 0 & H_K
\end{pmatrix}~,
\nonumber
\end{align}
where
\[
H_k:=\sum_{t=(k-1)\ell+1}^{k\ell-1}\barbv_t(\barbu_{t+1}-\barbu_t)^\top
-\bone_{\{k<K\}}\barbv_{k\ell}\barbu_{k\ell}^\top~.
\]

The vectors
$\bu_t,\bv_t\in\reals^n$ are chosen for $t\in[L]$ iteratively as follows.
In each block $k\in[K]$, initialize $U_k^{+1}=V_k^{+1}=\{1,\dots,M/2\}$ and $U_k^{-1}=V_k^{-1}=\{M/2+1,\dots,M\}$.
For $t\leq T:$

\begin{itemize}
\item Receive an algorithm's query $(\bx_t,\by_t)$.
\item
For $k=\lceil t/\ell\rceil$, and each
$\sigma\in\{-1,1\}$, choose a vertex of the polytope
\[
[-1,1]^n\cap\mathrm{span}\{\be_{k,m}:m\in U_k^\sigma\}\,\cap\,\mathrm{span}(\by_1,\dots,\by_t,\bone)^\perp~,
\]
and let $\bu_t$ be the sum of these two vertices.
Replace $U_k^\sigma$ by $\{m\in U_k^\sigma:(\bu_t)_{(k,m)}=\sigma\}$ for each $\sigma$.

\item
For each $\sigma\in\{-1,1\}$, choose a vertex of the polytope
\[
[-1,1]^n\cap\mathrm{span}\{\be_{k,m}:m\in V_k^\sigma\}\,\cap\,\mathrm{span}(\bx_1,\dots,\bx_t,\bone)^\perp~,
\]
and let $\bv_t$ be the sum of these two vertices. Replace $V_k^\sigma$ by $\{m\in V_k^\sigma:(\bv_t)_{(k,m)}=\sigma\}$ for each $\sigma$.

\item Let $\mA_t:=\frac{1}{K}\bone\bone^{\top}+4\sum_{k\in [K]}\bb_k(\bb_k+\bb_{k+1})^\top+\frac{1}{2K}\bone(\bb_1+\bu_1)^\top+\sum_{k\in[t-1]}\bv_k(\bu_{k+1}-\bu_k)^\top$,
and feed the algorithm with oracle response $(\bx_t^\top \mA_t,\mA_t\by_t)$.
\end{itemize}
After the algorithm returns $\by_{T+1}$, set $\bx_{T+1}=\bzero$ and $\bx_t=\by_t=\bzero$ for $T+1<t\leq L$, and complete the remaining vector choices as above without further oracle calls.
Note that since $T+1\le L$, the orthogonality constraints imposed when choosing $\bu_L$ include the algorithm’s output $\by_{T+1}$, and so $\bu_L^\top \by_{T+1}=0$.

The lemma below provides basic properties of the constructed hard instance family that will be used in the lower bound proof. We remark here that we do not require the query vectors $(\bx_1, \by_1),\dots,(\bx_t, \by_t)$ to be supported on the probability simplex, or even nonnegative. 

\begin{lemma} \label{lem: ut,vt}
For all $t\in[L]$ it holds that:
\begin{enumerate}[label=(\roman*)]
    \item $\bu_t,\bv_t$ are both supported on block $\lceil t/\ell\rceil$;
    \item $\bu_t,\bv_t$ both have coordinates that sum to zero;
    \item $\bu_t,\bv_t$ satisfy $\bu_t\perp (\by_1,\dots,\by_t),~\bv_t\perp (\bx_1,\dots,\bx_t)$.
    \item $\bu_t,\bv_t\in[-1,1]^n$.
    \item In each block $k$ and for each $\sigma\in\{-1,1\}$, the final sets $U_k^\sigma,V_k^\sigma$ are nonempty, and all vectors $\bu_t$ in that block equal $\sigma$ on $U_k^\sigma$, while all vectors $\bv_t$ in that block equal $\sigma$ on $V_k^\sigma$.
\end{enumerate}
\end{lemma}

\begin{proof}
Items $(i), (iii)$ are immediate by construction, item $(ii)$ follows by orthogonality to $\bone$, and item $(iv)$ is clear since the two summands defining $\bu_t,\bv_t$ have disjoint supports.

To see item $(v)$, fix a block $k\in[K]$ and a sign
$\sigma\in\{-1,1\}$. We prove by induction on $j=0,\ldots,\ell$
that, after the first $j$ vectors $\bu_t$ in block $k$ have been
chosen, the current set $U_k^\sigma$ has cardinality at least
$M/(2\cdot4^j)$. The case $j=0$ follows from initialization.
For the inductive step, let $j\in[\ell]$, set
$t=(k-1)\ell+j$, and let $W$ denote $U_k^\sigma$ immediately
before choosing $\bu_t$. By the induction hypothesis and the
definition of $M$,
\[
    |W|\ge \frac{M}{2\cdot4^{j-1}}
    \ge 16(L+1)\ge 4(t+1)~.
\]

Let $\bw$ be the vertex chosen for the fixed sign $\sigma$, which is supported on $W$, and let
$\mathcal L:=\mathrm{span}\{\be_{k,m}:m\in W\}\,\cap\,\mathrm{span}(\by_1,\dots,\by_t,\bone)^\perp$ and $S:=\{(k,m):m\in W,~|(\bw)_{(k,m)}|=1\}$.
Note that
\begin{equation} \label{eq: S>M2}
|S|\geq \dim(\mathcal L)
\geq |W|-(t+1)>\frac{|W|}{2}~.
\end{equation}
To see why the first inequality above holds, note that if it did not hold, then there would have been a nonzero $\ba\in\mathcal L$ with $\ba_S=\bzero$, thus both $\bw+\delta\ba,\bw-\delta\ba\in\mathcal L\cap[-1,1]^n$ for sufficiently small $\delta>0$ which contradicts $\bw$ being a vertex.

Let $c_+,c_-$ count the coordinates of $\bw$ equaling $1,-1$, respectively, and let $r:=|W|-|S|\leq t+1$.
Since the coordinates of $\bw$ sum to zero, $|c_+-c_-|\leq r$, and hence
\[
\min\{c_+,c_-\}\geq\frac{|S|-r}{2}\geq\frac{|W|}{2}-(t+1)\geq\frac{|W|}{4}~.
\]
Restricted to $W$, the vector $\bu_t$ coincides with $\bw$. Hence the updated
set $U_k^\sigma$ has cardinality at least
$|W|/4\ge M/(2\cdot4^j)$, completing the inductive step.

After all $\ell$ stages, each $U_k^\sigma$ has at least $M/(2\cdot4^\ell)=4(L+1)>0$ coordinates. By nesting, all $\bu_t$ in block $k$ equal $\sigma$ on this final set. The proof for $\bv_t$ follows the same argument. 
\end{proof}

\subsection{Proof of Theorem~\ref{thm: main lower}}\label{sec:lb-proof}

The overview of the proof is as follows. First, in \Cref{lem: A entry}, we argue that the constructed matrix $\mA$ is entrywise positive, thus it corresponds to a valid instance of \eqref{eq:packing-covering-game}. The same lemma also implies a $\Theta(K)$ bound on the width of $\mA.$  In \Cref{lem: val}, we then bound below the value of \eqref{eq:packing-covering-game} for the constructed $\mA$. The final step of the proof is to argue that any dual ($\by$-player) solution that can be output by any deterministic oracle-based algorithm based on $T$ queries would necessarily have the objective value smaller by a factor at least $1 + \Omega(1/K) = 1 + \Omega(\frac{\log(T)}{T}).$ For this last step,  
we first argue in \Cref{lemma:lb-consistency} that all queries returned by the oracle are consistent with the final matrix $\mA.$ Since at query $t,$ the oracle returns the answers for the ``hindsight'' matrix $\mA_t,$ this means that the oracle successfully hides all information about $\bu_{\tau + 1}, \bv_{\tau}$, for $\tau \geq t.$ We then finally argue, in \Cref{lem: output}, that this hiding of information is sufficient to exhibit the claimed multiplicative gap.

\begin{lemma} \label{lem: A entry}
All entries of $\mA$ are in the interval $[1/K,8+2/K]$. In particular, $\mA>0$ entrywise.
\end{lemma}

\begin{proof}

We start by examining the last summand in  \eqref{eq: B explicit}.
Fix a coordinate $j=(k,m)$, and note that only the vectors
$\bu_{(k-1)\ell+1},\ldots,\bu_{k\ell}$ can be nonzero at this
coordinate. Let $\sigma=1$ if $m\le M/2$, and $\sigma=-1$
otherwise. The update rule for $U_k^\sigma$ implies that if
$(\bu_t)_j\ne\sigma$ for some $(k-1)\ell+1\le t\le k\ell$,
then $(\bu_s)_j=0$ for every $t<s\le k\ell$.
Thus, among these $\ell$ vectors, the coordinate equals
$\sigma$ until its first different value, and is zero in
every subsequent vector.

Denoting $\bu_0=\bu_{L+1}=\mathbf{0}$, note that the sequence
$(\bu_0)_j,\ldots,(\bu_{L+1})_j$ can change at most three
times: from zero to its first nonzero value, from $\sigma$
to a different nonzero value, and from its last nonzero
value to zero. The first and last changes have magnitude
at most $1$, and the middle change has magnitude at most $2$.
Therefore, using $\|\bv_t\|_{\infty}\le1$, every entry of the last
summand in \eqref{eq: B explicit} satisfies
\[
\bigg|
\sum_{t\in[L-1]}(\bv_t)_i
\big((\bu_{t+1})_j-(\bu_t)_j\big)
\bigg|
\le
\sum_{t=0}^{L}
|(\bu_{t+1})_j-(\bu_t)_j|
\le 1+2+1=4~,
\]
and so in particular, every entry of the last
summand is at least $-4$.
Further note that the last summand is supported only on the blocks $(k,k)$,
$k\in[K]$, and $(k,k+1)$, $k\in[K-1]$, while every matching-coordinate entry of
the second summand in \eqref{eq: B explicit} equals $4$. Hence, the sum of the second and the last
summands is entrywise nonnegative. The third summand is also
entrywise nonnegative, since $\bb_1+\bu_1\ge0$. Adding the first 
summand therefore results in any entry of $\mA$ having value of at least $1/K>0$.

The same argument also leads to the claimed entrywise upper bound. Indeed, the first three summands in \eqref{eq: B explicit} are entrywise bounded above by $1/K$, $4$ and $1/K$, respectively, while the calculation
above shows that the value of every entry of the last
summand is at most $4$. Summed up, this is at most $8+2/K$. 
\end{proof}

\begin{lemma} \label{lem: val}
$\val(\mA)\geq\frac{9}{K}+\frac{1}{K^2}$.    
\end{lemma}
\begin{proof}
For $k\in[K]$, let $m_k\in U_k^+$ be a coordinate at which all $\bu_t$ in block $k$ equal $1$, which exists by Lemma~\ref{lem: ut,vt}$(v)$.
Let $\by':=\frac{1}{K}\sum_{k\in[K]}\be_{k,m_k}$.
Note that $\by'\in\Delta_n$, and that for any $k\in[K]$ and $t\in[L]:~\bb_k^\top\by'=\bu_t^\top\by'=\frac{1}{K}$. Therefore, we get
\begin{align*}
\Big[\frac{1}{K}\bone\bone^\top\Big]&\by'=\frac{1}{K}\bone~,
\\
\Big[4\sum_{k\in [K]}\bb_k(\bb_k+\bb_{k+1})^\top\Big]&\by'=\sum_{k\in[K]}\bb_{k}\cdot\frac{8}{K}
=\frac{8}{K}\bone~,
\\
\Big[\frac{1}{2K}\bone(\bb_1+\bu_1)^\top\Big]&\by'=\frac{1}{2K}\bone\cdot \frac{2}{K}=\frac{1}{K^2}\bone~,
\\
\Big[\sum_{t\in[L-1]}\bv_t(\bu_{t+1}-\bu_t)^\top\Big]&\by'=\bzero.
\end{align*}
Hence, by summation
\(
\mA\by'=\big(\frac{9}{K}+\frac{1}{K^2}\big)\bone
\), 
and it follows that 
\(
\val(\mA)=\max_{\bx\in\Delta_n}\min_{i}(\mA\bx)_i\geq \min_{i}(\mA\by')_i=\frac{9}{K}+\frac{1}{K^2}.
\)
\end{proof}

\begin{lemma}\label{lemma:lb-consistency}
    The oracle responses are all consistent (in hindsight) with the fixed matrix $\mA$.
\end{lemma}
\begin{proof}
For any $t\leq T$, we see by Lemma~\ref{lem: ut,vt}$(iii)$ that
\begin{align*}
\bx_t^\top \mA-\bx_t^\top \mA_t&=\sum_{t\leq k\leq L-1}\bx_t^\top\bv_k(\bu_{k+1}-\bu_k)^\top
=\bzero,
\\ \mA\by_t-\mA_t\by_t
&=\sum_{t\leq k\leq L-1}\bv_k(\bu_{k+1}-\bu_k)^\top\by_t
=\bzero,
\end{align*}
and so $\bx_t^\top \mA_t=\bx_t^\top \mA$ and $\mA_t\by_t=\mA\by_t$.
\end{proof}

\begin{lemma} \label{lem: output}
If $\by\in\Delta_n$ satisfies $\bu_L^\top\by=0$, then $\min_{i}(\mA\by)_i\leq\frac{9}{K}+\frac{1}{2K^2}$.
\end{lemma}
\begin{proof}
For any $\bx\in\Delta_n$, it holds that
$\min_i(\mA\by)_i\leq\bx^\top\mA\by$, and so it suffices
to construct $\bx\in\Delta_n$ for which the claimed bound holds.
To that end, for $k\in[K]$, denote
\[
\lambda_k:=
\begin{cases}
    \frac{15K+1}{16K^2}, & k\text{ is odd},\\
    \frac{17K+1}{16K^2}, & k\text{ is even},
\end{cases}
\qquad \lambda_0:=\lambda_K~.
\]
For each block $k\in[K]$, choose one coordinate
$m_k^+\in V_k^+$ and one coordinate $m_k^-\in V_k^-$
from the final sets, which are nonempty by
Lemma~\ref{lem: ut,vt}$(v)$. Define
\[
\bx:=
\sum_{k\in[K]}
\left[
    \frac{\lambda_k+\frac{1}{2K}}{2}\,\be_{k,m_k^+}
    +
    \frac{\lambda_k-\frac{1}{2K}}{2}\,\be_{k,m_k^-}
\right].
\]
Thus $\bx$ is zero at every other coordinate.
Since $\lambda_k\geq\frac{1}{2K}$, we have $\bx\geq0$.
Moreover, $K$ is odd, so there are $(K+1)/2$ odd indices
and $(K-1)/2$ even indices, hence
\[
\bone^\top\bx
=
\sum_{k\in[K]}\lambda_k
=
\frac{K+1}{2}\cdot\frac{15K+1}{16K^2}
+
\frac{K-1}{2}\cdot\frac{17K+1}{16K^2}
=1~,
\]
which proves that $\bx\in\Delta_n$.

For every $k\in[K]$, we have
\[
\bx^\top\bb_k
=
\frac{\lambda_k+\frac{1}{2K}}{2}
+
\frac{\lambda_k-\frac{1}{2K}}{2}
=\lambda_k~.
\]
For every $t\in[L]$, let $k=\lceil t/\ell\rceil$.
By Lemma~\ref{lem: ut,vt} (items $(i)$ and $(v)$), $\bv_t$ is supported on
block $k$ and equals $1$ at $(k,m_k^+)$ and $-1$ at
$(k,m_k^-)$. Therefore,
\[
\bx^\top\bv_t
=
\frac{\lambda_k+\frac{1}{2K}}{2}
-
\frac{\lambda_k-\frac{1}{2K}}{2}
=\frac{1}{2K}~.
\]
Consequently,
\begin{align*}
\bx^\top\Big[\frac{1}{K}\bone\bone^\top\Big]
&=\frac{1}{K}\bone^\top~,
\\
\bx^\top\Big[
    4\sum_{k\in[K]}\bb_k(\bb_k+\bb_{k+1})^\top
\Big]
&=4\sum_{k\in[K]}\lambda_k(\bb_k+\bb_{k+1})^\top~,
\\
\bx^\top\Big[
    \frac{1}{2K}\bone(\bb_1+\bu_1)^\top
\Big]
&=\frac{1}{2K}(\bb_1+\bu_1)^\top~,
\\
\bx^\top\Big[
    \sum_{t\in[L-1]}\bv_t(\bu_{t+1}-\bu_t)^\top
\Big]
&=\frac{1}{2K}\sum_{t\in[L-1]}(\bu_{t+1}-\bu_t)^\top
=\frac{1}{2K}(\bu_L-\bu_1)^\top~.
\end{align*}
Since $1$ and $K$ are both odd, whereas $k-1$ and $k$
have opposite parity for $2\leq k\leq K$, summing these
equalities gives
\begin{align*}
\bx^\top\mA
&=
\frac{1}{K}\bone^\top
+4\sum_{k\in[K]}\lambda_k(\bb_k+\bb_{k+1})^\top
+\frac{1}{2K}\bb_1^\top
+\frac{1}{2K}\bu_L^\top
\\
&=
\frac{1}{K}\bone^\top
+\underbrace{4(\lambda_1+\lambda_K)}_{
    =\frac{15}{2K}+\frac{1}{2K^2}}\bb_1^\top
+\sum_{2\leq k\leq K}
\underbrace{4(\lambda_k+\lambda_{k-1})}_{
    =\frac{8}{K}+\frac{1}{2K^2}}\bb_k^\top
+\frac{1}{2K}\bb_1^\top
+\frac{1}{2K}\bu_L^\top
\\
&=
\left(\frac{1}{K}+\frac{8}{K}+\frac{1}{2K^2}\right)
\bone^\top
+\frac{1}{2K}\bu_L^\top~.
\end{align*}
Hence, for $\by\in\Delta_n$ satisfying $\bu_L^\top\by=0$,
we obtain
\(
\min_i(\mA\by)_i
\leq\bx^\top\mA\by
=\frac{9}{K}+\frac{1}{2K^2}.
\)
\end{proof}

To complete the proof of \Cref{thm: main lower}, combining Lemmas~\ref{lem: val} and \ref{lem: output}, we get that the algorithm's output $\by_{T+1}\perp \bu_L$ has multiplicative error of at least
\[
\frac{\val(\mA)}{\min_{i}(\mA\by_{T+1})_i}\geq
\frac{9/K+1/K^2}{9/K+1/(2K^2)}
=1+\frac{1}{18K+1}
~.
\]
By the choice of $K$ it holds that $\ell=\Theta(\log(T))$, $K=\Theta(T/\log(T))$, and $L=O(T)$. Also $4^\ell\leq4K$, so $n=8K(L+1)4^\ell=O(T^3)$. Thus the multiplicative error is $1+\Omega(\log(T)/T)$, which gives the claimed $\Omega(\log(1/\eps)/\eps)$ lower bound.

\begin{remark}[Packing LP lower bound]
\label{rem:packing-lower-bound}
Having proved Theorem~\ref{thm: main lower}, we now argue that the same proof also yields the claimed lower bound for
packing LPs. Consider $\max_{{\by \geq \bzero_n}:{\mA \by \le \bone_n}} \bone_n^\top \by,$
namely the packing LP with square matrix $\mA^\top$, and note that transposition only exchanges the sides of the two-sided matrix-vector products.

Consider the same hard instance as in \Cref{sec:instance-fam}, where we also include the algorithm's output $\widehat{\by}$ as
$\by_{T+1}$ when completing the construction, so that
$\bu_L^\top\widehat{\by}=0$; this only imposes one additional orthogonality
constraint (and does not require
$\widehat{\by}\in\Delta_n$).
The vector $\bx\in\Delta_n$ constructed in the proof of
\Cref{lem: output} then gives, for any feasible output,
\[
    \left(\frac9K+\frac{1}{2K^2}\right)
    \bone_n^\top\widehat{\by}
    = \bx^\top\mA\widehat{\by}
    \le 1.
\]
On the other hand,
by the same calculation as in
\Cref{lem: val},
choosing one coordinate from each final set
$U_k^-$ and assigning these coordinates weight $1/K$ gives
$\by^-\in\Delta_n$ with
$\mA\by^-=(9/K)\bone_n$.
Thus $(K/9)\by^-$ is feasible and
$\max_{{\by \geq \bzero_n}:{\mA \by \le \bone_n}} \bone_n^\top \by\ge K/9$. Consequently,
\[
\bone_n^\top\widehat{\by}
\le \frac{1}{9/K+1/(2K^2)}
\le \frac{\max_{{\by \geq \bzero_n}:{\mA \by \le \bone_n}} \bone_n^\top \by}{1+1/(18K)}.
\]
Since $K=\Theta(T/\log T)$, this yields the same
$\Omega(\log(1/\epsilon)/\epsilon)$ oracle lower bound.
\end{remark}

\section{Conclusion}\label{sec:conclusion} 

We gave a deterministic width-independent algorithm for positive matrix games \eqref{eq:packing-covering-game}, which capture packing and covering linear programs. Our algorithm simultaneously achieves $\frac{\mathrm{polylog(m, n, 1/\epsilon)}}{\epsilon}$ matrix-vector queries and parallel depth, with $\widetilde{O}(\nnz(\mA)/\epsilon)$ total work. Our approach is based on constructing a scale-invariant variational inequality and using coupled entropic geometry to construct a distance generating function for which the constructed variational inequality can be efficiently solved using mirror-prox-type iterations. We complemented this result with an oracle lower bound, which shows that every deterministic algorithm which accesses the matrix $\mA$ through two-sided matrix-vector products requires $\Omega(\epsilon^{-1}\log(1/\epsilon))$ queries in the worst case to output a $(1+\epsilon)$-multiplicative solution to \eqref{eq:packing-covering-game}. In particular, this lower bound shows that our algorithm is near-optimal (optimal up to a logarithmic factor) in the considered oracle model. Remaining questions for future work include closing the logarithmic gap in the oracle complexity and extending the work-depth guarantee to general mixed packing-covering feasibility.

\section*{Acknowledgments}

This work was supported in part by an AFOSR YIP award, under the contract number FA9550-24-1-0076. Any opinions, findings and conclusions or recommendations expressed in this material are those of the author(s) and do not necessarily reflect the views of the U.S.\ Department of Defense.

\bibliographystyle{alpha}
\bibliography{cleaned-up-bib}

\appendix

\section{Deferred Proofs}\label{sec:appx}

In this section, we provide complete proofs for auxiliary lemmas and claims stated in the main body. 

\subsection{Proof of Lemma~\ref{lemma:relaxed-rel-Lip-of-aux-reg}}\label{appx:proof-lemma:relaxed-rel-Lip-of-aux-reg}

Fix $\bu=(\bx^u,\by^u)$,
$\bz=(\bx^z,\by^z)$, and
$\bw=(\bx^w,\by^w)$ in $\mathcal{Z}$,
and assume that $\bx^z,\bx^w>\bzero$.
Denote
\begin{align*}
D:=&
D_{\varphi}(\bw,\bz)
+D_{\varphi}(\bu,\bw)
+D_{\varphi}(\bu,\bz),\\
C:=&
\left(
    \bx^w\oslash(\mA\by^w)
    -
    \bx^z\oslash(\mA\by^z)
\right)^\top
\mA(\by^u-\by^w).
\end{align*}
By definitions \eqref{eq:mon-op-def} and
\eqref{eq:aux-dgf}, the following identity holds:
\begin{align*}
(\mF(\bw)-\mF(\bz))^\top
(\bw-\bu)=
2(\nabla h(\bx^w)-\nabla h(\bx^z))^\top
(\bx^w-\bx^u)-
(\nabla\varphi(\bw)-\nabla\varphi(\bz))^\top
(\bw-\bu)+2C.
\end{align*}
Applying the three-point Bregman identity to $h$ and $\varphi$
then gives
\begin{align*}
(\mF(\bw)-\mF(\bz))^\top
  (\bw-\bu)&=
2\bigl[
    D_h(\bx^w,\bx^z)
    +D_h(\bx^u,\bx^w)
    -D_h(\bx^u,\bx^z)
\bigr]
\\&~~~~-D_{\varphi}(\bw,\bz)
-D_{\varphi}(\bu,\bw)
+D_{\varphi}(\bu,\bz)
+2C\\
&\le D+2C,    
\end{align*}
where the inequality follows from 
noting by direct differentiation that
\[
D_{\varphi}\bigl((\bx',\by'),(\bx,\by)\bigr)
=D_h(\bx',\bx)
+
D_h\left(\bx',\bx\odot\bigl((\mA\by')\oslash(\mA\by)\bigr)
\right)
\geq
D_h(\bx',\bx),
\]
for all $\bx,\bx',\by,\by'$.

It remains to bound $C$. It holds that
\begin{align*}
C&=
\left\|\sqrt{\bx^w\odot(\mA\by^u)\oslash(\mA\by^w)}
-\sqrt{\bx^z\odot(\mA\by^w)\oslash(\mA\by^z)}\right\|_2^2
-\underset{\leq 0}{\underbrace{\left\|\sqrt{\bx^w}-\sqrt{\bx^z\odot(\mA\by^u)\oslash(\mA\by^z)}\right\|_2^2}}
\\
&\leq\left\|\sqrt{\bx^w\odot(\mA\by^u)
\oslash(\mA\by^w)}- \sqrt{\bx^u}+\sqrt{\bx^u}-\sqrt{\bx^w}+\sqrt{\bx^w}
-\sqrt{\bx^z\odot(\mA\by^w)\oslash(\mA\by^z)}
\right\|_2^2
\\&\leq
3\left\|\sqrt{\bx^w\odot(\mA\by^u)
\oslash(\mA\by^w)}- \sqrt{\bx^u}
\right\|_2^2
+3\left\|\sqrt{\bx^u}-\sqrt{\bx^w}
\right\|_2^2
+3\left\|\sqrt{\bx^w}
-\sqrt{\bx^z\odot(\mA\by^w)\oslash(\mA\by^z)}
\right\|_2^2
\\&\leq
3D_h\left(
    \mathbf{x}^u,
    \mathbf{x}^w\odot(\mA\mathbf{y}^u)
    \oslash(\mA\mathbf{y}^w)
\right)+
3D_h(\mathbf{x}^u,\mathbf{x}^w)+
3D_h\left(
    \mathbf{x}^w,
    \mathbf{x}^z\odot(\mA\mathbf{y}^w)
    \oslash(\mA\mathbf{y}^z)
\right)
\\&=3\bigl[
    D_{\varphi}(\mathbf{u},\mathbf{w})
    +D_{\varphi}(\mathbf{w},\mathbf{z})
    -D_h(\mathbf{x}^w,\mathbf{x}^z)
\bigr]
\\&\le 3D,
\end{align*}
where the third inequality uses
the fact that $(\sqrt{\bx}-\sqrt{\by})^2\leq D_h(\bx,\by)$ (see, e.g., \citep[Lemma~2.4]{tsybakov2009introduction}), and the final inequality follows from the nonnegativity
of the entropy and Bregman divergences.
Consequently,
\[
(\mF(\bw)-\mF(\bz))^\top
(\bw-\bu)
\le D+2C\le 7D
=
7\left(D_{\varphi}(\bw,\bz)
+D_{\varphi}(\bu,\bw)
+D_{\varphi}(\bu,\bz)\right),
\]
which proves the claim for all vectors with strictly positive $\bx$-components.
Finally, the boundary case follows from a standard continuity argument by replacing each
$\bx$-component with
$(1-\delta)\bx+\delta\mathbf{1}_m/m$ and letting
$\delta\downarrow 0$.

\subsection{Proof of Claim~\ref{claim:xy-vs-pi-relationship}}\label{appx:proof-claim:xy-vs-pi-relationship}

Fix $t$ and $j\in[n]$, and let
$\bp_t:=\bpi_t^\top\bone_m$.
Throughout the proof, we omit from all sums over $i$ the rows with $x_{t,i}=0$,
since these satisfy $\Pi_{t,ik}=0$ for every $k$.
On the remaining rows, note that $\Pi_{t,ij}>0$ whenever $A_{ij}>0$,
as established in the proof of Lemma~\ref{lemma:dgf-properties}.
Thus, if $p_{t,j}=0$, the left-hand side of the claim is zero.
We may therefore assume $p_{t,j}>0$.
Note that by \eqref{eq:grad-psi-wrt-x}, for each $i$ such that $A_{ij}>0$ it holds that
\[
\log\Big(\frac{\Pi_{t, ij}(p_{t,j} + \epstr)^\beta}{A_{ij}y_{t,j}}\Big) =
(\nabla_\bx\psi(\bx_t, \by_t))_i - 1 - \beta,
\]
so in particular the right hand side is independent of $j$.
Consequently,
\begin{align*}
\sum_i x_{t,i}
\frac{A_{ij}y_{t,j}}{(\mA\by_t)_i}
&=
(p_{t,j}+\epstr)^\beta
\sum_i
\frac{\Pi_{t,ij}}{
    \sum_{k=1}^n
    \frac{\Pi_{t,ik}}{x_{t,i}}
    (p_{t,k}+\epstr)^\beta
}\\
&\le
(p_{t,j}+\epstr)^\beta
\sum_i \Pi_{t,ij}
\left(
    \sum_{k=1}^n
    \frac{\Pi_{t,ik}}{
        x_{t,i}(p_{t,k}+\epstr)
    }
\right)^\beta\\
&\le
(p_{t,j}+\epstr)^\beta
p_{t,j}^{1-\beta}
\left(
    \sum_i \frac{\Pi_{t,ij}}{x_{t,i}}
    \sum_{k=1}^n
    \frac{\Pi_{t,ik}}{p_{t,k}+\epstr}
\right)^\beta,
\end{align*}
where we applied Jensen's inequality twice: first with the convex function $s\mapsto s^{-\beta}$,
weights $\Pi_{t,ik}/x_{t,i}$, and arguments
$(p_{t,k}+\epsilon_{\mathrm{tr}})^{-1}$, followed by taking reciprocals;
then with the concave function $s\mapsto s^\beta$,
weights $\Pi_{t,ij}/p_{t,j}$, and arguments
$\sum_{k=1}^n \frac{\Pi_{t,ik}}{x_{t,i}(p_{t,k}+\epsilon_{\mathrm{tr}})}$.
Further note that $\Pi_{t,ij}\le x_{t,i}$ gives $\sum_i \frac{\Pi_{t,ij}}{x_{t,i}}
\sum_{k=1}^n
\frac{\Pi_{t,ik}}{p_{t,k}+\epstr}
\le
\sum_{k=1}^n
\frac{p_{t,k}}{p_{t,k}+\epstr}
\le n$,
and so
\[
(p_{t,j}+\epstr)^\beta
p_{t,j}^{1-\beta}
\left(
    \sum_i \frac{\Pi_{t,ij}}{x_{t,i}}
    \sum_{k=1}^n
    \frac{\Pi_{t,ik}}{p_{t,k}+\epstr}
\right)^\beta
\le
(p_{t,j}+\epstr)^\beta
p_{t,j}^{1-\beta}n^\beta
\leq
e(p_{t,j}+\epstr),
\]
where the last inequality used $0<\beta\le1$, $p_{t,j}\le p_{t,j}+\epstr$
and $n^\beta
    =\exp\left(\frac{\log(n)}{1+\log(n)}\right)
    \le e$,
thus completing the proof.

\subsection{Proof of Claim~\ref{claim:psi_1''-ub-y}}\label{appx:proof-claim:psi_1''-ub-y}

Denote once again $\bp_t:=\bpi_t^\top\bone_m$,
so that $\dot{\bp}_t=\dot{\bpi}_t^\top\bone_m$.
For $t\in(0,1)$, any row with $x_{t,i}=0$ is identically
zero along the segment and may be omitted from the sums below.
For each $j\in[n]$, Cauchy--Schwarz gives
\begin{align*}
\left(
    \dot p_{t,j}
    -(p_{t,j}+\epstr)
      \frac{\dot y_{t,j}}{y_{t,j}}
\right)^2
&=
\bigg(
    \sum_{i:A_{ij}\ne0}
    \left(
        \dot\Pi_{t,ij}
        -\Pi_{t,ij}\frac{\dot y_{t,j}}{y_{t,j}}
    \right)
    -\epstr
     \frac{\dot y_{t,j}}{y_{t,j}}
\bigg)^2\\
&\le
(p_{t,j}+\epstr)
\bigg[
    \sum_{i:A_{ij}\ne0}
    \frac{
        \big(
            \dot\Pi_{t,ij}
            -\Pi_{t,ij}\dot y_{t,j}/y_{t,j}
        \big)^2
    }{\Pi_{t,ij}}
    +
    \epstr
    \frac{\dot y_{t,j}^2}{y_{t,j}^2}
\bigg]. 
\end{align*}

Since $y_{t,j}\ge\epstr$, we have
$\dot y_{t,j}^2/y_{t,j}
\ge\epstr\dot y_{t,j}^2/y_{t,j}^2$.
Thus, \eqref{eq:psi_1''} and the preceding inequality imply
\begin{align*}
\psi_1''(t)
&\ge
\sum_{j=1}^n
\frac{
    \left(
        \dot p_{t,j}
        -(p_{t,j}+\epstr)
          \dot y_{t,j}/y_{t,j}
    \right)^2
    +\beta\dot p_{t,j}^2
}{
    p_{t,j}+\epstr
}\\
&=
\sum_{j=1}^n
\frac{1+\beta}{p_{t,j}+\epstr}
\left(
    \dot p_{t,j}
    -\frac{p_{t,j}+\epstr}{1+\beta}
     \frac{\dot y_{t,j}}{y_{t,j}}
\right)^2+
\frac{\beta}{1+\beta}
\sum_{j=1}^n
(p_{t,j}+\epstr)
\frac{\dot y_{t,j}^2}{y_{t,j}^2}\\
&\ge
\frac{\beta}{1+\beta}
\sum_{j=1}^n
(p_{t,j}+\epstr)
\frac{\dot y_{t,j}^2}{y_{t,j}^2},
\end{align*}
where the equality follows by completing the square.
This is the desired inequality.

\subsection{Proof of \Cref{claim:splitting-solve}}\label{appx:alt-min}

To prove the claim, for convenience, we first restate the optimization problem \eqref{eq:separable-prob} that the claimed algorithm solves efficiently to target additive error:
\begin{equation}\label{eq:separable-prob-appx}
    \min_{\bq\in\Delta_n,\ \by\in\Deltatr_n}f(\bq, \by), \; \text{ where }\; f(\bq, \by):= \sum_{j=1}^n q_j\log\frac{q_j}{c_jy_j}
        +\beta\, h(\bq+\epstr\bone_n)+h(\by)+\bd^\top\by.
\end{equation}
The algorithm we use can be seen as entropic mirror descent applied to $\widehat{f}(\bq):= \min_{\by \in \Deltatr_n} f(\bq, \by),$ where minimization over $\by$ is done approximately, so we need to account for the approximation error. To analyze this entropic mirror descent update, we use that $\widehat{f}$ is both relatively smooth and relatively strongly convex with respect to $h$ (with condition number $\Theta(1/\beta)$), so, with suitably small error in the minimization over $\by,$ it converges geometrically fast. See, e.g., \cite[Theorem 3.1]{lu2018relatively} for convergence under exact access to $\widehat{f}$ and its gradient. We note here that while there exist results addressing similar settings \cite{stonyakin2021inexact}, their model of ``inexact oracle'' for $\widehat{f}$ is different than what we obtain here, thus we provide the full argument, for completeness. 

The main algorithm is summarized in \Cref{algo:alt-min}. It uses, as a subroutine, a known result for separable generalized allocation problems from \cite{hochbaum1994lower}, summarized in the following fact. We note here that it is possible to use alternative solvers for approximating $\widehat{f}.$ For instance, it is possible to argue quadratic convergence (in loglog steps) for a constrained version of Newton's method. Nevertheless, we opt for the near-linear time solver below because it is simpler to state and use as a black box. 

\begin{fact}[\cite{hochbaum1994lower}]\label{fact:y-minimization}
    Let $M \geq 1$, $n \geq 2$ be integers and let functions $F_j : \{0, \dots, M\} \to \R$, $j \in [n],$ have nondecreasing increments $F_j(\ell) - F_j(\ell-1),$ $\ell \in \{1,\dots, M\}.$ With constant-time access to function evaluations $F_j$, an exact minimizer of 
    \[
    \min\Big\{\sum_{j=1}^n F_j(k_j): \bk \in \Z_+^n, \sum_{j=1}^n k_j = M\Big\}
    \]
    can be computed in $O(n\log(M+n))$ work. There is also an implementation with $O(n \log^2(nM))$ work and poly-logarithmic depth in $nM$. Neither implementation enumerates the $nM$ increments.
\end{fact}
We briefly remark here that while \cite{hochbaum1994lower} does not explicitly discuss parallel implementations of her algorithm, parallel implementation is straightforward, using standard techniques. In particular, the parallel bound follows by viewing allocation in \cite{hochbaum1994lower} as selection from \(n\) implicit sorted arrays of marginal costs and using weighted-median pruning with parallel binary searches: there are \(O(\log(nM+2))\) rounds, each requiring \(O(n\log(nM+2))\) work and polylogarithmic depth.

Define $\by^{\bq} := \argmin{\by \in \Deltatr_n} f(\bq, \by)$ (unique by strict convexity), so that $\widehat{f}(\bq) = f(\bq, \by^{\bq}).$ To apply \Cref{fact:y-minimization}, define $F^{\bq}_j(y_j) = y_j \log(y_j) + d_j y_j - q_j \log(y_j),$ so that $\sum_{j=1}^n F^{\bq}_j(y_j)$ contains all $\by$-dependent terms from $f(\bq, \by),$ and, as a consequence, $\by^{\bq} =  \argmin{\by \in \Deltatr_n} \sum_{j=1}^n F^{\bq}_j(y_j).$ 

We briefly argue that the application of the results from  \Cref{fact:y-minimization} is sufficient to guarantee the required approximation error and runtime needed for the proof of \Cref{claim:splitting-solve}, summarized in the following claim.

\begin{algorithm}[ht]
    \caption{Y-minimization Solve for \Cref{algo:alt-min}}\label{algo:y-min}
    \DontPrintSemicolon
\small
\KwIn{$n, \bq, \bd, \varepsilon \in (0, 1/32], \epstr$.}
\KwOut{$\by \in \Deltatr_n$ such that $\|\log(\by\oslash \by^{\bq})\|_\infty \leq \varepsilon$ and $f(\bq, \by) - \widehat{f}(\bq) \leq \varepsilon^2.$}

Set $M = \big\lceil \frac{n}{\epstr^2\varepsilon} \big\rceil,$ $w = \frac{1-n\epstr}{M}$\;

Use the algorithm from \Cref{fact:y-minimization} to find an exact solution $\bk_*$ of $\min_{\bk \in \Z_+^n: \sum_{j=1}^n k_j = M} \sum_{j=1}^n F^{\bq}_j(\epstr + w k_j)$\;

\Return{$\by = \epstr \bone_n + w \bk_*$}
\end{algorithm}

\begin{algorithm}[ht]
\caption{Alternating Mirror Descent solver for \Cref{claim:splitting-solve}}
\label{algo:alt-min}
\DontPrintSemicolon
\small
\KwIn{$\epstr, n, \beta, \bc, \bd, \delta \in (0, \epstr)$.}
\KwOut{$\bq, \by$ such that $f(\bq, \by) - \min_{\bq' \in \Delta_n, \by' \in \Deltatr_n}f(\bq', \by') \leq \delta$.}

\tcp{Set parameters.}

$\bq_0 = (1/n)\bone_n$, $\varepsilon = \frac{\beta\delta}{16(1 + \beta)},$ $K = \max\{1,\, \lceil \frac{\log(2(1+ \beta)\log(n)/\delta)}{\log(1+\beta)} \rceil\}$

\For{$k=0,\ldots,K-1$}{
    \tcp{Approximate mirror descent update}

    Let $\by_k$ be the output of \Cref{algo:y-min} with input $n, \bq_k, \bd, \varepsilon, \epstr$\;

    $q_{k+1, j} \gets \frac{(c_{j} y_{k, j})^{1/(1+\beta)}(q_{k, j}/(q_{k, j} + \epstr))^{\beta/(1+\beta)}}{\sum_{j' = 1}^n (c_{j'} y_{k, j'})^{1/(1+\beta)}(q_{k, j'}/(q_{k, j'} + \epstr))^{\beta/(1+\beta)}}$, $j \in [n]$
    }

\Return{$\bq_K, \by_{K-1}$}\;

\end{algorithm}

\begin{claim}\label{claim:y-solve}
    Given $\bq \in \Delta_n$ and $\varepsilon \in (0, 1/32],$ \Cref{algo:y-min} returns $\by \in \Deltatr_n$ such that $\|\log(\by\oslash \by^{\bq})\|_\infty \leq \varepsilon$ and $f(\bq, \by) - \widehat{f}(\bq) \leq \varepsilon^2.$ The total work that the algorithm takes is $O(n\log(n/(\epstr \varepsilon)))$ and it is implementable in parallel with  $O(n\log^2(n/(\epstr \varepsilon)))$ work and polylogarithmic depth in $n, 1/\epstr, 1/\varepsilon.$ 
\end{claim}
\begin{proof}[Proof of \Cref{claim:y-solve}]
    Under the notation of \Cref{algo:y-min}, the objective $\sum_{j=1}^n F^{\bq}_j(\epstr + w k_j)$ is separable and convex, so \Cref{fact:y-minimization} applies, and we immediately get the claimed runtime for computing the discrete minimizer $\bk_*$. Thus, we only need to bound the error for the output $\by  = \epstr \bone_n + w \bk_*.$ Recall that $\by^{\bq} = \argmin{\bv \in \Deltatr_n}\sum_{j=1}^n F^{\bq}_j(v_j).$ For notational convenience, let $F(\bv) := \sum_{j=1}^n F^{\bq}_j(v_j).$   

    Let $\overline{\by}^{\bq}$ denote $\by^{\bq}$ with its coordinates rounded to the grid $\epstr + w k_j$, $k_j \in \{0, \dots, M\},$ (to one of the closest two points on the grid, if not already on the grid)  such that $\sum_j k_j = M.$
    This is possible because $\by^{\bq} \in \Deltatr_n$ and $wM  = 1 - n \epstr.$ It is immediate that $\bone_n^\top \overline{\by}^{\bq} = 1,$ $\|\overline{\by}^{\bq} - \by^{\bq}\|_\infty \leq w$, and $y^{\bq}_j = \epstr$ implies $\overline{y}^{\bq}_j = \epstr.$ 

    Now consider $\nabla F(\by^{\bq}).$ For all coordinates $j$ such that $y^{\bq}_j > \epstr$, due to the simplex constraint on $\by^{\bq}$, the value of $\nabla_j F(\by^{\bq})$ is the same. Since, by construction, $y^{\bq}_j = \epstr$ implies $\overline{y}^{\bq}_j = \epstr$ and both $\bone_n^\top\overline{\by}^{\bq} = 1$ and $\bone_n^\top\by^{\bq} = 1$, we can conclude that
    \begin{equation}\label{eq:rounded-sol-first-order-zero}
        \inner{\nabla F(\by^{\bq}), \overline{\by}^{\bq} - \by^{\bq}} = 0. 
    \end{equation}
    Observe that for any $\bv \in \Deltatr_n,$ $1 \leq (F_j^{\bq}(v_j))'' = \frac{1}{v_j} + \frac{q_j}{v_j^2} \leq \frac{2}{\epstr^2}.$ As a consequence, using \eqref{eq:rounded-sol-first-order-zero}, we have $F(\overline{\by}^{\bq}) - F(\by^{\bq}) \leq \frac{\norm{\overline{\by}^{\bq} - \by^{\bq}}_2^2}{\epstr^2} \leq \frac{n\, w^2}{\epstr^2}.$ Since the output $\by$ minimizes $F$ on the grid, we have $F(\by) \leq F(\overline{\by}^{\bq})$, and so, by the choice of parameter $w,$ 
    \begin{equation}\label{eq:min-y-gap-guarantee}
       f(\bq, \by) - \widehat{f}(\bq) = F(\by) - F(\by^{\bq}) \leq \frac{n\, w^2}{\epstr^2} \leq \varepsilon^2. 
    \end{equation}
    It remains to argue that $\|\log(\by\oslash \by^{\bq})\|_\infty \leq \varepsilon$. Because both $\by, \by^{\bq} \in \Deltatr_n$, the logarithm is entrywise $(1/\epstr)$-Lipschitz, thus $\|\log(\by\oslash \by^{\bq})\|_\infty \leq \frac{1}{\epstr} \|\by - \by^{\bq}\|_2$. Since $(F_j^{\bq}(v_j))'' = \frac{1}{v_j} + \frac{q_j}{v_j^2} \geq 1,$ $F$ is $1$-strongly convex with respect to the $\ell_2$-norm, thus $\frac{1}{2} \|\by - \by^{\bq}\|_2^2 \leq F(\by) - F(\by^{\bq})$ as $\by^{\bq}$ minimizes $F.$ Combining the two inequalities with \eqref{eq:min-y-gap-guarantee}, we get
    \begin{equation}\notag
        \|\log(\by\oslash \by^{\bq})\|_\infty \leq \frac{\sqrt{2n}w}{\epstr^2} \leq \varepsilon,
    \end{equation}
    completing the proof. 
\end{proof}
To complete the proof of \Cref{claim:splitting-solve}, it remains to argue that: (1) $\widehat{f}$ is $\beta$-relatively strongly convex and $(1+\beta)$-relatively smooth with respect to $h$ and (2) the update in Line 4 of \Cref{algo:alt-min} is an entropic mirror descent update for $\widehat{f}.$ The rest of the analysis is then a minor extension of \cite[Theorem 3.1]{lu2018relatively}. 

For (1), the claim that $\widehat{f}$ is $\beta$-relatively strongly convex with respect to $h$ already follows from \Cref{claim:f-rel-str-cvx}, using \Cref{fact:inf-convexity} to minimize $f$ over $\by.$ Note that by \Cref{fact:danskin} we also have $\nabla \widehat{f}(\bq) = \nabla_{\bq}f(\bq, \by^{\bq}).$ For $(1+\beta)$-relative smoothness, observe that for fixed $\by,$ the Bregman divergence of $f$ at $\bq', \bq$ is $D_h(\bq', \bq) + \beta D_{h}(\bq' + \epstr \bone_n, \bq + \epstr \bone_n) \leq (1+\beta)D_h(\bq', \bq).$ Thus $\widehat{f}$ is $(1+\beta)$-relatively smooth with respect to $h$.

For (2), by a direct calculation, Line 4 in \Cref{algo:alt-min} is equivalent to the entropic mirror descent update given by:
\begin{equation}\label{eq:q-MD-update}
    \bq_{k+1} = \argmin{\bu \in \Delta_n}\Big\{\inner{\nabla_{\bq}f(\bq_k, \by_k), \bu - \bq_k} + (1+\beta)D_h(\bu, \bq_k)\Big\}, 
\end{equation}
which is solved by $q_{k+1, j} \propto q_{k, j} \exp\big(-\frac{[\nabla_{\bq}f(\bq_k, \by_k)]_j}{1+\beta}\big)$ normalized to sum up to one. 

Using the three-point inequality (see \Cref{fact:3PI}) for \eqref{eq:q-MD-update} and $(1+\beta)$-relative smoothness of $f$ for fixed $\by = \by_k,$ we have
\begin{equation}\label{eq:MD-3PI}
    \begin{aligned}
      &\;   f(\bq_{k+1}, \by_k) - f(\bq_*, \bv_*) + (1+\beta)D_h(\bq_*, \bq_{k+1})\\
      \leq \; &f(\bq_k, \by_k) - f(\bq_*, \bv_*) + \inner{\nabla_q f(\bq_k, \by_k), \bq_* - \bq_k} + (1+\beta) D_h(\bq_*, \bq_k). 
    \end{aligned}
\end{equation}
Now, by \Cref{claim:y-solve}, $f(\bq_k, \by_k) - \widehat{f}(\bq_k) \leq \varepsilon^2$ and $\|{\nabla_q f(\bq_k, \by_k) - \nabla \widehat{f}(\bq_k)}\|_\infty = \norm{\log(\by_k \oslash \by^{\bq_k})}_{\infty} \leq \varepsilon$. We thus have:
\begin{align*}
    f(\bq_k, \by_k) - f(\bq_*, \bv_*) + \inner{\nabla_q f(\bq_k, \by_k), \bq_* - \bq_k} &\leq \widehat{f}(\bq_k) - \widehat{f}(\bq_*) + \inner{\nabla  \widehat{f}(\bq_k), \bq_* - \bq_k} + \varepsilon^2 + 2\varepsilon\\
    &= -D_{\widehat{f}}(\bq_*, \bq_k) + \varepsilon^2 + 2\varepsilon\\
    &\leq -\beta D_h(\bq_*, \bq_k) + \varepsilon^2 + 2\varepsilon,
\end{align*}
where we have used H\"{o}lder's inequality to bound $\inner{{\nabla_q f(\bq_k, \by_k) - \nabla \widehat{f}(\bq_k)}, \bq_* - \bq_k}$ along with $\norm{\bq_* - \bq_k}_1 \leq 2$ in the first inequality, and $\beta$-relative strong convexity of $\widehat{f}$ with respect to $h$ in the second inequality. Plugging the last displayed inequality back into \eqref{eq:MD-3PI} and crudely bounding $\varepsilon^2 + 2\varepsilon \leq 3\varepsilon,$ we now have
\begin{equation}\label{eq:MD-recursion}
    f(\bq_{k+1}, \by_k) - f(\bq_*, \bv_*) + (1+\beta)D_h(\bq_*, \bq_{k+1}) \leq D_h(\bq_*, \bq_k) + 3\varepsilon. 
\end{equation}
Unrolling the recursion in \eqref{eq:MD-recursion}, since $f(\bq_{k+1}, \by_k) - f(\bq_*, \bv_*) \geq 0$ and $D_h(\bq_*, \bq_0) \leq \log(n)$ (by the chosen initialization), we finally obtain for any $k \geq 0$ (and, in particular, for $k = K-1$):
\begin{equation}\notag
    f(\bq_{k+1}, \by_k) - f(\bq_*, \bv_*) \leq \frac{\log(n)}{(1+\beta)^k} + \frac{3\varepsilon}{1 - 1/(1+\beta)}.
\end{equation}
Plugging in $k = K-1$ and the values of $\beta$ and $\varepsilon,$ the above bound is at most $\delta,$ as claimed.

Finally, there are $K = O(\log(n)\log(n/\delta))$ iterations in \Cref{algo:alt-min}. Each iteration calls \Cref{algo:y-min} once and updates $\bq_k$ using $O(n)$ additional arithmetic operations. The latter computation can be implemented in parallel in $\log(n)$ depth (to compute the normalizing factor), with $O(n \log(n))$ work. This completes the proof of \Cref{claim:splitting-solve}.

\end{document}